\input amstex

\define\scrO{\Cal O}
\define\Pee{{\Bbb P}}
\define\Zee{{\Bbb Z}}
\define\Cee{{\Bbb C}}
\define\Ar{{\Bbb R}}
\define\Q{{\Bbb Q}}
\define\Pic{\operatorname{Pic}}
\define\Ker{\operatorname{Ker}}

\define\Sym{\operatorname{Sym}}
\define\Hom{\operatorname{Hom}}
\define\Aut{\operatorname{Aut}}

\define\Id{\operatorname{Id}}

\define\ad{\operatorname{ad}}
\define\Ad{\operatorname{Ad}}
\define\dbar{\bar{\partial}}

\define\proof{\demo{Proof}}
\define\endproof{\qed\enddemo}
\define\endstatement{\endproclaim}
\define\theorem#1{\proclaim{Theorem #1}}
\define\lemma#1{\proclaim{Lemma #1}}
\define\proposition#1{\proclaim{Proposition #1}}
\define\corollary#1{\proclaim{Corollary #1}}
\define\claim#1{\proclaim{Claim #1}}

\define\section#1{\specialhead #1 \endspecialhead}
\define\ssection#1{\medskip\noindent{\it  #1}\medskip\nobreak}


\documentstyle{amsppt}

\topmatter
\title Holomorphic principal bundles over elliptic curves
\endtitle
\author {Robert Friedman and  John W. Morgan}
\endauthor
\address Department of Mathematics, Columbia University, New York, NY
10027, USA\endaddress
\email rf\@math.columbia.edu, jm\@math.columbia.edu  \endemail
\thanks The first author was partially supported by NSF grant
DMS-96-22681. The second author was partially supported by NSF grant
DMS-94-02988. 
\endthanks

\endtopmatter

\nologo

\document

\section{Introduction.}

This paper, the first in a projected series of three, is concerned with
the classification of holomorphic principal $G$-bundles over an
elliptic curve, where $G$ is a reductive linear algebraic group. The
motivation for this study comes from physics. The  F-theory/heterotic
string duality predicts that, given an elliptically fibered Calabi-Yau
manifold $M$ of dimension $n$ over a base space $B$, together with a
stable $E_8\times E_8$-bundle over $M$ and a so-called complexified
K\"ahler class, there should be an associated Calabi-Yau manifold
$W$ of dimension
$n+1$, fibered over the same base  $B$, where the  fibers are elliptic
$K3$ surfaces. In some sense, this paper and its two sequels are an
attempt to understand this prediction in purely mathematical terms. The
physical approach to this construction has been discussed in \cite{14},
as well as in the many references in that paper. For several
mathematical reasons it has seemed worthwhile to consider not just
$E_8\times E_8$ but the case of a general reductive group $G$. In this
paper, we shall consider the problem of classifying holomorphic
$G$-bundles over a single smooth elliptic curve $E$, as well as proving
a number of auxiliary results which we shall need later.  The methods
of this paper  are not well-suited to dealing with singular curves of
arithmetic genus one or with families. Moreover, they do not work well
in defining universal bundles, even locally in the case of a single
smooth elliptic curve. In fact, all of these problems are already
evident in the case of vector bundles with trivial determinant, i\.e\.
in case $G = SL_n(\Cee)$. In the first sequel to this paper, we will
describe another method for constructing bundles, at least when
$G$ is simple, by considering deformations of certain minimally
unstable $G$-bundles corresponding to special maximal parabolic
subgroups of $G$. It turns out that this construction overcomes the
problems outlined above, and it can be used to give a new proof and
generalization of a theorem of Looijenga and Bernshtein-Shvartsman
\cite{24, 4} concerning the global structure of the moduli space of
$G$-bundles. In Part III of this series, we specialize to the case
motivated by physics and relate the case where the group
$G$ is $E_6, E_7, E_8$ to del Pezzo surfaces and simple elliptic
singularities. The study  of flat $G$-bundles over $E$ also leads to a
somewhat different set of questions, involving normal forms for two or
more commuting elements in a compact Lie group. This question is studied
from the point of view of Lie groups in \cite{6}.

The study of holomorphic principal bundles over compact Riemann
surfaces has a long history, dating back in modern times to
Grothendieck's 1956 paper \cite{18} on the structure of principal
$G$-bundles over $\Pee^1$. As is very well-known, Grothendieck showed
there that every holomorphic vector bundle over $\Pee^1$ is isomorphic
to a direct sum of line bundles, thereby showing that there is a
holomorphic reduction of the structure group to the diagonal subgroup
of $GL_n(\Cee)$. What is less well-known is that, in the same paper, he
went on to prove the analogous statement for principal $G$-bundles, for
an arbitrary reductive algebraic group $G$, by showing that in this
case the structure group reduces holomorphically to a Cartan subgroup
of $G$. A year later, Atiyah classified all holomorphic vector bundles
over an elliptic curve
\cite{1}, without however considering more general structure groups. In
this case, the structure group does not in general reduce to a diagonal
subgroup of $GL_n(\Cee)$. Since the notion of a semistable bundle was
not available at that time, Atiyah was not able to discuss the moduli
space in a meaningful way and could not attempt to construct universal
bundles.

The study of vector bundles over a Riemann surface was greatly
clarified by Mumford's introduction in 1961 of the notion of stability
or semistability of a vector bundle. Shortly thereafter, Narasimhan and
Seshadri \cite{25} showed that a vector bundle of degree zero is stable
if and only if it arises as the flat vector bundle associated to an
irreducible unitary representation of $\pi_1(C,*)$. A slight
modification, using a  central extension of $\pi_1(C,*)$, is needed to
handle the case of arbitrary degree. The  definition of stability and
the analogue of the Narasimhan-Seshadri theorem for arbitrary reductive
groups $G$ were worked out by Ramanathan \cite{27,28}. Finally,
Atiyah-Bott \cite{2} showed that there was a correspondence between
flat vector bundles, or more generally flat $K$-bundles where $K$ is a
compact semisimple group,  and Yang-Mills theory. There is also a
slightly technical generalization to arbitrary compact groups involving
representations of the central extension of $\pi_1(C,*)$ discussed
above. One basic reason to work with semistable bundles is that the set
of all semistable bundles can be parametrized by a coarse moduli space,
whereas this is impossible for the set of all bundles. However, in
order to get a separated moduli space, one must identify two semistable
bundles if they are S-equivalent, which is a weaker equivalence
relation than isomorphism in general. It follows easily from the
definitions and Atiyah's classification theorem that the coarse moduli
space of vector bundles with a fixed determinant over an elliptic curve
is in fact simply a projective space.

Atiyah's classification, in the case of rank two and trivial
determinant, was  analyzed in a relative setting in \cite{11, 12, 13}.
The major motivation for this analysis came from Donaldson theory,
where it was necessary to understand the moduli space of stable rank
two vector bundles over an elliptic surface.  In order to do this, one
constructs relative moduli spaces of vector bundles over families of
elliptic curves and tries to write down relative universal vector
bundles over these relative moduli spaces. The point is that a vector
bundle whose restriction to every fiber is semistable defines a section
of the corresponding relative moduli space. Conversely, given a section
of the relative moduli space, one can pull back a universal bundle via
the section to obtain a bundle over the elliptic surface, and then see
how all bundles can be obtained from this basic construction. These
ideas were re-examined for arbitrary rank over elliptic fibrations over
an arbitrary base in
\cite{16}.

Given Grothendieck's theorem and Atiyah's classification of vector
bundles, it is natural to try to understand the structure of the moduli
space of $S$-equivalence classes of semistable
$G$-bundles over an elliptic curve for more general groups
$G$. However, this  question does not seem to have arisen
mathematically in the intervening 40 years since the appearance of
Grothendieck's paper, with the exception of the theorem of Looijenga 
and Bernshtein-Shvartsman, which as stated does not even appear at
first glance to involve
$G$-bundles. We will discuss this theorem in more detail in Section 5.
By contrast, the case of bundles over curves of higher genus has
received a great deal of attention, mainly due to the efforts to prove
the Verlinde formula (see for example \cite{3}, \cite{10}). 

We turn now to a more detailed description of the contents of this
paper. Let
$G$ be a connected reductive (complex linear algebraic) group,
with maximal compact subgroup $K$. The  main result of this paper
(Theorem 4.5) generalizes Atiyah's theorem  by giving a classification
up to isomorphism of holomorphic principal $G$-bundles $\xi$ over an
elliptic curve $E$. The classification is in terms of three pieces of
data.  The first is the Levi factor $L$ of a parabolic subgroup $P$ of
$G$, the Harder-Narasimhan parabolic of the bundle $\xi$. The bundle
$\xi$ reduces to a semistable bundle $\xi_L$ over $L$ and $L$ is maximal
among Levi factors of parabolic subgroups with this property. The
second datum is a representation $\rho$ of a central   extension of
$\pi_1(E,*)$ by $\Ar$ into  the maximal compact subgroup  $K_L$ of $L$
This representation is    required to send the copy of $\Ar$  into the
center of $K_L$. These representations are considered up to conjugation
in $K_L$. Such a representation is equivalent to a Yang-Mills
connection on the principal $K_L$-bundle whose curvature is central in
$K_L$. The third datum is the conjugacy class of a nilpotent element in
the Lie algebra of the centralizer in $L$ of the image of $\rho$. We
show how this classification generalizes Atiyah's theorem for vector
bundles. The remainder of the paper is concerned with the local and
global structure of the moduli space of S-equivalence classes of
semistable $G$-bundles. We also describe special representatives for
the S-equivalence classes and their automorphism groups. 

The paper is organized as follows. We begin by reviewing Atiyah's
classification of vector bundles in Section 1.  Following the statement
of his theorem, we recall the definitions of semistability,
S-equivalence, and the Harder-Narasimhan filtration for vector bundles.
Using these definitions and some of the basic results, we outline
another proof of the classification theorem. Then we discuss the
relation between semistable bundles and representations of the
fundamental group of an elliptic curve (or a central extension of it),
or equivalently Yang-Mills  connections. (In the case that the image
of the representation lies in $SL_n(\Cee)$, the representation factors
through the fundamental group of
$E$ and the resulting Yang-Mills connection is flat.) Bundles with
Yang-Mills connections have maximal dimensional automorphism group in
their S-equivalence class. For many applications to moduli questions,
it is important to classify semistable vector bundles which are  at
the opposite extreme from vector bundles with Yang-Mills connections,
in the sense that their automorphism groups have minimal possible
dimension in their S-equivalence class.  Such bundles are called {\sl
regular\/}, and already play an important role in studying vector
bundles over elliptic fibrations \cite{11, 12, 13, 16}. The plan then
for the rest of the paper will be to use Section 1 as a guide for
generalizations to other reductive groups
$G$. Unlike the case of vector bundles, however, we cannot give
completely elementary proofs of the main results, but must rely on the
machinery of Yang-Mills and deformation theory.

In Section 2, we reduce the case of unstable $G$-bundles to the
semistable case via the analogue of the Harder-Narasimhan filtration.
For a general reductive group $G$ and a general compact Riemann surface
$C$, this construction says that every unstable
$G$-bundle over $C$ has a canonical reduction to a bundle with
structure group a parabolic subgroup of $G$, in such a way that the
associated bundle over the Levi factor is semistable. In case $g(C)=0$
or $1$, one can further reduce the structure group to the Levi factor.
For the case of genus zero, we use this fact to give a quick proof of
Grothendieck's theorem on $G$-bundles over $\Pee^1$. Over a general
curve, there is an invariant which  Atiyah and Bott have associated to
an unstable $G$-bundle. We use this invariant, in the case of curves of
genus $1$, to determine when a reduction of structure from $G$ to a
reductive subgroup yields the Harder-Narasimhan reduction described
above.

The results of Section 2 enable us to concentrate attention on
semistable $G$-bundles. By general theory, the S-equivalence classes of
semistable 
$G$-bundles form a coarse moduli space. Thus, the classification of
semistable $G$-bundles divides into two parts: First describe all
bundles in a given S-equivalence class, and then describe the set of
all S-equivalence classes. Sections 3 and 4 are concerned with the
first problem. The first step is to choose a canonical representative
for each S-equivalence class. It turns out that there are two possible
such representatives, the Yang-Mills representative and the regular
representative; for an open dense set of $G$-bundles, these
representatives will coincide. In order to describe the possible
S-equivalence classes, we will single out the Yang-Mills
representatives as the canonical choice. They are constructed via
Yang-Mills connections, and  correspond to conjugacy classes of
homomorphisms $\rho$ from a central extension of $\pi_1(E,*)$ into the
maximal compact  subgroup $K$ of $G$. Using this description, we
identify the set of semistable $G$-bundles S-equivalent to a Yang-Mills
bundle arising from the representation $\rho$ with the set of conjugacy
classes of nilpotent elements $X$ in the subalgebra $\frak z_{\frak
g}(\rho)$ of the Lie algebra of $G$ which centralizes $\rho$. To make
this identification, we give normal forms for the transition functions
of a semistable $G$-bundle in Section 3, depending only on the pair
$(\rho, X)$, where $\rho$ is the representation corresponding to the
Yang-Mills representative and $X$ is a nilpotent element in $\frak
z_{\frak g}(\rho)$. In Section 4, we show that two pairs $(\rho, X)$
and $(\rho', X')$ lead to isomorphic bundles exactly when there is an
element of $G$ conjugating $\rho$ to $\rho'$ and $X$ to $X'$.  The
generalized form of Atiyah's theorem is then a recapitulation of the
results in Sections 2, 3, and 4. We also describe  the structure of the 
automorphism group of a semistable bundle.

In Section 5, we turn to the local and global structure of the moduli
space. We  relate the singularities of the moduli space, viewed as an
orbifold, to the components of the centralizer of the corresponding
representation $\rho$ from $\pi_1(E,*)$ into $K$. Quoting results from
\cite{6}, we are able to identify these singularities explicitly. Next
we turn to the general global description of the moduli space. If $G$
is semisimple, it is described as the quotient of a product of copies
of the elliptic curve $E$ by a finite group. In case
$G$ is simple and simply connected, the theorem of Looijenga and
Bernshtein-Shvartsman says that the moduli space is a weighted
projective space. Analogous results hold in the general case. 

In Section 6, we define the notion of a regular bundle. These are
bundles whose automorphism groups have minimal possible dimension. We
show that each S-equivalence class contains a unique regular
representative. The automorphism groups of regular bundles are studied
in detail, and we make the link between the deformation theory of a
regular bundle and the local structure of the moduli space. Finally,
Section 7 deals with the example of the classical groups and gives
elementary constructions for all regular semistable bundles in each
case.

Throughout this paper we adopt the following notation: $G$ is a
connected, reductive complex linear algebraic group with $K$ as a
maximal compact subgroup. The Lie algebra of $G$ is $\frak g$, and
$\frak h$  denotes a Cartan subalgebra of $\frak g$. Also, $C$ denotes
a smooth projective complex curve of arbitrary genus whereas $E$ is
always an elliptic curve. We always assume that we have chosen a fixed
origin
$p_0$ for the elliptic curve $E$.

\medskip
\noindent {\bf Acknowledgements.} Our interest in
this circle of problems grew out of our joint work with Ed Witten
as described in \cite{14}. Indeed, many of the
results in this paper were established jointly with him during the
course of that work. It is a pleasure to
thank him for many stimulating conversations on these subjects, and
especially for his patience in explaining to us
the physics behind many of his ideas and questions.
 We would also
 like to thank A. Borel, P. Deligne, W. Schmid, and M.
 Thaddeus for helpful
discussions and correspondence on various aspects of this paper.

\section{1.  Holomorphic vector bundles over an elliptic curve.} 

\ssection{1.1. Statement of Atiyah's theorem.}

Let us begin by giving  Atiyah's theorem as originally stated in
\cite{1}. Recall that a vector bundle $V$ on a projective variety is
{\sl indecomposable\/} if $V$ is not the direct sum of two proper
subbundles. An easy argument using the finite-dimensionality of
cohomology shows that every vector bundle on a projective variety can
be written as a direct sum of indecomposable bundles, where the
summands and their multiplicities are uniquely determined up to
isomorphism. Thus, to classify all holomorphic vector bundles over
$E$, it suffices to describe the indecomposable ones. To do so, we
begin with the following:

\lemma{1.1}  For any $n\ge 1$, there is a unique indecomposable vector
bundle $I_n$  of rank $n$ over
$E$ all of whose Jordan-H\"older constituents are isomorphic to
$\scrO_E$. The bundle  $I_n$ has a canonical filtration
$$\{0\} \subset F^1\subset \cdots \subset F^n = I_n,$$ with $F^i\cong
I_i$ and $F^{i+1}/F^i\cong \scrO_E$. 
\qed
\endstatement

We shall describe explicit transition functions for $I_n$, in a
somewhat more general context, in Section 3 of this paper.

For a line bundle $\lambda$ of degree zero on $E$, we let
$I_n(\lambda) = I_n\otimes \lambda$. As we shall see, the
indecomposable vector bundles of degree zero over $E$ are exactly the
bundles $I_n(\lambda)$. To describe bundles of nonzero degree, we have
the following, whose proof is also in \cite{1}:

\lemma{1.2} Let $n$ be a positive integer. For each integer $a$
relatively prime to
$n$ and each line bundle
$\lambda$ over $E$ of degree $a$, there exists up to isomorphism a
unique indecomposable bundle
$W_n(a; \lambda)$ over $E$ of rank $n$ and such that
$\det W_n(a; \lambda) =\lambda$. Moreover, $ W_n(a; \lambda)$ is simple.
\qed
\endstatement

Note that, if $a=0$, then $ W_n(a; \lambda)$ must have rank one and so
is isomorphic to $\lambda$. In general, for $\lambda =\scrO_E(ap_0)$,
we let $W_n(a;
\scrO_E(ap_0)) = W_n(a)$ and set $W_n(1) = W_n$. Define
$I_{d}(W_{n}(a;\lambda))$ to be $I_{d}\otimes W_{n}(a;\lambda)$; it has
a filtration, all of whose successive quotients are isomorphic to
$W_{n}(a;\lambda)$, and in fact it is the unique indecomposable bundle
of rank $dn$ with this property.

We can now state Atiyah's theorem.

\theorem{1.3} The bundle  $I_{d}(W_{n}(a;\lambda))$, where  $a$ is
relatively prime to $n$, is indecomposable, and every indecomposable
holomorphic bundle is isomorphic to $I_{d}(W_{n}(a;\lambda))$ for a
suitable choice of $d,n, \lambda$. Every holomorphic vector bundle
$V$ over $E$ is a direct sum of vector bundles of the form
$I_{d_i}(W_{n_i}(a_i; \lambda_i))$, where the $n_i$ and $d_i$ are
positive integers and the
$\lambda_i$ are line bundles of degree $a_i$, where $a_i$ is relatively
prime to $n_i$. Moreover, the triples $(n_i, d_i,\lambda_i)$ are
uniquely specified up to permutation by the isomorphism type of
$V$.   
\qed 
\endstatement

We shall outline a proof of this theorem below.

\ssection{1.2. Semistable holomorphic vector bundles.}

Let us restate the above results from a more modern perspective, which
can then be generalized to the case of groups other than $GL_n$. We
begin   by recalling the standard definitions:

\definition{Definition 1.4}  Let $C$ be a smooth projective curve and
$V$ be a nonzero holomorphic vector bundle over
$C$. We define the {\sl slope\/} $\mu(V)$ by
$$\mu(V) = \deg V/\operatorname{rank}V.$$ The bundle $V$ is {\sl
semistable\/} if, for every holomorphic subbundle $W$ of $V$   with $0
< \operatorname{rank} W < \operatorname{rank} V$, we have
$\mu(W) \leq \mu (V)$. Stability of $V$ is defined similarly, by
replacing the above inequality by a strict inequality.
\enddefinition

For example, $\mu(I_{d}(W_{n}(a;\lambda))) = a/n$. As we shall see
below, the bundles $W_{n}(a;\lambda)$ are stable, and in fact they are
exactly the stable bundles over $E$. It then follows easily by
induction that
$I_{d}(W_{n}(a;\lambda))$ is semistable for all $d\geq 1$.

\theorem{1.5} Let $V$ be a holomorphic vector bundle over a smooth
projective curve 
$C$. Then there is a unique filtration 
$$\{0\}\subset F_1\subset \cdots \subset F_k = V$$ by subbundles $F_i$
of $V$, with the property that the sequence
$\mu(F_i/F_{i-1})$ is  strictly decreasing, and such that $F_i/F_{i-1}$
is semistable.
\qed
\endstatement

The filtration above is called the {\sl Harder-Narasimhan filtration\/}
of $V$. For example, $V$ is semistable if and only if the filtration is
trivial, i\.e\. $F_1=V$ above. 

In the case of an elliptic curve, it follows  by an easy direct
argument that the Harder-Narasimhan filtration is split. As a
consequence, we see the first part of the next theorem, which is a
variant of Theorem 1.3: 

\theorem{1.6} Let $E$ be an elliptic curve. 
\roster
\item"{(i)}" Every holomorphic vector bundle
$V$ over $E$ is isomorphic to a direct sum $\bigoplus_iV_i$ of
semistable bundles, where $\mu(V_i)>\mu(V_{i+1})$. 
\item"{(ii)}" Let $V$ be a semistable bundle over $E$ with slope
$\mu(V) = a/n$, where $n$ is a positive integer and $a$ is an integer
relatively prime to $n$. Then
$V$ is a direct sum of bundles of the form $I_d(W_n(a; \lambda))$,
where $\lambda$ is a line bundle of degree $a$. \qed
\endroster
\endstatement

The $V_i$ are determined up to isomorphism by the isomorphism type of
$V$, but the splitting is not canonical.

Given (i) of Theorem 1.6, the proof of Theorem 1.3 reduces to the
classification of semistable bundles given in (ii). An outline of this
classification is as follows: first, one checks that for each $n>0$,
each integer $a$ relatively prime to $n$, and each line bundle
$\lambda$  of degree $a$, there exists a stable vector bundle
$W_n(a;\lambda)$ as in Lemma 1.2. Next, suppose that
$V$ is semistable with slope $\mu(V) = a/n$ for some positive integer
$n$ and some integer
$a$ relatively prime to $n$. An argument with the
Grothendieck-Riemann-Roch theorem shows that there exists a line bundle
$\lambda$ of degree zero such that
$\Hom(W_n(a;\lambda), V) \neq 0$. By stability, every nonzero map from
$W_n(a;\lambda)$ to $V$ must be an injection onto a subbundle. In
particular, $V$ has a filtration whose successive quotients are of the
form $W_n(a;\lambda_i)$ for appropriate choices of $i$, and examining
the possibilities leads to the final classification.
\medskip

We return for a moment to a general curve $C$. Associated to a
holomorphic vector bundle $V$ of rank $n$ over $C$ and its
Harder-Narasimhan filtration as in Theorem 1.5 is the decreasing
sequence of integers
$\{\mu(F_i/F_{i-1})\}$.  Following Atiyah-Bott \cite{2}, we replace this
strictly decreasing sequence of rational numbers by the following: to
each $i$, if 
$n_i$ is the rank of $F_i/F_{i-1}$, then we associate to $F_i/F_{i-1}$
the
$n_i$-vector
$\nu_i = (\mu(F_i/F_{i-1}),
\dots, \mu(F_i/F_{i-1}))$. To $V$ itself we associate the $n$-vector
$$\mu(V)=\mu = (\mu_1, \dots, \mu_n) = (\nu_1, \dots, \nu_k).$$ Note
that the sequence $\mu_i$ is weakly decreasing.

\corollary{1.7} Let $V$ be a vector bundle of rank $n$ over $C$ and let
$\mu(V)$ be the vector defined above. Then\rom:
\roster
\item"{(i)}"  $\sum _i\mu _i=\deg V$, and hence  is an integer\rom;
\item"{(ii)}" If $d(\mu_i)$ is the multiplicity of $\mu_i$ in $\mu$,
i\.e\. the number of $j$ such that
$\mu_j=\mu_i$, then $\mu_i\cdot d(\mu_i) \in \Bbb Z$.
\endroster Conversely, given a weakly decreasing sequence
$\mu_1,\ldots,\mu_n$ of rational numbers satisfying the above two
conditions,  if we set $\mu=(\mu_1,\ldots,\mu_n)$, then there is a
vector bundle $V$ for which $\mu(V)=\mu$.
\qed
\endstatement

In general, we can view the vector $\mu = (\mu_1, \dots, \mu_n)$ as
lying in the standard Cartan subalgebra $\frak h$ of the Lie algebra
$\frak{gl}_n$ of 
$GL_n(\Cee)$, where $\frak h$ is the set of diagonal matrices. The fact
that the $\mu_i$ are decreasing means that $\mu$ lies in the standard
Weyl chamber defined by the usual   ordering of the roots for
$\frak{sl}_n$. Moreover, $V$ is semistable if and only if $\mu$ lies in
the center of $\frak{gl}_n$.

\ssection{1.3. S-equivalence and Yang-Mills connections on vector
bundles.}

Let us begin this section by recalling a standard definition.

\definition{Definition 1.8} A {\sl family\/} of vector  bundles over
the curve $C$ parametrized by a complex space (or scheme) $S$ is a
holomorphic vector bundle
$\Cal V$ over $C\times S$. The family $\Cal V$ is a family of {\sl
semistable\/}  vector bundles over $C$ if $\Cal V|C\times
\{s\}=\Cal V _s$ is semistable for all $s\in S$.  Finally, let
$V$ and $V'$ be two semistable bundles over
$C$. We say that $V$ and
$V'$ are {\sl S-equivalent\/} if there exists a family of semistable
bundles
$\Cal V$ parametrized by  an irreducible $S$ and a point $s\in S$ such
that, for
$t\neq s$, $\Cal V_t \cong V$ and $\Cal V_s\cong V'$. More generally,
we let S-equivalence be the equivalence relation generated by the above
relation.
\enddefinition

One can describe S-equivalence quite explicitly. Every semistable
vector bundle $V$ has a filtration by subbundles $V_i$ such that the
quotients $V_i/V_{i-1}$ are stable bundles with $\mu(V_i/V_{i-1})
=\mu(V)$. Such a filtration, which is not in general canonical, is
called a {\sl Jordan-H\"older filtration}. The associated graded bundle
$\operatorname{gr}V$ is then a direct sum of stable bundles of the same
slope, and its isomorphism class is independent of the choice of
Jordan-H\"older filtration. In general, a direct sum of stable bundles
of the same slope is called {\sl polystable}. We then have the
following characterization of S-equivalence:

\theorem{1.9} Two semistable bundles $V_1$ and $V_2$ are S-equivalent
if and only  if the associated graded bundles $\operatorname{gr}V_1$ and
$\operatorname{gr}V_2$ are isomorphic.
\qed
\endstatement

In general, S-equivalence can be coarser than isomorphism. On the other
hand, we must factor out by S-equivalence if we are to have a separated
moduli space, and indeed it is a standard result that there is a coarse
moduli space for the set of semistable vector bundles modulo
S-equivalence. One goal in studying vector bundles is to single out a
best representative in each S-equivalence class. For simplicity, we
shall begin by describing one such choice in the case where $V$ is
semistable and has trivial determinant, in other words the case of an
$SL_n(\Cee)$-bundle. As we have seen, in this case $V$ is isomorphic to
$\bigoplus _iI_{d_i}(\lambda_i)$, where the $\lambda_i$ are line
bundles of degree zero. By induction on $d$, $I_d$ is S-equivalent to
$\scrO_E^d$, and hence 
$$\bigoplus _iI_{d_i}(\lambda_i) \, \text{ is S-equivalent to }
\,\bigoplus _i(\underbrace{\lambda_i\oplus \cdots \oplus
\lambda_i}_{\text{$d_i$ times}}).$$ In particular, every semistable
holomorphic vector bundle
$V$ with trivial determinant is S-equivalent to a direct sum of line
bundles of degree zero. It is easy to see that the line bundles and
their multiplicities are uniquely determined by $V$.  Conversely, $V$
is determined by the line bundles $\lambda_i$ and the integers $d_i$,
which can be thought of as Jordan blocks corresponding to the summand
$I_{d_i}(\lambda_i)$.

The above picture can be described in terms of flat
$SU(n)$-connections as well. By the theorem of Narasimhan-Seshadri, a
semistable vector bundle $V$ with trivial determinant is S-equivalent
to a vector bundle given by a flat $SU(n)$-connection, or equivalently
by a representation $\rho\: \pi_1(E) \to SU(n)$. Since $\pi_1(E)$ is
abelian and every commuting pair of matrices in $SU(n)$ can be
simultaneously diagonalized, a vector bundle given by a flat
$SU(n)$-connection is isomorphic to a direct sum of $n$ line bundles
$\lambda_1, \dots, \lambda_n$ on $E$ such that $\lambda_1\otimes
\cdots \otimes \lambda_n \cong \scrO_E$. It is easy to see that the set
of the $\lambda_i$, together with the multiplicities, is unique up to
order.

If $V$ is a semistable vector bundle whose determinant is of non-zero
degree we cannot hope to have an S-equivalent bundle with a flat $U(n)$
connection, because the degree of the determinant does not change under
S-equivalence and no bundle with non-zero degree can have a flat
connection. Thus, the notion of flat connection must be generalized in
the case of $U(n)$. This was done by Atiyah-Bott \cite{2}. They
describe semistable vector bundles whose determinant is not necessarily
trivial via central Yang-Mills connections on $U(n)$-bundles. Let $C$
be a smooth projective curve, together with a Hermitian metric. Fix a
$C^\infty$ vector bundle $V$ of rank $n$ over $C$, or equivalently a
$C^\infty$ $U(n)$-bundle. A {\sl Yang-Mills connection\/} $A$ on $V$ is
a $U(n)$-connection $A$ whose curvature $F_A$ is covariantly constant.
If moreover $*F_A$ is a constant multiple of the identity matrix in
$Hom (V,V)$, we call $A$ a {\sl central Yang-Mills connection}. A
Yang-Mills connection $A$ defines a holomorphic structure on $V$, by
taking the $(0,1)$-part of $A$.

A  version of the Narasimhan-Seshadri theorem then says that there is a
canonical representative, up to isomorphism, in each  S-equivalence
class, namely the central curvature Yang-Mills representative.  

\theorem{1.10} Let $C$ be a smooth projective curve.
\roster 
\item"{(i)}" If $V$ is a holomorphic vector bundle whose holomorphic
structure is defined by a central Yang-Mills connection, then $V$ is
semistable, and in fact it is polystable, i\.e\. a direct sum of stable
bundles with the same slope.
\item"{(ii)}" Conversely, if $V$ is a holomorphic semistable vector
bundle on $C$, then $V$ is S-equivalent to a holomorphic vector bundle
$V_0$  defined by a central Yang-Mills connection $A$, unique up to
$U(n)$-gauge equivalence.
\item"{(iii)}" Two holomorphic semistable vector  bundles $V_1$ and
$V_2$ are S-equivalent if and only if the corresponding central
Yang-Mills connections are $U(n)$-gauge equivalent. \qed
\endroster
\endstatement

There is also a formulation of this result in terms of representations.
Let $\alpha$ and $\beta$ be a pair  of oriented generators for
$\pi_1(E)$. There is a unique central extension
$$0 \to \Ar \to \Gamma _\Ar \to \pi_1(E) \to 0,$$ where if $A$ and $B$
are any two lifts of $\alpha, \beta$ to $\Gamma_\Ar$, then the
commutator $[A, B] = 1\in \Bbb R$. Of course, there is a similar
construction for a curve of arbitrary genus. In particular, taking the
quotient of $\Gamma_\Ar$ by $\Bbb Z\subset \Ar$, we see that
$\Gamma_\Ar/\Zee \cong U(1) \times \pi_1(E)$. In any case,
$\Gamma_\Ar$ is a Lie group. A {\sl central representation\/} of
$\Gamma_\Ar$ to a Lie group $G$ is a (continuous) representation
$\rho\: \Gamma_\Ar \to G$ such that $\rho(\Ar)$ is contained in the
center of $G$.  With this said, we have the following theorem of
Atiyah-Bott:

\theorem{1.11} Fix once and for all a Yang-Mills connection $A_0$ on
the line bundle $L$ over $C$ with $c_1(L) = 1$. Then the choice of
$A_0$ induces a bijection from the set of central Yang-Mills
connections on a $C^\infty$ vector bundle of rank $n$, modulo gauge
equivalence, and the set of isomorphism classes of central
representations $\rho\: \Gamma _\Ar \to U(n)$.  \qed
\endstatement

\corollary{1.12} Given the fixed Yang-Mills connection $A_0$ on the line
bundle $L$ over $C$ with $c_1(L) = 1$, the  following three sets of
data are equivalent\rom:
\roster
\item"{(i)}" A holomorphic semistable vector bundle over $C$ of rank
$n$, up to S-equi\-va\-lence\rom; 
\item"{(ii)}" A central Yang-Mills connection on a $C^\infty$ vector
bundle of rank $n$ over $C$, up to gauge equivalence\rom;
\item"{(iii)}" A  central representation 
$\rho\: \Gamma _\Ar \to U(n)$, up to conjugation.\qed
\endroster
\endstatement

We note that, in the above correspondence between isomorphism classes
of central representations and semistable bundles, irreducible
representations correspond to stable bundles.  

To make the link between this discussion and the previous
classification, note that to classify central representations 
$\rho\: \Gamma _\Ar \to U(n)$, up to conjugation, we may assume that
$\rho$ is irreducible. Given such a representation, let $A\in U(n)$ be
the image of a lift of
$x\in \pi_1(E)$ to $\Gamma_\Ar$ and let $B$ be the image of a lift of
$y\in \pi_1(E)$. Then $[A,B] = ABA^{-1}B^{-1} = \lambda\Id$ for some
nonzero complex number $\lambda$. It follows that $\lambda^n = 1$, and
$\rho$ is irreducible if and only if $\lambda = \exp(2\pi\sqrt{-1}a/n)$
for an integer $a$ relatively prime to
$n$. Thus $\rho$ must be given on the central subgroup $\Ar$ by
$$\rho(t) = \exp(2\pi\sqrt{-1}at/n)\cdot \Id,$$ where $a$ is a
well-defined integer (not just modulo $n$). It then follows that, up to
conjugation, we can assume that $B$ is diagonal with entries $\mu,
\lambda
\mu, \cdots \lambda^{n-1}\mu$ and that $Ae_1 = \nu e_2$, $Ae_i=
e_{i+1}$ for
$i>1$, $Ae_n = e_1$, for uniquely specified nonzero complex numbers
$\mu, \nu$ of absolute value $1$. This leads to another description of
the bundles
$W_n(a; \lambda)$.

\ssection{1.4. Regular vector bundles.} 

In the preceding discussion, we found a best representative for the
S-equivalence class of a holomorphic vector bundle by looking at a flat
representative, or a bundle with a Yang-Mills connection with central
curvature  in the case where the determinant was not necessarily
trivial.  It turns out that such representatives have holomorphic
automorphism groups which are as large as possible, and which can be
quite complicated. For example, if
$V=\scrO_E^n$ is the trivial bundle, then $\Aut V = GL_n(\Cee)$. Here
we describe a choice whose automorphism group is as small as possible.

\definition{Definition 1.13} A semistable holomorphic vector bundle
$V$ is {\sl regular\/} if, for all vector bundles $V'$ which are
S-equivalent to $V$, $\dim \Aut V \leq \dim \Aut V'$.
\enddefinition

We shall see that, given an S-equivalence class of semistable bundles,
a regular representative is unique up to isomorphism.

We begin by determining the automorphism groups of the bundles we have
described in \S 1.1. A proof is given in
\cite{16}, as well as in Lemma 4.3 of this paper.

\lemma{1.14} If $I_d$ is the bundle of Lemma \rom{1.1}, then $\Hom
(I_d, I_d)
\cong \Cee[T]/(T^d)$, with $t\in \Cee$ acting as $t\cdot\Id$, and a
choice of $T$ is given by choosing a surjection $I_d \to I_{d-1}$,
followed by an inclusion of $I_{d-1}$ in
$I_d$. In particular  $\Hom (I_d, I_d)$ is an abelian $\Cee$-algebra of
dimension $d$.
\qed
\endstatement

A similar argument handles the case of the bundles $I_d(W_n(a;
\lambda))$. In fact, we have the following general result:

\lemma{1.15} Suppose that $W$ is a stable bundle on $E$ and that $I$
and $I'$ are bundles of the form $\bigoplus _iI_{d_i}$. Then the map 
$\varphi \mapsto \varphi\otimes \Id$ defines an isomorphism $$\Hom
(I,I')
\to \Hom (I\otimes W, I'\otimes W).$$   In particular, $\Hom
(I_d(W_n(a; \lambda)), I_d(W_n(a;
\lambda))$ is an abelian $\Cee$-algebra of dimension $d$, isomorphic to
$\Cee[T]/(T^d)$.
\endstatement

\proof Clearly the map from $\Hom (I,I')$ to $\Hom (I\otimes W,
I'\otimes W)$ is an inclusion. To see that it is onto, note that there
is an isomorphism of vector bundles
$$Hom (I\otimes W, I'\otimes W) \cong Hom (I,I') \otimes Hom (W,W).$$
Using $\frac1{n}\operatorname{trace}\: Hom (W,W) \to \scrO_E$, there is
an induced map 
$$Hom (I\otimes W, I'\otimes W) \cong Hom (I,I') \otimes Hom (W,W) \to
Hom(I,I'),$$ which leads to the corresponding map on global sections
$\tau\:
\Hom (I\otimes W, I'\otimes W)\to
\Hom (I,I')$. Note that $\tau$ is functorial with respect to morphisms
both of $I$ and
$I'$. The composition of
$\tau$ with the inclusion of $\Hom (I,I')$ in
$\Hom (I\otimes W, I'\otimes W)$ is the identity. Thus, it suffices to
prove that if $\varphi\: I\otimes W \to I'\otimes W$ is such that
$\tau(\varphi) = 0$,  then $\varphi = 0$. The proof is by double
induction on the ranks of $I$ and
$I'$. Note that if $I$ and $I'$ both have rank one, then $\Hom
(I\otimes W, I'\otimes W) \cong \Hom (W,W) \cong \Cee$ and the
statement is clear. Suppose that the statement has been shown for all
$I'$ of a given rank and bundles $I$ of rank less than $k$. Given $I$
of degree $k$, if
$I_0$ is a proper degree zero subbundle of
$I$, then
$\tau (\varphi|I_0\otimes W) = 0$ by functoriality. Thus by induction
$\varphi|I_0\otimes W = 0$ and $\varphi$ induces a homomorphism $\bar
\varphi$ from
$(I/I_0)\otimes W$ to $I'\otimes W$ with $\tau(\bar \varphi) = 0$, and
again by induction $\bar \varphi = 0$. Thus $\varphi =0$. A similar
induction handles the case where the rank of $I'$ varies but that of
$I$ is held fixed, by looking at degree zero quotients of $I'$. This
completes the inductive step and the proof of the lemma.
\endproof

\corollary{1.16} Suppose that $W$ and $W'$ are two stable bundles over
$E$ with $\mu(W) = \mu(W')$ and that $I$ and $I'$ are as in Lemma
\rom{1.15}. Then $\Hom (I\otimes W, I'\otimes W')\cong \Hom (I,I')$ if
$W\cong W'$, and $\Hom (I\otimes W, I'\otimes W') = 0$ otherwise.
\endstatement

\proof The first statement follows from Lemma 1.15. To see the second,
it suffices to show that $\Hom (W,W') =0$ if $W$ is not isomorphic to
$W'$, since both $I\otimes W$ and $I'\otimes W'$ have filtrations whose
associated gradeds are isomorphic to a direct sum of a number of copies
of $W$ (resp\. $W'$). But if $W$ and $W'$ are stable and $\mu(W) =
\mu(W')$, then every nonzero homomorphism from $W$ to $W'$ is
necessarily an isomorphism. Thus either $\Hom (W,W') =0$ or $W$ is
isomorphic to $W'$. This concludes the proof of Corollary 1.16.
\endproof

We can now give an explicit description of regular bundles.

\theorem{1.17} Let $V$ be a semistable vector bundle of rank $n$ and
degree $d$. Let 
$e = \gcd(d,n)$ and set $k=n/e, a= d/e$. Suppose that $V\cong
\bigoplus _iI_{d_i}(W_k(a;\lambda _i))$ as in Theorem \rom{1.3}, so
that $e=\sum _id_i$.  Then\rom:
\roster
\item"{(i)}" With the above notation,  $\dim \Aut V \geq e$, with
equality if and only if, for $i\neq j$, then $\lambda _i\neq
\lambda _j$.
\item"{(ii)}" If $\dim \Hom (V, V) = e$, then $\Hom (V, V)$ is abelian
and in fact in the above notation $\Hom (V,V) \cong \bigoplus
_i\Cee[t]/(t^{d_i})$.
\endroster
\endstatement

\proof If $V\cong \bigoplus _iI_{d_i}(W_k(a;\lambda _i))$, then
$$\Hom (V,V) \cong \bigoplus _{i,j}\Hom(I_{d_i}(W_k(a;\lambda _i)),
I_{d_j}(W_k(a;\lambda _j))).$$ But by Corollary 1.16,
$$\Hom(I_{d_i}(W_k(a;\lambda _i)), I_{d_j}(W_k(a;\lambda _j)))=\cases 0
, &\text{if $\lambda_j\not\cong\lambda _i$;}\\ 
\Hom(I_{d_i}, I_{d_j}) &\text{if $\lambda_j\cong\lambda _i$.}
\endcases$$

Thus, fixing  a line bundle
$\lambda_i$, the nonzero factors in
$Hom(V,V)$ are of  the form
$$\Hom \Big(\bigoplus _{\lambda_j\cong \lambda_i}I_{d_j}, \bigoplus
_{\lambda_j\cong \lambda_i}I_{d_j}\Big).$$  So it suffices to show that
the dimension of such a space is strictly bigger than the sum over all
$d_j$ such that $\lambda_j\cong
\lambda_i$ unless $\lambda_j\cong \lambda_i$ implies
$j=i$, i\.e\. unless there is just one summand in the direct sum. We
leave this straightforward argument to the reader.
\endproof 

\corollary{1.18}  A semistable vector bundle is regular if and only if
it is isomorphic to
$$\bigoplus_i \Big(I_{d_i}\otimes W_k(a;\lambda_i)\Big)$$ where, for
$i\neq j$ we have $\lambda_i\not\cong \lambda_j$. Thus, each
S-equivalence class of semistable vector bundles has a unique regular
representative, up to isomorphism.
\qed
\endstatement

Our goal in this paper is now to generalize the results of this section
to principal holomorphic $G$-bundles over $E$, where $G$ is now an
arbitrary complex reductive group.

\section{2. The Harder-Narasimhan reduction for holomorphic principal
bundles.}

\ssection{2.1. The space of all holomorphic structures on a principal
  bundle.}

We begin with some generalities concerning $C^\infty$ and holomorphic
principal
$G$-bundles. For the moment, let $\xi$ denote a fixed $C^\infty$
principal $G$-bundle over the smooth projective curve $C$.   The
topological types of such bundles are classified by
$H^2(C;
\pi _1(G)) \cong \pi  _1(G)$.

For a $C^\infty$ principal $G$-bundle $\xi$, a holomorphic structure on
$\xi$ is given by an integrable $\dbar$-operator on $\xi$, i\.e\. an
integrable
$(0,1)$-connection on $\xi$. Note that $C^\infty$
$(0,1)$-connections always exist, and that every $(0,1)$-connection is
integrable since $\dim C = 1$. The set of all $(0,1)$-connections is an
affine space $\Cal A$ over $\Omega^{0,1}(\ad \xi)$, the space of
$C^\infty$ $(0,1)$-forms on $C$ with values in the vector bundle $\ad
\xi$  associated to the adjoint representation of $G$. The complex
gauge group $\Cal G^\Cee$ acts on
$\Cal A$ and the quotient is the set of holomorphic structures on
$\xi$. In particular, for a fixed $C^\infty$-type of a principal
$G$-bundle, the space of all complex structures is connected and
nonempty. For the rest of this paper, we shall denote by $\xi$ a
holomorphic principal $G$-bundle, in other words a $C^\infty$
$G$-bundle together with some fixed choice of holomorphic structure.

\ssection{2.2. Semistable principal bundles.}

The quotient $\Cal A/\Cal G^\Cee$ of all holomorphic structures on a
given $C^\infty$-principal $G$-bundle is highly non-separated, and to
find a separated space we need to restrict to the set of semistable
principal $G$-bundles. For the purposes of this paper, we shall define
semistability as follows:

\definition{Definition 2.1}  Let
$G$ be a connected reductive complex Lie group and let $\xi \to C$ be a
holomorphic principal
$G$-bundle. Then
$\xi$ is {\sl semistable\/} if  the associated vector bundle $\ad
\xi$ is a semistable vector bundle. The principal $G$-bundle $\xi
\to C$ is {\sl unstable\/} if it is not semistable.
\enddefinition
 
The above definition differs from that given in Ramanathan
\cite{27}, but the two definitions are equivalent. In fact, the
following  is shown in essentially contained in \cite{27} (see for
example Proposition 7.1) and
\cite{2} (e.g. Proposition 10.6 and 10.9):

\theorem{2.2} Let $G$ be a connected reductive group and let $\xi \to
C$ be a principal $G$-bundle. The following conditions are
equivalent\rom:
\roster
\item"{(i)}" The bundle $\xi$ is semistable.
\item"{(ii)}"The principal $G/\Cal C(G)$-bundle associated to $\xi$ is
  semistable, where $\Cal C(G)$ is the identity component of the center
of $G$ and the quotient
  group $G/\Cal C(G)$ is semisimple, possibly trivial. 
\item"{(iii)}" For every  irreducible representation
$\phi\: G \to GL(V)$, where $V$ is a finite-dimensional vector space,
the vector bundle
$V(\xi) = \xi \times _GV$ associated to
$\xi$ via $\phi$ is semistable.
\item"{(iv)}" There exists a finite-dimensional vector space $V$ and an
irreducible  representation
$\phi\: G \to GL(V)$ whose kernel is finite modulo the center such that
the  associated vector bundle $V(\xi)$ is semistable. 
\item"{(v)}" For every parabolic subgroup $P$ of $G$, and every dominant
character $\chi\: P\to \Cee^*$, if the structure group of $\xi$ reduces
to $P$, then the associated line bundle over $C$ defined by the
character $\chi$ has nonpositive degree.  \qed
\endroster
\endstatement

Stability of a $G$-bundle is defined via (v) above, by requiring that,
for every  parabolic subgroup $P$ of $G$, and every dominant character
$\chi\: P\to \Cee^*$,  if the structure group of $\xi$ reduces to $P$,
then the associated line bundle over $C$ defined by the character
$\chi$ has strictly negative degree. Even if $\xi$ is stable in this
sense, the vector bundle
$\ad \xi$ may only be strictly semistable. In  case $G$ is simply
connected and $C$ is an elliptic curve $E$, there are essentially no  
properly stable $G$-bundles over $E$.

In case $G$ is semisimple, we have the following variant of Theorem 2.2:

\theorem{2.3} Suppose that $G$ is semisimple. Let $\xi\to C$ be a
principal $G$-bundle. Then the following are equivalent\rom:
\roster
\item"{(i)}" $\xi$ is semistable.
\item"{(ii)}" For every finite-dimensional vector space $V$ and
representation
$\phi\: G \to GL(V)$, the vector bundle $V(\xi)$  associated to
$\xi$ via $\phi$ is semistable.
\qed
\endroster
\endstatement

\ssection{2.3. Families of principal bundles and S-equivalence.}
 
Recall the following standard terminology: 

\definition{Definition 2.4} A {\sl family\/} of principal 
$G$-bundles over the curve $C$ parametrized by a complex space (or
scheme) $S$ is a principal $G$-bundle
$\Xi$ over $C\times S$. The family $\Xi$ is a family of {\sl
semistable\/}   principal $G$-bundles over $C$ if $\Xi|C\times
\{s\}=\Xi _s$ is semistable for all $s\in S$.  Finally, let
$\xi$ and $\xi'$ be two semistable bundles over
$C$. We say that $\xi$ and
$\xi'$ are {\sl S-equivalent\/} if there exists a family of semistable
bundles
$\Xi$ parametrized by  an irreducible $S$ and a point $s\in S$ such
that, for
$t\neq s$, $\Xi |C \times \{t\} \cong \xi$ and $\Xi |C \times
\{s\} \cong
\xi'$. More generally, we let S-equivalence be the equivalence relation
generated by the above relation.
\enddefinition

\ssection{2.4. The Harder-Narasimhan reduction.}

Let $\xi$ be an unstable $G$-bundle over the smooth curve $C$. In this 
case, the structure group of $\xi$ reduces canonically to a parabolic
subgroup $P$ of $G$, the {\sl Harder-Narasimhan parabolic\/} associated
to $\xi$ \cite{26} or \cite{2}, pp\. 589--590. The construction is as
follows.  On $V =\ad \xi$, there is the nondegenerate skew bilinear
form $[\cdot, \cdot]$ induced by Lie bracket. The Harder-Narasimhan
filtration of the vector bundle $V$ is self-dual in the obvious sense,
under the isomorphism from $V$ to
$V\spcheck$ induced by the bracket. Thus, if
$$\{0\} \subset V_0 \subset \cdots \subset V_k = V$$ is the
Harder-Narasimhan filtration, then $V_i = V_{k-i}^\perp$. After
re-indexing, we can write the filtration as 
$$\{0\} \subset V_{-n} \subset \cdots \subset V_n = V,$$  where if the
filtration has odd length we set $V_0 = V_{-1}$ (we will see that the
filtration cannot in fact have odd length in a minute). Here
$\mu (V_k/V_{k-1})=
\mu _k$ is a strictly decreasing sequence of integers, with
$\mu _{-k} = -\mu _k$, and  $V_{-i}=V_{i-1}^\perp$. Since, for $i,
j\geq 0$, 
$$\mu (V_{-i}\otimes V_{-j}) = \mu _{-i}+\mu _{-j}\geq \mu _{-i},$$ it
follows that $V_0$ is a bundle of subalgebras and that $V_{-1}$ is a
subbundle  of nilpotent ideals of $V_0$. Note that if $V_0 = V_{-1}$,
then $V_0$ is a nilpotent subbundle of $V$ of rank $\frac12\dim \frak
g$. But every nilpotent subalgebra of $\frak g$ has dimension strictly
less that $\frac12\dim
\frak g$ (the maximal possible such dimension is the dimension of the
nilpotent radical of   a Borel subalgebra, which is the number of
positive roots). Thus $V_0
\neq V_{-1}$. Moreover $V_{-1}^\perp = V_0$. 

Fixing a point $x$ of $C$, we can choose an identification of $\xi_x$
with $G$ and hence we can identify
$V_x$ with
$\frak g$. Via this identification,
$(V_0)_x$ is identified with a subalgebra $\frak p$ of $\frak g$, and
$(V_{-1})_x$ with a nilpotent ideal $\frak u = \frak u_{-1}$ of $\frak
p$. By
\cite{7, Chap\. 7,8, VIII, \S 10, Th\'eor\`eme 1}, $\frak p$ is a
parabolic subalgebra of $\frak g$ and the proof of the theorem shows
moreover that $\frak u$ is the nilpotent radical of $\frak p$.  Let
$P$ be the parabolic subgroup of $G$ corresponding to $\frak p$. The
map $t\mapsto (V_0)_t$ defines a holomorphic section of the bundle
associated to $\xi$ via the action of $G$ on the space of all parabolic
subalgebras conjugate to $\frak p$, namely $G/N_G(\frak p)$, and thus a
reduction of the structure group of $\xi$ to $N_G(\frak p)$. The
normalizer of $\frak p$ is contained in the normalizer
$N_G(P)$ of $P$ in $G$. By \cite{19}, $N_G(P)=P$ and thus $N_G(\frak p)
= P$. This construction thus reduces  the structure group of $\xi$ to
$P$ (compare \cite{18}). We call $P$ (which is well-defined up to
conjugation) the {\sl Harder-Narasimhan parabolic\/} associated to $P$
and the 
$P$-bundle $\xi _P$ the {\sl Harder-Narasimhan reduction\/} of $\xi$.
In particular, every unstable bundle has a canonical reduction of
structure group to a parabolic subgroup of $G$. In fact, it is easy to
check that, in the above notation, for every $i> 0$, $\det
(V_{-i}/V_{-i-1})$ is a line bundle of positive degree over $E$
associated to a dominant character $\chi$ of $P$, and thus $\xi$ is
unstable by (v) of Theorem 2.2. In case $\xi$ is a semistable bundle
over $G$, we set $P=G$ and $\xi_P = \xi$.

Let $\frak u_{-i}$ correspond to $(V_{-i})_x$. For $i\geq 0$, $\frak
u_{-i}$ is an ideal in $\frak u$ and in $\frak p$. Furthermore, $\frak
u_{-n}$ lies in the center of $\frak u$. Likewise $\frak u_{-i}/\frak
u_{-i-1}$ lies in the center of $\frak u/\frak u_{-i-1}$. Let $U_{-i}$
be the unipotent subgroup of $P$ corresponding to $\frak u_{-i}$.
Following our earlier notation, we let $U = U_{-1}$; it is the
unipotent radical of $P$. Let $L = P/U$. Then $P$ is isomorphic to a
semidirect product of $L$ and $U$. A subgroup of $P$ isomorphic to $L$
under the projection (well-defined up to conjugation) is called a {\sl
Levi
  factor\/} and also denoted by $L$. The group $L$ acts by conjugation
on $U$ and $U_{-i}$. Thus there is an induced action of $L$ on
$U_{-i}/U_{-i-1}$. Let $\eta$ be the principal $L$-bundle induced from
$\xi_P$ by taking the projection from $P$ to $L$. We can summarize the
above as follows:

\theorem{2.5} Let $\xi$ be an unstable $G$-bundle. Then associated to
$\xi$ is a parabolic subgroup $P$ of $G$, canonically defined  up to
conjugation, and a principal
$P$-bundle $\xi_P$ such that $\xi_P\times _PG \cong \xi$. The bundle
$\xi_P$ is uniquely determined up to isomorphism of $P$-bundles. If
$U$ is the unipotent radical of $P$, $L=P/U$ is a Levi factor of
$P$, and $\xi_L=\xi_P/U$ is the associated $L$-bundle, then $\xi_L$ is
a semistable principal $L$-bundle for the reductive group $L$, uniquely
determined up to isomorphism of
$L$-bundles. Finally, the identity component of
$\Aut_G\xi$ is equal to the identity component of $\Aut_P\xi_P$.
\endstatement
\proof We have already discussed the construction of $\xi_P$ and hence
of $\xi_L$. To see that $\xi_L$ is semistable, note that $\ad \xi_L$ is
the $L$-bundle associated to the representation of $L$ on $\frak
p/\frak u$, which is the vector bundle $V_0/V_{-1}$. By construction,
$V_0/V_{-1}$ is semistable, and hence $\xi_L$ is semistable as well.

To see the final statement, since the bundles $V_k/V_{k-1}$ are
semistable with negative slope for $k > 0$, $h^0(E; \ad \xi) = h^0(E;
V_0) = h^0(E; \ad _P\xi_P)$, and this implies the statement about the
identity components.
\endproof

\ssection{2.5. The case of genus zero and one.}

In case $C$ has genus $0$ or $1$, we can further reduce the structure
group of the Harder-Narasimhan reduction $\xi_P$ to that of a Levi
factor $L$, viewed as a subgroup of
$P$. This statement for genus
$1$ is then a preliminary analogue of (i) of Theorem 1.6.

\proposition{2.6} Let $\xi$ be an unstable $G$-bundle over a smooth
projective curve $C$ and let $P$ be the Harder-Narasimhan parabolic
subgroup of $G$ associated to $\xi$. Suppose that         the genus of
$C$ is $0$ or $1$. Then the structure group of $\xi$ reduces from $P$
to a Levi factor $L$. In fact, there exists a maximal parabolic
subgroup $P_0$ of $G$ such that the structure group of $\xi$ reduces to
a Levi factor $L_0$ of $P_0$, and such that, for every nonzero dominant
character $\chi\: P_0 \to \Cee^*$, the line bundle associated to the
$L_0$-bundle $\xi$ and the character $\chi$ has positive degree.  
\endstatement 
\proof Let $\xi_P$ be the Harder-Narasimhan reduction, and let $\eta$
be the $L$-bundle induced by $\xi_P$. To prove the first part of the
proposition, it suffices to show that the set of isomorphism classes of
lifts of the $L$-bundle $\eta$ to a
$P$-bundle consists of a single element, which then must be the bundle
induced by the inclusion of $L$ as a subgroup of $P$. Such lifts are
classified by the nonabelian cohomology set $H^1(C; U(\eta))$, where
$U(\eta)$ is the sheaf of groups associated to the principal bundle
$\eta$ via the action of $L$ on $U$. Likewise, the set of isomorphism
classes of lifts of $\eta$ to a $LU_{-i}$-bundle is classified by the
nonabelian cohomology set $H^1(C; U_{-i}(\eta))$. There is an exact
sequence
$$0 \to U_{-i-1} \to U_{-i} \to U_{-i}/U_{-i-1} \to 0,$$ where
$U_{-i}/U_{-i-1}$ is an abelian unipotent group, and thus a vector
space. By taking the associated long exact sequence of pointed sets (see
\cite{17}), it suffices to show that
$H^1(C; (U_{-i}/U_{-i-1})(\eta)) = \{0\}$ for all $i\geq 1$, where
$(U_{-i}/U_{-i-1})(\eta))$ is the vector bundle $\eta\times
_L(U_{-i}/U_{-i-1})$ associated to $\eta$ and to the representation of
$L$ on the vector group $U_{-i}/U_{-i-1}$. But it is easy to see  that
the vector bundle
$(U_{-i}/U_{-i-1})(\eta)$ is simply the vector bundle
$V_{-i}/V_{-i-1}$, which is  a semistable bundle with slope $\mu _i >
0$. By Serre duality, 
$$h^1(C; V_{-i}/V_{-i-1}) = h^0(C;
\left(V_{-i}/V_{-i-1}\right)\spcheck \otimes K_C).$$ Since $g(C) \leq
1$, $\deg K_C \leq 0$, and thus
$\left(V_{-i}/V_{-i-1}\right) \spcheck \otimes K_C$ is a semistable
bundle of negative degree. Hence $ h^0(C;
\left(V_{-i}/V_{-i-1}\right)\spcheck \otimes K_C) = 0$, and so $h^1(C;
V_{-i}/V_{-i-1}) = 0$ as well. It follows that the unique lift of
$\eta$ to a $P$-bundle is given by the bundle induced by the inclusion
of $L$ in $P$. This proves that the structure group of $\xi$  reduces
to $L$.

Let $\xi _P$ be the Harder-Narasimhan reduction of $\xi$. By the
construction of $\xi _P$, there exists a dominant character
$\chi\: P\to \Cee^*$ such that the degree of the line bundle
corresponding to $\xi _P$ is positive. Suppose that $P$ contains a
Cartan subgroup $H$ and a Borel subgroup $B$ corresponding to a choice
of positive roots, or equivalently a basis $\Delta$ of simple roots.
Then the parabolic subgroups containing $B$ are in $1-1$ correspondence
with proper subsets $\Sigma$ of $\Delta$, and, if $P$ corresponds to
$\Sigma$, then the character group of $P$ is the group generated by the
fundamental weights $\varpi_\alpha$ for $\alpha \in
\Delta - \Sigma$. Since $\chi$ is a positive combination of the
fundamental weights $\varpi_\alpha$ with $\alpha \in \Delta -\Sigma$,
there exists a fundamental weight $\varpi_\alpha$, $\alpha \in \Delta
-\Sigma$, such that the line bundle associated to the character
$\varpi_\alpha$ of $P$ has positive degree. Let $P_0$ be the maximal
parabolic corresponding to the subset $\Delta -\{\alpha\}$. Then $P_0$
contains $P$, and there exists a Levi factor $L_0$ of $P_0$ containing
$L$. Every dominant character of $P_0$ is a positive multiple of
$\varpi_\alpha$, and thus the extension of $\xi _P$ to a $P_0$-bundle
is as claimed.  \endproof

We remark that the genus zero case of Proposition 2.6 leads to a quick
proof of Grothendieck's theorem that, for $G$ a connected reductive
group, every holomorphic
$G$-bundle
$\xi$ over $\Pee^1$ has a holomorphic reduction of structure group to a
Cartan subgroup. The proof is by induction on the rank of the derived
subgroup of $G$. If this rank is zero, then $G$ is an algebraic torus
and there is nothing to prove. If the rank is positive and $\xi$ is not
semistable, then we can reduce to the case of smaller rank by
Proposition 2.6. Thus we may assume that $\xi$ is semistable.
Furthermore, it suffices to show that the induced bundle on the
quotient of $G$ by the identity component of the center is trivial, for
then by the long exact cohomology sequence the structure group of $\xi$
reduces to the identity component of the center. Thus we may assume
that $\xi$ is semistable and that $G$ is semisimple, and will show that
$\xi$ is the product bundle. If $G= SL_n(\Cee)$, then the result is
immediate from Riemann-Roch, which guarantees that the rank $n$ vector
bundle associated to $\xi$ has a section, and induction on $n$. In
general, there is a faithful representation of $G$ in $SL_N(\Cee)$ for
some
$N$, and by (ii) of Theorem 2.3 the associated bundle $\xi \times
_GSL_N(\Cee)$ is trivial. By the long exact cohomology sequence, $\xi$
is in the image of $H^0(\Pee^1; \underline{SL_N(\Cee)/G})$, where
$\underline{SL_N(\Cee)/G})$ denotes the sheaf of sets on $\Pee^1$
whose sections over an open set $U$ are the  holomorphic maps from
$U$ to the quotient variety $SL_N(\Cee)/G$. But since $SL_N(\Cee)$ is
affine and $G$ is reductive, the quotient $SL_N(\Cee)/G$ is affine and
hence every map from $\Pee^1$ to $SL_N(\Cee)/G$ is constant. It then
follows easily that the image of $H^0(\Pee^1; \underline{SL_N(\Cee)/G})$
in
$H^1(\Pee^1; \underline{G})$ is trivial, so that $\xi$ is the trivial
bundle.
\medskip

Let us summarize what we showed so far:

\theorem{2.7} Let $\xi$ be a holomorphic $G$-bundle over an elliptic
curve. If $\xi$ is not semistable, then there is a holomorphic
reduction of  $\xi$ to a semistable $L$-bundle $\xi_L$, where $L$ is a
Levi factor of the Harder-Narasimhan parabolic subgroup associated to
$\xi$. The parabolic subgroup $P$ and hence $L$ are determined up to
conjugation in $G$, and the bundle $\xi_L$ is determined up to
isomorphism of
$L$-bundles. Finally, the identity component of $\Aut_G \xi$ is equal
to the identity component of $\Aut _P\xi_P$, where $\xi_P$ is the
Harder-Narasimhan reduction of $\xi$.
\qed
\endstatement

In case $G=GL_n$, the reduction of the bundle $\xi$ to an $L$-bundle
corresponds to the direct sum decomposition $V = \bigoplus
_{i=1}^kV_i$, where each $V_i$ is semistable and
$\mu(V_i) > \mu (V_{i+1})$. In this case $L = GL_{d_1} \times \cdots
\times GL_{d_k}$, where $d_i$ is the rank of $V_i$.

\ssection{2.6. Some comments on the Jordan-H\"older filtration.}

If $V$ is a semistable vector bundle over the smooth algebraic curve
$C$, then $V$ also has a canonical filtration. Recall that a semistable
vector bundle is {\sl polystable\/} if it is a direct sum of stable
bundles (necessarily all of the same slope). It is easy to check that,
if $V$ is semistable, then there is a largest polystable subbundle
$F^0$ of $V$. (If $W_1$ is a polystable subbundle of $V$ and $W_2$
is a stable subbundle of $V$ with $\mu(W_1) = \mu (W_2) = \mu (V)
=\mu$, then $\mu (W_1+W_2) \leq \mu$ and thus $\mu (W_1\cap W_2) \geq
\mu(W_i)$. By stability, either $W_1 \cap W_2 =W_2$ or $W_1\cap W_2 =
0$. Thus, if $W_2$ is not contained in $W_1$, the map $W_1\oplus W_2
\to V$ is injective, and its image is a subbundle.) Applying the same
construction to
$V/F^0$ and using induction on the rank produces a canonical filtration
$F^i$ of
$V$ with
$F^{i+1}/F^i$ polystable. Moreover $V$ is S-equivalent to
$\operatorname{gr}F^i$, which is polystable. Of course, if $V$ was
polystable to begin with, the filtration is trivial. This construction
is generalized to an arbitrary reductive group $G$ in \cite{28}.

In case $C=E$ is an elliptic curve, it follows from Theorem 1.3 that
every polystable vector bundle of degree zero is a direct sum of line
bundles. Thus the above construction reduces the structure group of a
semistable $SL_n$-bundle to the Borel subgroup. In fact, the structure
group can be reduced to a much smaller subgroup. We will see in the
next section how to do this for a general simply connected group.

\ssection{2.7. The Atiyah-Bott point.}

In this subsection, we introduce an invariant of a $G$-bundle $\xi$,
due to Atiyah-Bott. Then we apply this invariant to solve the following
problem: suppose that $L$ is a Levi factor of a parabolic subgroup of
$G$, and that $\xi_L$ is a holomorphic reduction of $\xi$ to $L$ which
is semistable as an $L$-bundle. When is $\xi_L$ the $L$-bundle arising
from the Harder-Narasimhan reduction of $\xi$? Throughout this
subsection, we work over an elliptic curve $E$, although much of the
discussion generalizes to curves of higher genus.

Let $L$ be a reductive group and let $\Cal C$ be the identity component
of the  center of
$L$. Let $\overline {\Cal C}$ be the quotient of
$L$ by its derived subgroup $DL$. Then $\Cal C\to \overline {\Cal C}$
is a finite covering. Let $\frak c$ be the Lie algebra of $\Cal C$,
which we can view as the universal cover of $\Cal C$ and of $\overline
{\Cal C}$ and let
$\Lambda_{\Cal C}\subseteq \Lambda_{\overline {\Cal C}}\subseteq \frak
c$ be the lattices $\pi_1(\Cal C)$ and $\pi_1(\overline {\Cal C})$
respectively. Let $\xi$ be  a semistable $L$-bundle. We define the {\sl
Atiyah-Bott point} $\mu(\xi)\in \Lambda_{\overline {\Cal C}}\subset
\frak c$ as follows. We have the $\overline {\Cal C}$-bundle $\det
\xi =\xi/DL$. We define
$$\mu(\xi)=c_1(\det \xi )\in H^2(E; \Lambda_{\overline {\Cal C}}) \cong
\Lambda_{\overline {\Cal C}}\subset
\frak c,$$ where $c_1$ is the analogue of the first Chern class arising
from the coboundary map in the long exact sequence
$$0 \to \Lambda_{\overline {\Cal C}} \to \frak c \otimes _\Cee\scrO_E
\to
\underline{\overline {\Cal C}} \to 0.$$
 Notice that
$\mu(\xi)$ is an element of the Lie algebra of $L$ whose image under
the exponential map is  of  finite order.

\lemma{2.8} Suppose that $L\subseteq L'$ is an inclusion of reductive
groups, that $\xi$ is a semistable $L$-bundle and that
$\xi'=\xi\times_LL'$ is also semistable. Then under the inclusion
$\frak l\subseteq \frak l'$  the image of $\mu(\xi)\in \frak l'$  lies
in the Lie algebra of the center of $L'$ and is the Atiyah-Bott point
$\mu(\xi')$.
\endstatement

\proof Let ${\Cal C}$, resp\. $\Cal C'$, be the center of $L$, resp\.
$L'$, and let
$\overline {\Cal C}$, resp\. $\overline {\Cal C}'$ be the quotient of
$L$, resp\. $L'$, by its derived subgroup $DL$, resp\. $DL'$. We
decompose $\frak l'$ under the action of the semisimple element
$\mu(\xi)\in \frak l$:
$$\frak l'=\bigoplus_q(\frak l')^q,$$ where $(\frak l')^q$ is the
$q$-eigenspace  for the adjoint action of $\mu(\xi)$ and $q\in \Q$
(since $\mu(\xi)$ is a rational point of the Lie algebra of $\overline
{\Cal C}'$). We have the associated vector bundle
$(\frak l')^q(\xi)$, and its degree is given by $\dim (\frak l')^q\cdot
q$. By definition, since
$\xi'$ is semistable, so is 
$\ad(\xi')=\xi\times_L\frak l'$. This means that all of the nonzero
summands $(\frak l')^q(\xi)$ must have the same slope, but since $\frak
l \subseteq(\frak l')^0$, there is a nonzero summand with slope zero.
Thus the only possibility  is that
$(\frak l')^0=\frak l'$. Hence,
$\mu(\xi)$ is contained in the center of $\frak l'$, proving the first
statement.

Clearly, the inclusion $L\subseteq L'$ induces an inclusion $DL\subseteq
DL'$ and hence there is an induced mapping $\overline {\Cal C}\to
\overline {\Cal C}'$. Under this mapping the bundle $\det \xi$ maps to
$\det(\xi')$. This proves that the Atiyah-Bott points for
$\xi$ and $\xi'$ correspond under this mapping. 
\endproof

\lemma{2.9} Let $\xi$ be a holomorphic $G$-bundle. Suppose that $\xi$
reduces to a semistable bundle $\xi_L$, where $L$ is a reductive
subgroup of $G$. Then there exists a Harder-Narasimhan parabolic
subgroup $P$ for $\xi$  and a Levi factor $L'$ of $P$ which contains
$L$, and such that the Harder-Narasimhan reduction $\xi_{L'}$ of $\xi$
is
$\xi_L\times _LL'$.
\endstatement

\proof Consider the adjoint action of $\mu(\xi_L)$ on $\frak g$. It is a
diagonalizable action all of whose eigenvalues are rational numbers. We
define 
$\frak g_0(\mu(\xi_L))\subseteq \frak g$ to be the subspace which is
the direct sum of the eigenspaces corresponding to  nonnegative
eigenvalues of this action. It follows from the proof of Lemma 2.8 that
$\frak g_0(\mu(\xi_L))$  is the Lie algebra of the Harder-Narasimhan
parabolic for
$\xi_L\times_LG$ and that the subalgebra of $0$-eigenspaces is the Lie
algebra of the Levi factor $L'$ of this parabolic. (This also includes
the case where $P=G$.) Clearly,
$L\subseteq L'$. It follows from the Harder-Narasimhan construction
that $\xi_{P}=\xi_L\times_LP$  and hence that
$\xi_{L'}=\xi_L\times_LL'$.
\endproof

\corollary{2.10} Let $\xi$ be a holomorphic $G$-bundle. Suppose that
there are two reductions of $\xi$  to semistable bundles $\xi_L$ and
$\xi_{L'}$ over reductive subgroups $L,L'$ of $G$. Then the Atiyah-Bott
points
$\mu(\xi_L)\in\frak l\subseteq \frak g$ and $\mu(\xi_{L'})\in \frak
l'\subseteq \frak g$ are conjugate points of
$\frak g$ under the adjoint action of $G$.
\endstatement
\proof By Lemma 2.9, there exists a Harder-Narasimhan parabolic $P$ and
a  Levi factor $L''$ of $P$ which contains $L$. By Lemma 2.8,
$\mu(\xi_L) = \mu(\xi_{L''})$, where $\xi_{L''} =
\xi
\times _LL''$. The corollary now follows from the uniqueness of the
Harder-Narasimhan reduction up to conjugation in
$G$.
\endproof

Thus, we define the Atiyah-Bott point of any holomorphic $G$-bundle for
any reductive group $G$, by choosing a reductive subgroup
$L\subseteq G$ and a semistable $L$-bundle $\xi_L$ reducing $\xi$ and
taking the image of the Atiyah-Bott point of $\xi_L$ under the
inclusion $\frak l\subseteq \frak g$. For example, we could always
choose the reduction of
$\xi$ to the Levi factor of a Harder-Narasimhan parabolic. By the above
corollary the result is well-defined up to conjugation by an element of
$G$, and the conjugacy class will be denoted $\mu(\xi)$. Of course,
fixing a Cartan subalgebra
$\frak h$ for $G$ and a Weyl chamber $C_0$ of $\frak h$, we can
identify $\mu(\xi)$ once and for all with a point of
$C_0$. Note that the point $\mu(\xi)$ satisfies an integrality
condition similar to that of Corollary 1.7 (ii), namely that $\mu(\xi)
\in \Lambda_{\overline{\Cal C}} \subseteq \Lambda _{\Cal C}\otimes \Bbb
Q$.

Quite generally, suppose that $\mu$ is a  semisimple element of the Lie
algebra $\frak g$, all of whose eigenvalues are rational. We can break
$\frak g$ into a direct sum $\bigoplus _q\frak g^q$ of eigenspaces for
$\ad\mu$. Let
$P(\mu)$ be the connected group whose Lie algebra corresponds to the
direct sum of the eigenspaces for nonnegative eigenvalues
$q$, and let $L(\mu)$ be the subgroup of $P(\mu)$ corresponding to the
eigenspace for $q=0$. It is easy to see that $P(\mu)$ is a (closed)
parabolic subgroup of $G$ and that $L(\mu)$ is a Levi factor of
$P(\mu)$. Clearly $\mu$ is contained in the center of $L(\mu)$.

We may then characterize the $L$-bundles arising from the
Harder-Narasimhan reduction as  follows:

\corollary{2.11} Let $\xi$ be a holomorphic $G$-bundle and let
$\mu(\xi)$ be the  associated Atiyah-Bott point. Let $L$ be a reductive
subgroup of 
$G$. Suppose that there is a holomorphic reduction of $\xi$ to a 
semistable $L$-bundle $\xi_L$.  Then
$\xi_L$ is the
$L$-bundle arising from the Harder-Narasimhan reduction of $\xi$ if and
only if $L =  L(\mu(\xi))$, if and only if the roots of $G$ vanishing on
$\mu(\xi)$ are exactly the roots of $L$.
\qed
\endstatement

\definition{Definition 2.12} Let $L$ be a reductive subgroup of $G$,
let $\frak c$ be the Lie algebra of the center of $L$, and let
$\mu\in \frak c$. We say that the pair
$(L, \mu)$ is {\sl saturated\/} if the roots  of $G$ vanishing on
$\mu$ are exactly the roots of $L$. With $L$ as above, let 
$\xi_L$ be a semistable $L$-bundle. We say that the pair $(L, \xi_L)$
is {\sl saturated\/} if the pair $(L, \mu(\xi_L))$ is saturated. In
this case, by Corollary 2.11, $L$ is a Levi factor of the
Harder-Narasimhan parabolic associated to $\xi_L\times _LG=\xi$, and
$\xi_L$ is the $L$-bundle associated to the Harder-Narasimhan reduction
of $\xi$.
\enddefinition

For example, for $G=GL_n$ and $L = GL_{d_1} \times \cdots \times
GL_{d_k}$, a reduction of a vector bundle $V$ to an $L$-bundle is a
direct sum decomposition 
$V = \bigoplus _{i=1}^kV_i$ with $\operatorname{rank}V_i = d_i$. The
$L$-bundle $\xi_L$ is semistable if and only if  each
$V_i$ is semistable. Finally, in this case the pair $(L, \xi_L)$ is
saturated if and only if
$\mu(V_i) \neq \mu(V_j)$ for $i\neq j$.

\section{3. A normal form for the transition functions of a semistable
bundle.}

This section is concerned with the description of all semistable
bundles within a fixed S-equivalence class. Let $\xi$ be a fixed
semistable $G$-bundle. It is a well-known result that, for every smooth
projective curve $C$, there is a unique central Yang-Mills
representative $\xi_0$ S-equivalent to $\xi$. This representative is
associated to a representation $\rho$ from a central extension of
$\pi_1(C)$ to $K$. In case $C=E$, to find a normal form for the
transition functions of $\xi$, we begin by finding a  canonical form
for the transition functions of the corresponding Yang-Mills bundle
$\xi_0$ (which we call normalized transition functions). The canonical
form for the transition functions of $\xi$ is then obtained by
multiplying the normalized transition functions for t$\xi_0$  by
elements of the form $\{\exp(f_{ij}X)\}$ where $\{f_{ij}\}$ is a
$\scrO_E$-valued $1$-cocycle giving the standard generator for
$H^1(E;\scrO_E)$ and $X$ is a nilpotent element in 
$\frak z_{\frak g}(\rho)$, the subalgebra of $\frak g$ consisting of
elements which commute with the representation $\rho$.

We shall use the  following notation throughout the rest of this paper.
For any subset
$S\subseteq G$ we denote by
$Z_G(S)$ the centralizer of
$S$ in $G$ and by ${\frak z}_{\frak g}(S)$ its Lie algebra. Clearly
${\frak z}_{\frak g}(S)$ is the subalgebra of $\frak g$ invariant under
$\operatorname{Ad} S$. Given a representation
$\rho\:
\Gamma _\Ar \to G$, we define $Z_G(\rho)$ to be
$Z_G(\operatorname{Im} \rho)$, and similarly for $\frak z_{\frak
g}(\rho)$.

\ssection{3.1. Flat and Yang-Mills connections.}

Let $C$ be a smooth projective curve and let $K$ be a compact,
semisimple Lie group. Given a representation
$\rho\: \pi _1(C)\to K$, we can form the associated principal
$K$-bundle $\xi_0 = (\tilde C \times K)/\pi _1(C) \to C$, where $\pi
_1(C)$ acts on $\tilde C$, the universal cover of $C$, in the usual
way, and on $K$ via $\rho$.  We shall call such a
$K$-bundle a {\sl flat\/} $K$-bundle.  The inclusion of $K$ in $G$ then
yields a  holomorphic $G$-bundle $\xi _0\times _KG$. We shall refer to a
$G$-bundle which arises in this way as a flat $K$-bundle or, in the
case where $C=E$ is an elliptic curve and $G$ is simply connected, as a
{\sl split\/} bundle. (We shall also blur the distinction between the
$K$-bundle $\xi_0$ and the $G$-bundle $\xi _0\times _KG$.) The flat
structure on a flat $K$-bundle 
$\xi _0$ is equivalent to a  $K$-connection $A$ on
$\xi _0$ with curvature $F_A=0$, modulo gauge equivalence. There is an
induced connection on
$\xi _0\times _KG$, which we shall also denote by $A$, and the
holomorphic structure on
$\xi _0\times _KG$ is induced by the $(0,1)$-connection which is the
$(0,1)$ part of $A$. 

More generally, suppose that $K$ is compact but not necessarily
semisimple. In general we will let
$K =\Cal C\times _FD$, where $\Cal C = \Cal CK$ is the identity
component of the center of $K$ and $D$ is the derived group. Let $G$ be
the complexification of $K$; it is a reductive algebraic group, and
conversely every reductive group arises in this way. Let
$\Gamma _\Ar (C)$ be the central extension of
$\pi_1(C,*)$ defined in \S 6 of Atiyah-Bott \cite{2}, generalizing the
construction we have given in Section 1 for $C=E$. Thus there is a
central extension
$$0 \to \Ar \to \Gamma _\Ar (C) \to \pi_1(C,*) \to \{1\},$$ and it is
unique specified by the condition that, if $\alpha_i, \beta_i$ are the
standard generators of $\pi_1(C,*)$ and $A_i, B_i$ are any lifts of the
$\alpha_i,
\beta_i$ to elements of $\Gamma _\Ar (C)$, then
$$\prod _i[A_i, B_i] = 1\in \Ar.$$ It is straightforward to check that
the class of this central extension is the element $1\in H^2(C; \Ar)
\cong \Ar$. Let
$\nu\to C$ be the unique $C^\infty$
$U(1)$-bundle with Chern class
$1$ and let $\tilde C\to C$ be the universal cover of $C$. Then  $\nu
\times _C\tilde C$ is a $\pi_1(C) \times U(1)$-bundle, and one can
check that there is a unique lift of this bundle to a bundle
$\boldsymbol \Gamma_\Ar \to C$ with structure group $\Gamma _\Ar (C)$.
Fix a Yang-Mills (i.e. harmonic) connection $A_\nu$ on $\nu$. Then
$A_\nu$ lifts to a connection on $\boldsymbol \Gamma_\Ar \to C$, which
we denote by $A_0$. 

Given the compact Lie group $K$, a representation $\rho\:\Gamma_\Ar \to
K$ will be called {\sl central\/} if $\rho(\Ar)$ is contained in the
center of $K$. Recall that a Yang-Mills connection on a $K$-bundle is a
$K$-connection $A$ whose curvature
$F_A$ is covariant constant. Given a central representation
$\rho\:\Gamma _\Ar (C) \to K$,  there is the induced $K$-bundle
$\boldsymbol \Gamma_\Ar
\times _{\Gamma _\Ar (C)}K$, and it has a  Yang-Mills connection
$\rho_*(A_0)=A$. Note that
$*F_A$ has values in the Lie algebra of the center of
$K$, not just the Lie algebra of $K$. We will refer to such a
connection as a {\sl central Yang-Mills connection}. Thus, the image of
the fixed Yang-Mills connection
$A_0$ on $\boldsymbol \Gamma_\Ar$ under a central representation is a
central Yang-Mills connection, and conversely,  up to gauge
equivalence, every central Yang-Mills connection arises in this way.
Given a central  Yang-Mills connection $A$, we shall somewhat
inaccurately refer to the corresponding representation
$\rho\: \Gamma _\Ar (C) \to K$ as the {\sl holonomy representation}.
Note that the identification of central representations with central 
Yang-Mills connections depends on the choice of the connection $A_0$.

The inclusion $\Ar \subset \Cee$ defines an inclusion of groups
$\Gamma_\Ar(C) \subset
\Gamma _\Cee(C)$. Clearly there is an exact sequence
$$0 \to \Zee \to \Gamma _\Cee(C) \to \Cee^*\times \pi_1(C) \to 0,$$
compatible with the exact sequence
$$0 \to \Zee \to \Gamma _\Ar(C) \to U(1)\times \pi_1(C) \to 0.$$ The
fixed Yang-Mills connection $A_0$ on $\boldsymbol \Gamma_\Ar$ induces a
holomorphic structure on $\boldsymbol \Gamma_\Cee = \boldsymbol
\Gamma_\Ar\times _{\Gamma_\Ar(C)}\Gamma _\Cee(C)$, and we shall denote
this bundle, together with its holomorphic structure, by  $\boldsymbol
\Gamma_\Cee \to C$, with structure group
$\Gamma _\Cee (C)$. A central representation $\rho\: \Gamma _\Ar(C) \to
K$ induces a holomorphic representation $\rho_\Cee \:\Gamma _\Cee(C)
\to G$. Hence there is an associated holomorphic $G$-bundle
$\boldsymbol \Gamma_\Cee  \times _{\Gamma _\Cee(C)} G$. Clearly
$$\boldsymbol \Gamma_\Cee  \times _{\Gamma _\Cee(C)} G \cong
(\boldsymbol \Gamma_\Ar
\times _{\Gamma _\Ar (C)}K) \times _KG,$$ and the holomorphic structure
on $\boldsymbol \Gamma_\Cee  \times _{\Gamma _\Cee(C)} G$ is easily
seen to be the same as the one on $(\boldsymbol \Gamma_\Ar
\times _{\Gamma _\Ar (C)}K) \times _KG$ induced by taking the
$(0,1)$-component of the Yang-Mills connection $\rho_*(A_0)=A$. We
shall refer to a holomorphic $G$-bundle $\xi = \xi _K\times _KG$  whose
holomorphic structure arises from the $(0,1)$-part of a central
Yang-Mills connection on the $K$-bundle $\xi_K$ as a {\sl Yang-Mills
bundle\/} and we shall also refer to the associated central
representation $\rho$ as the holonomy of $\xi$. We emphasize that, in
our definition, a Yang-Mills bundle arises from a {\sl central\/}
Yang-Mills connection, since we shall never need to consider
non-central Yang-Mills connections.

We will need the following version of the theorem of
Narasimhan-Seshadri \cite{25} and Ramanathan \cite{27} (see also
Atiyah-Bott \cite{2} and Donaldson \cite{8}).

\theorem{3.1} Suppose that $G$ is reductive and that $K$ is a maximal
compact subgroup of $G$.
\roster
\item"{(i)}" Let $\xi_0$ be a Yang-Mills $K$-bundle. Then the induced
$G$-bundle $\xi _0\times_KG$ is semistable.
\item"{(ii)}" Let $\xi\to C$ be a semistable principal
$G$-bundle. Then there is a Yang-Mills $K$-bundle $\xi_0$ whose induced
$G$-bundle is S-equivalent to  
$\xi$.  More precisely, there is a family of semistable principal
$G$-bundles $\Xi$ over $C\times \Cee$, such that, for
$t\neq 0$, $\Xi_t\cong \xi$, and such that $\Xi _0$ is the $G$-bundle
induced by the   Yang-Mills $K$-bundle $\xi _0$. 
\item"{(iii)}" The bundle $\xi_0$ is  unique  up to isomorphism of
Yang-Mills $K$-bundles. More precisely, suppose that $\xi$ and
$\xi'$ are two semistable principal $G$-bundles associated to Yang-Mills
$K$-bundles $\xi_0, \xi_0'$. Then $\xi$ and $\xi'$ are S-equivalent if
and only if $\xi_0$ and $\xi_0'$ are isomorphic as Yang-Mills
$K$-bundles, if and only if their associated holonomy representations
$\Gamma _\Ar\to K$ are conjugate in $K$. \qed
\endroster
\endstatement

The proof of the result, in the form stated here, is essentially
contained in Ramanthan's thesis 
\cite{28, Prop\. 3.12 and 3.15} as well as in Atiyah-Bott \cite{2, Lemma
10.12}. (Although the result is stated for curves of genus at least $2$
in \cite{28}, the proof goes through with minor modifications in the 
case of genus $1$.)

Suppose that $K$ is compact, with $\Cal C$ the identity component of the
center of $K$, and let $\frak c$ be the Lie algebra of $\Cal C$. If
$\xi_0$ is a Yang-Mills bundle corresponding to a central
representation $\rho\: \Gamma_\Ar \to K$, then we have defined
$\mu(\xi) \in \frak c$ in \S 2.5. The following is a straightforward
calculation:

\lemma{3.2} With the above notation, let $d\rho\: \Ar \to \frak c$ be
the differential of the homomorphism $\rho$. Then $\mu(\xi) =
d\rho(1)$. \qed
\endstatement

\ssection{3.2. Semistable $G$-bundles with  Yang-Mills connections over
an elliptic curve.}

  From now on, we restrict attention to an elliptic curve $E$, and set
$\Gamma _\Ar =
\Gamma_\Ar(E)$ and $\Gamma _\Cee = \Gamma_\Cee(E)$. We have the
$C^\infty$
$\Gamma _\Ar$-bundle $\boldsymbol \Gamma _\Ar\to E$ and its
complexification
$\boldsymbol \Gamma _\Cee\to E$. The choice of a fixed harmonic
connection on the $U(1)$-bundle
$\nu \to E$ with Chern class $1$ defines a holomorphic structure on
$\boldsymbol \Gamma _\Cee$, which we also fix once and for all. Let
$\{U_i\}$ be an open Stein cover of
$E$ such that each 
$U_i$ is simply connected and such that $U_i\cap U_j$ is either
connected or empty for all $i\neq j$. Then we can find a holomorphic
trivialization of the bundle $\boldsymbol \Gamma _\Cee$ over the open
cover $\{U_i\}$. Let
$\gamma_{ij}\: U_i\cap U_j \to \Gamma _\Cee$ be the holomorphic
transition functions for this choice of trivialization.

Now suppose that $K = \Cal C\times _FD$ is a compact group, where as
usual $D$  is the derived group and $\Cal C$ is the identity component
of the center. Let
$\Cal C_\Cee$ be the algebraic torus which is the complexification of
$\Cal C$. Suppose that
$\xi_0$ is a holomorphic bundle corresponding to a Yang-Mills
connection with holonomy 
$\rho\colon\Gamma _\Ar\to K$. Then there is the associated holomorphic
representation $\rho_\Cee\: \Gamma _\Cee \to G$. Let $h_{ij} =
\rho_\Cee(\gamma_{ij})$. Clearly, the $h_{ij}$ are holomorphic 
transition  functions  for the bundle $\xi_0$, contained in the image
of 
$\rho$. We shall call such transition functions {\sl normalized\/}
(with respect to the fixed cover
$\{U_i\}$ and the choice of trivializations $\{\gamma_{ij}\}$). Since
the choice of the open cover and trivializations will be made once and
for all, we will suppress this choice from the notation. Note that the
construction of normalized transition functions is functorial with
respect to central homomorphisms of compact Lie groups, i\.e\.
homomorphisms $\pi\: K_1\to K_2$ such that $\pi(\Cal CK_1) \subseteq
\Cal CK_2$, where $\Cal CK_i$ is the identity component of the center
of $K_i$. We then have the following lemma on normalized transition
functions:  

\lemma{3.3} With notation as above, if $\xi_0$ is a holomorphic bundle
corresponding to a Yang-Mills connection with holonomy
$\rho\colon \Gamma _\Ar\to K$, then the structure group of $\xi_0$
reduces to $\operatorname{Im}\rho_\Cee$.  Moreover:
\roster
\item"{(i)}" The elements $h_{ij}$ generate
$\operatorname{Im}\rho_\Cee$ modulo
$\rho_\Cee(\Cee)
\subseteq \Cal C_\Cee$, and thus generate $\operatorname{Im}\rho_\Cee$
modulo $\Cal C_\Cee$.
\item"{(ii)}" If $\{h_{ij}\}$ and $\{h_{ij}'\}$ are two normalized 
transition functions \rom(with respect to the cover $\{U_i\}$\rom)
defining isomorphic bundles, then there exists $g\in K$ such  that
$h_{ij}' = gh_{ij}g^{-1}$ for all $i,j$.
\item"{(iii)}" An element
$X\in \frak g$ satisfies $\operatorname{Ad}(h_{ij})(X) = X$ for all
$i,j$  if and only if $X\in \frak z_{\frak g}(\rho)$. 
\item"{(iv)}" There is a natural inclusion of
$Z_G(\rho)$ in $\Aut _G(\xi _0)$, whose image is the set of
automorphisms of $\xi _0$ defined by functions which are constant with
respect to the given local trivializations over $\{U_i\}$. 
\endroster
\endstatement

\proof Clearly the structure group of $\xi_0$ reduces to
$\operatorname{Im}(\rho_\Cee)$. 

Next we show (i). Using the projection $\Gamma_\Cee \to \pi_1(E)$, we
can take the image $\bar \gamma_{ij}\in \pi_1(E)$ of the transition
functions $\gamma_{ij}$, and  the functions $\bar \gamma_{ij}$ are
transition functions for the $\pi_1(E)$-bundle
$\tilde E
\to E$. Thus, the $\bar\gamma_{ij}$ generate $\pi_1(E)$. It follows
that the images
$\bar h_{ij}$ of the $h_{ij}$ generate the image of the homomorphism
$\bar \rho_\Cee\:
\pi_1(E)
\to DG/F$ induced by $\rho$. From this, (i) is clear.

To see (ii), suppose that $\{h_{ij}\}$ and $\{h_{ij}'\}$ are two
normalized  transition functions defining isomorphic bundles. Let
$\rho$ and
$\rho'$ be the respective holonomy maps. Then $\rho$ and $\rho'$ are
conjugate under $K$: say $\rho' = g\rho g^{-1}$ for some $g\in K$. The
same will be true for $\rho_\Cee$ and $\rho_\Cee'$. Thus $h_{ij}' =
\rho_\Cee'(\gamma_{ij}) = g\rho_\Cee(\gamma_{ij}) g^{-1}
=gh_{ij}g^{-1}$.

To see (iii), by (i), $X\in
\frak g$ is fixed by $\operatorname{Ad}(h_{ij})$ for every
$i,j$ if and only if $X$ is fixed by
$\operatorname{Ad}(\operatorname{Im} \rho)$, if and only if
$X\in \frak z_{\frak g}(\rho)$.

As for (iv), for  $h\in Z_G(\rho)$, define an automorphism
$U_i\times G\to U_i\times G$ by sending $(u,g)$ to $(u,hg)$. Since
$h\in Z_G(\rho)$, $h$ commutes with $h_{ij}$ for all $i,j$. This means
these local automorphism fit together to define an automorphism of the
$G$-bundle $\xi_0$ which clearly is constant with respect to the given
trivialization over $\{U_i\}$. The converse is also easy to see, namely
that any constant automorphism of $\xi_0$ lies in the image of
$Z_G(\rho)\to
\Aut_G(\xi_0)$.
\endproof

We have not shown  that $Z_G(\rho) = \Aut _G(\xi_0)$, since {\it a
priori\/} there might exist nonconstant holomorphic functions
$g_i\colon U_i\to G$ with
$h_{ij}g_ih_{ij}^{-1} = g_j$. We will  return to this point in the next
section and verify that indeed $Z_G(\rho) =
\Aut _G(\xi_0)$.

If $K$ is semisimple, then  $\rho$ factors through a homomorphism
$\pi_1(E) \to K$. Since $\pi_1(E)$ is abelian,
$\operatorname{Im}\rho_\Cee = \operatorname{Im}\rho$ is contained in
$Z_G(\rho)$. In general, however, $\operatorname{Im}\rho$ need not be
contained in $Z_G(\rho)$.

\ssection{3.3. Construction of all semistable $G$-bundles.}

Next we give a general method for constructing
$G$-bundles. If $X\in \frak g$ and $f$ is a holomorphic function on
some complex manifold $M$, then $\exp (fX)$ is a holomorphic map from
$M$ to $G$, and, fixing $f$, $X\mapsto \exp (fX)$ is a holomorphic map
from
$\frak g$ to the space of holomorphic mappings from $M$ to $G$. We will
apply this as follows: for the given Stein cover $\{U_i\}$, let
$\{f_{ij}\}$ be a $1$-cocycle for $\scrO_E$ for the cover $\{U_i\}$
such that
$\{f_{ij}\}$ represents the same element as $d\overline z$ in the
$1$-dimensional vector space $H^1(E; \scrO_E)$. For a fixed $X\in
\frak g$, the function $\exp (f_{ij}X)$ defines an element of $H^1(E;
\underline{G})$ and thus a holomorphic principal bundle over $E$. Doing
this for all $X\in \frak g$, we define a holomorphic family of bundles
over
$E\times \frak g$, such that at 
$0\in \frak g$, the corresponding $G$-bundle is the trivial $G$ bundle
$E\times G$.

More generally, suppose that $\xi_0$  has a Yang-Mills connection with
holomony $\rho\colon \Gamma_\Ar\to K$ and in fact is defined   by
normalized transition functions $\{h_{ij}\}$  with values in
$Z_G(\rho)$. For $X\in\frak z_{\frak g}(\rho)$, we define  transition
functions $h_{ij}\exp (f_{ij}X)\:  U_i\cap U_j \to G$. Since the
$\{h_{ij}\}$ form a $1$-cocycle for
$K\subset G$, since each $h_{ij}$ centralizes $X$, and since the
$\{f_{ij}\}$ are a $1$-cocycle,  it is easy to see that $\{h_{ij}\exp
(f_{ij}X)\}$ is indeed a
$1$-cocycle with respect to the open cover $\{U_i\}$. Performing the
above construction for all
$X\in\frak z_{\frak g}(\rho)$, we obtain a holomorphic family of bundles
$\Xi_0$ over
$E\times  \frak z_{\frak g}(\rho)$, such that the bundle corresponding
to $0\in \frak z_{\frak g}(\rho)$ is $\xi_0$. 

\lemma{3.4} Let $\Xi_0 \to E\times  \frak z_{\frak g}(\rho)$ be the
family of $G$-bundles constructed above.
\roster
\item"{(i)}" The inclusion  $\frak z_{\frak g}(\rho) \subseteq
H^0(E;\ad \xi_0)$ of covariant constant sections is an isomorphism.
Likewise, there is a homomorphism  $H^1(E;
\scrO_E) \otimes  \frak z_{\frak g}(\rho)  \to H^1(E; \ad \xi_0)$
defined by cup product and it is also an isomorphism.
\item"{(ii)}"  The  Kodaira-Spencer map of the family $\Xi_0$ at the
origin is given by  
$$\frak z_{\frak g}(\rho) \to H^1(E; \scrO_E) \otimes  \frak z_{\frak
g}(\rho) = H^1(E; \ad \xi_0),$$  where the first map is defined by
tensoring with the element in
$H^1(E; \scrO_E)$ represented by $\{f_{ij}\}$ and the second isomorphism
is that of \rom{(i)}. In particular, the Kodaira-Spencer map is an
isomorphism, and hence $\Xi_0$ is a semiuniversal deformation of
$\xi _0$ at the origin.
\endroster
\endstatement

\proof To see (i), note that $\ad\xi_0$ is a Yang-Mills vector bundle
for the representation of the compact group $K$ acting on $\frak g$.
Thus, using the maximum principle, a standard argument due to
Narasimhan and Seshadri \cite{25} shows that the holomorphic sections
are the covariant constant sections with respect to the connection, and
these are given by constant elements of
$\frak g$ commuting with the holonomy representation $\rho$.  In
particular,
$H^0(E;\ad\xi_0)$ is identified with
$H^0(E;\scrO_E)\otimes \frak z_{\frak g}(\rho)$. A similar argument via
the Dolbeault isomorphism $H^1(E;\ad\xi_0) \cong H^{0,1}(E;\ad\xi_0)$
and using the existence of the nowhere vanishing $(0,1)$-form $d\bar z$
on $E$ shows that
$H^1(E;\ad\xi_0)$ is identified with $H^1(E;\scrO_E)\otimes \frak
z_{\frak g}(\rho)$ via cup product. 

To see (ii), we first note that, if $\eta$ is an arbitrary $G$-bundle
with transition functions $\{g_{ij}\}$ and if $\{g_{ij}\exp
(tX_{ij})\}$ is a first order deformation of $\eta$, then the $X_{ij}$
satisfy
$$\Ad(g_{kj})X_{ij} + X_{jk} = X_{ik},$$ in other words $\{X_{ij}\}$ is
a $1$-cocycle for $\ad \eta$, and by definition it is the
Kodaira-Spencer class. Applying this remark to $\xi_0$ and the bundle
$\Xi_0$ above, it follows  that the first order term for the transition
functions of $\Xi_0$ along the tangent direction
$X\in \frak z_{\frak g}(\rho)$ is the
$1$-cocycle  $\{f_{ij}X\}$, and clearly this  element projects to
$[\{f_{ij}\}]\otimes X\in H^1(E;\scrO_E)\otimes \frak z_{\frak
  g}(\rho)$, where $[\{f_{ij}\}]$ is the image of $\{f_{ij}\}$ in 
$H^1(E; \scrO_E)$. Thus the identification of $H^1(E; \scrO_E)$ with
$\Cee$ via $[\{f_{ij}\}]$ shows that the Kodaira-Spencer map is  the
map $\frak z_{\frak g}(\rho)\to H^1(E;\scrO_E)\otimes \frak z_{\frak
g}(\rho)$ given by tensoring with the class $[\{f_{ij}\}]$ followed by
the natural inclusion of this tensor product into
$H^1(E;\ad\xi_0)$. By (i), this map is an isomorphism from $\frak
z_{\frak g}(\rho)$ to $H^1(E;\ad\xi_0)$.

Thus, the Kodaira-Spencer map is an isomorphism, which means that
$\Xi_0$ is a semiuniversal deformation of $\xi_0$ at the origin.
\endproof

\corollary{3.5} Let $\xi$ be a semistable principal $G$-bundle. Let
$\xi_0$ be the unique bundle with a Yang-Mills $K$-connection
S-equivalent to $\xi$, and let $\{h_{ij}\}$ be normalized transition
functions for
$\xi_0$. Then there exists an  $X\in\frak z_{\frak g}(\rho)$, such that
$\xi$ is trivialized on the open cover $\{U_i\}$ with  transition
functions $\{h_{ij}\exp (f_{ij}X)\}$.
\endstatement

\proof  Let $S$ be a semiuniversal deformation of the bundle $\xi_0$,
with    
$0\in S$ corresponding to $\xi_0$. By the theorem of
Narasimhan-Seshadri, Ramanathan, Donaldson et al (Theorem 3.1), in
every neighborhood of
$0\in S$ there  exists a bundle isomorphic to $\xi$. Taking 
$S$ to be a neighborhood of the origin in $\frak z_{\frak g}(\rho)$ and 
the family to be the family $\Xi_0$  constructed above, the result now
follows immediately from Lemma 3.4.
\endproof 

The main result of this subsection is that we can further assume that 
$X\in
\frak z_{\frak g}(\rho)$ given above is {\sl nilpotent\/}.

\theorem{3.6} Let $\xi$ be a semistable principal $G$-bundle, let
$\xi_0$ be the unique bundle with a Yang-Mills $K$-connection
S-equivalent to $\xi$, and let $\{h_{ij}\}$ be  the normalized
transition functions for $\xi _0$ with respect to the open covering
$\{U_i\}$. Then there exists a holomorphic trivialization of $\xi$ over
the cover $\{U_i\}$ with  holomorphic transition functions of the form
$\{h_{ij}\exp (f_{ij}X)\}$, where
$X\in
\frak z_{\frak g}(\rho)$ is nilpotent.
\endstatement

\proof Begin with transition functions $\{h_{ij}\exp (f_{ij}X)\}$ for
$\xi$ as in  Corollary 3.5, where $X$ is an arbitrary element of $\frak
z_{\frak g}(\rho)$. Note that modifying the
$1$-cocycle
$\{h_{ij}\exp (f_{ij}X)\}$ by the $1$-coboundary which is given by a 
constant element $g$ of
$Z_G(\rho)$ has the effect of leaving the
$h_{ij}$ invariant and replacing $X$ by $(\operatorname{Ad}g)(X)$. Write
$X= X_s+X_n$, where $X_s,X_n\in \frak z_{\frak g}(\rho)$,
$\ad X_s$ is semisimple, $\ad X_n$ is nilpotent, and $[X_s, X_n] =0$.
Fix a Cartan subalgebra $\frak h(\rho)$ of $\frak z_{\frak g}(\rho)$,
which we assume to be the complexification of the  Lie algebra of a
maximal torus $T(\rho)$ in the centralizer of $\rho$ in $K$. After
conjugating by an element of
$Z_G(\rho)$, we may assume that $X_s$ lies in the Cartan subalgebra
$\frak h(\rho)$. Let $H(\rho)$ be the corresponding maximal torus of
$Z_G(\rho)$. Its unique maximal compact subgroup is $T(\rho)=K\cap
H(\rho)$. Let $H'(\rho)$ be  the subgroup of the complex torus $H(\rho)$
consisting of all elements which commute with 
$X_n$. Although $H'(\rho)$ is not necessarily connected, its identity
component  is a subtorus of $H(\rho)$. Let  $T'(\rho)$ be the maximal
compact subgroup of
$H'(\rho)$. Then $T'(\rho)$ is contained in the unique maximal compact
subgroup 
$T(\rho)=K\cap H(\rho)$ of $H(\rho)$, and in particular $T'(\rho)$ is
contained in
$K$.

Consider the principal bundle $\eta$ over $E$ with structure group
$H'(\rho)$ and transition  functions $\{\exp(f_{ij}X_s)\}$. This bundle
is topologically trivial (since $t\mapsto \{\exp(tf_{ij}X_s)\}$ 
defines a path from its transition functions to those of the trivial
bundle) and hence the complex structure on $\eta$ is defined by a flat
connection with values in the maximal compact subgroup $T'(\rho)$ of
$H'(\rho)$. Hence, with respect to the given open covering $\{U_i\}$, 
there  exist normalized transition functions $\{h'_{ij}\}$ with values
in
$T'(\rho)$ corresponding to a flat $H'(\rho)$-bundle $\eta_0$, where
$\eta_0$ is isomorphic as a holomorphic  $H'(\rho)$-bundle to
$\eta$. Let
$\varphi_i\colon U_i\to H'(\rho)$ be  local holomorphic functions such
that
$\varphi_i\exp(f_{ij}X_s) = h_{ij}' \varphi_j$. Consider now the bundle
$\hat\xi$ over $E$ whose transition functions relative to the open 
covering
$\{U_i\}$  are
$\{h_{ij}\exp(f_{ij}X_s)\}$. (Since $\operatorname{Ad}(h_{ij})X_s=X_s$ 
for all $i,j$, these functions satisfy the cocycle condition.) Changing
the local trivializations of $\hat\xi$ over the $U_i$ by the $\varphi_i$
changes the transition  functions for $\hat\xi$  to transition
functions $\{h_{ij}h'_{ij}\}$. Since the $h'_{ij}$ are normalized 
transition functions which commute with the $h_{k\ell}$, it is easy to
see that the functions $h_{ij}h'_{ij}$ are in fact normalized
transition functions for a Yang-Mills bundle. In fact, if
$\rho'\: \Gamma_\Ar \to K$ is the representation corresponding to the
$h'_{ij}$, then
$\rho$ and $\rho'$ commute, and thus $\rho\rho'$ is again a
homomorphism from $\Gamma _\Ar$ to $K$ whose normalized  transition
functions are $\{h_{ij}h'_{ij}\}$. In particular, the bundle $\hat \xi$
is a Yang-Mills bundle. 

Since by definition the holomorphic functions
$\varphi_i$ commute with $X_n$, changing the local trivializations of
$\xi$ by the
$\varphi_i$ produces new transition functions 
$$\{h_{ij}h'_{ij}\exp(f_{ij}X_n)\}$$ for $\xi$. Here $X_n$ is nilpotent
and is centralized by the
$h_{ij}$ and the $h'_{ij}$ and hence by the products
$h_{ij}h'_{ij}$.  

A standard argument shows that, since $X_n$ is nilpotent,  the bundle
defined by transition functions $\{h_{ij}h'_{ij}\exp (f_{ij}X_n)\}$ is
S-equivalent to that defined by transition functions
$\{h_{ij}h'_{ij}\}$, which  is $\hat\xi$. (It suffices to note that
there is a
$1$-parameter subgroup $\{\, h_\lambda: \lambda\in \Cee^*\,\}$ of $H$
such that
$\lim _{\lambda\to 0}\Ad(h_\lambda)(X) =0$.) Thus the bundle
$G$-bundle
$\xi_0$ with normalized transition functions $\{h_{ij}\}$ is
S-equivalent to the $G$-bundle $\hat\xi$ with normalized transition
functions $\{h_{ij}h'_{ij}\}$.   By Theorem 3.1, this is only possible
if the two bundles are isomorphic as Yang-Mills 
$K$-bundles. By (ii) of Lemma 3.3, there is an element $g\in K$ such
that
$h_{ij}=g(h_{ij}h'_{ij})g^{-1}$ for all $i,j$. Let $X =
\operatorname{Ad} (g)X_n$. After conjugating the local trivializations 
by
$g$, the bundle
 $\xi$ has transition functions of the form
$h_{ij}\exp (f_{ij}X)$, where $X$ is a nilpotent element of $\frak g$.
Since $X_n$  is centralized by the $h_{ij}h'_{ij}$, the element $X$  is
centralized by $h_{ij}$. By  Lemma 3.3, this implies that
$X \in \frak z_{\frak g}(\rho)$. Of course,
$X=\operatorname{Ad}(g)X_n$ is a nilpotent  element of $\frak g$, and
hence of $\frak z_{\frak g}(\rho)$ We have thus found transition
functions for $\xi$ as claimed in Theorem 3.6.
\endproof

We shall also refer to transition functions of the form  $\{h_{ij}\exp
(f_{ij}X)\}$,  where the $\{h_{ij}\}$ are normalized transition
functions for the S-equivalent Yang-Mills bundle with holonomy $\rho$
and
$X\in\frak z_{ \frak g}(\rho)$ is nilpotent, as {\sl normalized
transition  functions\/} for the semistable bundle $\xi$.

\corollary{3.7} Suppose that $K$ is semisimple. Let $\xi$ be a
semistable principal
$G$-bundle. Then the structure  group of $\xi$ reduces to an abelian
subgroup of $G$. More precisely, let
$\xi_0$ be the bundle with flat $K$-connection S-equivalent to $\xi$
and let $\rho\colon \pi_1(E)\to K$ be the holonomy representation of
this connection. Then there exists  a $1$-parameter unipotent subgroup
$U$ of $Z_G(\rho)$ such that the  structure group of $\xi$ reduces to
$\operatorname{Im}(\rho) \times U$.
\endstatement

\proof In the notation of the statement of Theorem 3.7, let $U$ be the
$1$-parameter unipotent subgroup of $Z_G(\rho)$ defined by $\exp
(tX_n)$. Then we  have defined a holomorphic reduction of structure to
the abelian subgroup
$\operatorname{Im}(\rho)\times U$ of $G$.
\endproof

\remark{Remark} Note however that we cannot necessarily reduce the
structure group to a {\sl connected\/} abelian subgroup, even when $G$
is simply connected.
\endremark
\medskip

\ssection{3.4. Defining holomorphic structures via $(0,1)$-connections.}

We can also give a description of the holomorphic structure on $\xi$ in
terms of a
$(0,1)$-connection. Fix once and for all a nonzero holomorphic  $1$-form
$dz$ on
$E$, corresponding to a choice of coordinate $z$ on the universal cover 
$\Cee$, and let $d\bar z$ be the corresponding $(0,1)$-form.  We fix a
good open cover
$\{U_i\}$ of
$E$ and normalized transition functions $\{h_{ij}\}$ for the
S-equivalent Yang-Mills $K$-bundle. Let $X\in
\frak z_{\frak g}(\rho)$. Over each
$U_i$, we can define a $(0,1)$-connection $\dbar _i$ by the formula
$$\dbar _i = \dbar + Xd\bar z,$$ where $\dbar$ denotes the trivial
connection on the product bundle over 
$U_i$. Since  $X$ is centralized by the $h_{ij}$, we can glue these
connections together via the transition functions $\{h_{ij}\}$ to
obtain a well-defined
$(0,1)$-connection 
$\dbar _{\xi _0,X}$ on $\xi _0$. Our result is then:

\proposition{3.8} In the above notation, the holomorphic structure
defined on $\xi_0$ by the $(0,1)$-connection $\dbar _{\xi _0,X}$ has
transition functions of the  form $\{h_{ij}\exp (f_{ij}X)\}$, where
$\{f_{ij}\}$ is a $1$-cocycle for $\scrO_E$  so that the induced
element of $H^1(E; \scrO_E)$ corresponds to
$-d\bar z$ under the Dolbeault isomorphism.
\endstatement

\proof We need to find local holomorphic cross sections $s_i$ over
$U_i$  and compare them over $U_i\cap U_j$. First, on $U_i$ there
exists a
$C^\infty$ function $\psi _i\:U_i \to \Cee$ such that $\dbar \psi _i =
d\bar z|U_i$. On
$U_i\cap U_j$, $\psi _i-\psi _j= -f_{ij}$ is holomorphic, and, via the
Dolbeault isomorphism,
$-f_{ij}$ is a generator for $H^1(E;\scrO_E)$, i\.e\. the
$1$-cocycle
$\{f_{ij}\}$ and $-d\bar z$ give the same element of $H^1(E;\scrO_E)$. 

Choose a faithful representation $G\to GL_n(\Cee)\subset \Bbb
M_n(\Cee)$,  the vector space of $n\times n$-matrices. Given a
$C^\infty$ section
$s_i\: U_i \to G$,  we identify $s_i$ with the corresponding function
to $\Bbb M_n(\Cee)$, and likewise we view the connection  $\dbar _{\xi
_0,X}$ as a connection on the corresponding vector bundle, defined by
the same formula. In this case, 
$s_i$ is holomorphic if and only if 
$$\frac{\partial s_i}{\partial \bar z} + X\cdot s_i = 0.$$ In
particular, we can choose $s_i$ to be of the form
$$s_i = \exp(-\psi _i(z)X).$$  Comparing on overlaps, we have
$$s_is_j^{-1} = h_{ij}\exp((-\psi _i + \psi _j)X) = h_{ij}\exp
(f_{ij}X),$$ as claimed.
\endproof

\remark{Remark} (1) If $\xi _0$ is a topologically trivial bundle over
$E$ with a flat
$K$-connection, then one can show that the holomorphic structure on
$\xi _0$ is given by a 
$(0,1)$-connection of the form $\dbar + hd\bar z$, where $h$ is a
constant element of a Cartan subalgebra $\frak h$ of $\frak g$.
However, in this description,
$Z_G(\rho)$ does not always act via constant elements of $G$.

(2) Ideally, it would be nice to have a version of (3.8) which said,
for example, that every holomorphic structure on the trivial bundle
arose  from  a connection on the trivial bundle with constant
coefficients, in other words of the form 
$\dbar + Yd\bar z$, where we could write $Y = h+X$ where $\ad h$ is 
semisimple,
$\ad X$ is nilpotent, and $[h, X] = 0$. While this is possible for 
$SL_n(\Cee)$ and $Sp(2n)$, one can show that, for every other simple and
simply connected group $G$, there exist  semistable holomorphic
structures on the $C^\infty$-trivial bundle which do not admit such a
description.
\endremark

\section{4. Isomorphism classification of semistable bundles.}

In the previous section, we gave a normal form for the transition
functions of a semistable bundle in terms of the transition functions
$\{h_{ij}\}$ of the associated Yang-Mills bundle and the choice of a
nilpotent element
$X\in \frak z_{\frak g}(\rho)$. Our goal in this section will be to
decide when two such normalized transition functions determine
isomorphic bundles. Related to this problem is that of determining the
group of holomorphic automorphisms of a semistable bundle. We give a
description in terms of $\{h_{ij}\}$ and $X$. At the end of this
section, we describe the Zariski tangent space of the deformations of a
semistable bundle in terms of the data $\rho$ and $X$.

\ssection{4.1. Isomorphisms and automorphisms of semistable bundles.}

Fix a good Stein cover $\{U_i\}$ of $E$ as in \S 3.2, and fix
holomorphic functions $f_{ij}$ on $U_i\cap U_j$ such that $\{f_{ij}\}$
is a
$1$-cocycle generating $H^1(E; \scrO_E)$. A semistable
$G$-bundle
$\xi$ is given by normalized transition functions
$\{h_{ij}\exp(f_{ij}X)\}$ and is thus determined by the representation
$\rho\: \Gamma_\Ar \to K$ and the nilpotent element $X \in \frak
z_{\frak g}(\rho)$. Thus, a pair $(\rho,X)$ as above determines a
semistable $G$-bundle. 

\theorem{4.1}
\roster
\item"{(i)}" Let $\xi$ be the semistable $G$-bundle determined by the
pair
$(\rho,X)$, where $\rho\: \Gamma_\Ar \to K$ is a representation and 
$X$ is a nilpotent element in $\frak z_{\frak g}(\rho)$. Then every
automorphism of
$\xi$ is given by a constant function $g$ with respect to the given
local trivializations of $\xi$ for the cover $\{U_i\}$, such that
$g\in Z_G(\rho)$ and $\operatorname{Ad}g(X) = X$. Thus, the group of
automorphisms of $\xi$ is identified with the subgroup of $Z_G(\rho)$
which centralizes $X$.
\item"{(ii)}" Two  pairs
$(\rho,X)$ and $(\rho',X')$ determine isomorphic holomorphic
$G$-bundles if and only if $\rho$ and $\rho'$ are conjugate in $G$ by
an element sending $X$ to $X'$. In this case, $\rho$ and $\rho'$ are
also conjugate by an element of $K$.
\endroster
\endstatement
\proof If $g\in G$ commutes with the image of $\rho$ and is such that
$\operatorname{Ad}g(X) = X$, then the automorphism of $\xi$ defined by
the constant function
$g$ commutes with the transition functions of $\xi$ and hence defines a
holomorphic automorphism of $\xi$. Likewise, if $g$ conjugates $\rho$
to $\rho'$ and $X$ to $X'$, then $g$ defines a holomorphic isomorphism
from $\xi$ to $\xi'$. Thus we need to prove the converse of these
statements. To prove (i), suppose that
$\varphi$ is a holomorphic automorphism of
$\xi$. We shall first show that $\varphi$ is given by a constant element
$g$ of $G$ and then argue that $g\in Z_G(\rho)$ and that
$\operatorname{Ad}g(X) = X$. The proof of (ii) will run along similar
lines. It will suffice to show that, for every irreducible
representation $\pi\: G \to GL(V)$, where $V$ is a finite-dimensional
vector space, $\pi(\varphi)$ is constant and commutes with the image of
$\pi\circ \rho$ and with $\pi_*(X)$. 

Thus, fix a finite-dimensional representation $\pi\:G \to GL(V)$. We
may assume that there is a Hermitian inner product on $V$ such that the
image of $K$ is contained in the corresponding unitary group $U(V)$. 
Let $\hat\rho\:\Gamma_\Ar\to U(V)$ be the representation induced from
$\rho$ by the homomorphism $K \to U(V)$. Then there is a decomposition
of $V$ as a $\Gamma_\Ar$-module: 
$$V=\bigoplus_\chi \chi\otimes V_\chi(\hat \rho),$$ where the $\chi$
are distinct irreducible representations of $\Gamma_\Ar$ and the
representation of $\Gamma_\Ar$ on the summand $\chi\otimes V_\chi(\hat
\rho)$ is the tensor product of the representation on
$\chi$ with  the trivial representation on $V_\chi(\hat \rho)$. We
shall always assume that
$V_\chi(\hat \rho)\neq 0$, i\.e\. that the $\chi$ in the direct sum
actually appear in $V$. Letting
$\xi_0$ be the Yang-Mills bundle determined by
$\rho$, there is then a decomposition of $V(\xi_0)=\xi_0\times_GV$ as a
direct sum of subbundles
$$V(\xi_0) \cong \bigoplus_\chi \Big(\lambda_\chi\otimes _\Cee
V_\chi(\hat
\rho_0)\Big),$$ where $\lambda_\chi$ is the vector bundle determined by
the representation $\chi$. Since $\chi$ is irreducible, $\lambda_\chi$
is stable. 

Next we must deal with the nilpotent element $X$.  Let
$Y\in \frak {gl}(V)$ be the image of
$X$ under $\pi$. It is a nilpotent transformation commuting with the
action of
$\Gamma_\Ar$. Thus, the action of $Y$ preserves the above decomposition
of $V$ and acts as $\Id \otimes Y_\chi$ on the summand $\chi \otimes
V_\chi(\hat\rho)$, where 
$Y_\chi$ is a nilpotent endomorphism of
$V_\chi(\hat\rho)$. Correspondingly, then, it is easy to check that $Y$
acts as a holomorphic automorphism of
$V(\xi_0)$, preserving the above decomposition of $V(\xi_0)$, and it
acts on each summand
$\lambda_\chi \otimes V_\chi(\hat\rho)$ as
$\Id \otimes Y_\chi$.

Quite generally, suppose that $V$ is a finite-dimensional vector space
and that
$Y$ is a nilpotent endomorphism of $V$. Define a holomorphic vector
bundle $V(Y)$ via the transition functions $\{\exp(f_{ij}Y)\}$. Note
that the transition functions for $V(Y)$ are in fact normalized
transition functions for the group
$SL(V)$, where $h_{ij} =1\in SU(V)$. In the notation of Lemma 1.1, if
$Y$ has $k$ Jordan blocks of size $d_i, i=1, \dots, k$, then 
$$V(Y) \cong \bigoplus _{i=1}^kI_{d_i}.$$ In particular, applying the
above construction to the  nilpotent endomorphism of
$V_\chi(\hat\rho)$ defined by the image
$Y$ of
$X$, we obtain a vector bundle $V_\chi(\hat\rho, Y)$ which is
isomorphic to a direct sum of bundles of the form $I_{d_i}$ for some
positive integers $d_i$. Clearly, as  a vector bundle,
$$V(\xi)=\bigoplus_\chi \Big(\lambda_\chi\otimes V_\chi(\hat\rho,
Y)\Big).$$

Now suppose that $\xi,\xi'$ are two holomorphic $G$-bundles determined
by two pairs $(\rho, X)$ and $(\rho', X')$ respectively. Let
$\varphi\:\xi\to \xi'$ be a holomorphic isomorphism. Let $Y$ and $Y'$
be the nilpotent elements of $\frak{gl}(V)$ determined by $X$ and $X'$,
respectively. There is an induced isomorphism
$$\hat\varphi\:V(\xi)\to V(\xi').$$

\claim{4.2} In the above notation, the isomorphism 
$$\hat\varphi \: V(\xi)=\bigoplus_\chi \Big(\lambda_\chi\otimes
V_\chi(\hat\rho, Y)\Big)\to V(\xi')=\bigoplus_\chi
\Big(\lambda_\chi\otimes V_\chi(\hat\rho', Y')\Big)$$  is compatible
with the direct sum decompositions. Moreover, the restriction of
$\hat\varphi$ to a homomorphism $\lambda_\chi\otimes V_\chi(\hat\rho,
Y)\to \lambda_\chi\otimes V_\chi(\hat\rho', Y')$ is of the form
$\Id\otimes M_\chi$, where 
$M_\chi\: V_\chi(\hat\rho)\to V_\chi(\hat\rho')$ is a  vector space
isomorphism such that $M_\chi Y = Y'M_\chi$.
\endstatement

\proof The vector bundle $V_\chi(\hat\rho, Y)$ is filtered by
subbundles such that each successive quotient is isomorphic to
$\scrO_E$. Thus, 
$\lambda_\chi\otimes V_\chi(\hat\rho, Y)$ is filtered by subbundles such
that each succesive quotient is isomorphic to $\lambda_\chi$. Hence,
the slope of $\lambda_\chi\otimes V_\chi(\hat\rho, Y)$ is equal to the
slope of
$\lambda_\chi$. It follows by the semistability of $V(\xi)$ that all of
the vector bundles $\lambda_\chi$ which appear in the above direct sum
have the same slope, and similarly for the bundles appearing in the
direct sum decomposition for $V(\xi')$. Since $V(\xi)\cong V(\xi')$,
all of these common slopes are equal. But if $\lambda_\chi$ and
$\lambda_{\chi'}$ are two such stable bundles with the same slope,
every holomorphic map from
$\lambda_\chi$ to
$\lambda_{\chi'}$ is  zero if $\chi$ and $\chi'$ are distinct
irreducible representations of
$\Gamma_\Ar$, and it is multiplication by a scalar if $\chi = \chi'$. In
particular, we see that $\hat \varphi$ preserves the direct sum
decomposition.

Since, by Lemma 1.15, the map $M\mapsto \Id \otimes M$ is an isomorphism
$$\Hom (V_\chi(\hat\rho, Y) ,V_\chi(\hat\rho', Y')) \cong  \Hom
(\lambda_\chi\otimes  V_\chi(\hat\rho, Y), \lambda_\chi\otimes
V_\chi(\hat\rho', Y')),$$  we see that Claim 4.2 follows from: 
\enddemo

\lemma{4.3} For $i=1,2$, let $V_i$ be a finite-dimensional vector space
and let
$Y_i\: V_i \to V_i$ be a nilpotent endomorphism. Let $\varphi\: V_1(Y_1)
\to V_2(Y_2)$ be a holomorphic map. Then, in the given local
trivializations of the bundles $V_i(Y_i)$,
$\varphi$ is given by a constant element $M$ of $\Hom (V_1, V_2)$ such
that $MY_1 = Y_2M$.
\endstatement
\proof We begin with the following special case:

\lemma{4.4} Let $V$ be a finite-dimensional vector space with a
nilpotent endomorphism $Y$, and let $s$ be a holomorphic section of
$V(Y)$. Then, in the given local trivializations of the bundle $V(Y)$,
$s$ is given by a constant element $s$ of $V$ such that $Y(s) = 0$.
\endstatement
\proof It is clear from the description of the transition functions
that every element $s\in V$ such that $Ys = 0$ defines a holomorphic
section of $V(Y)$. Conversely, let $s$ be a holomorphic section of
$V(Y)$. We must show that $s$ is a constant section in $\Ker Y$. As we
have seen in
\S 3.4, the holomorphic structure on $V(Y)$ is defined by the
$(0,1)$-connection $\dbar + Yd\bar z$ on the $C^\infty$ trivial bundle
$E\times V$. Putting
$Y$ into Jordan blocks decomposes $V(Y)$ as a direct sum, and it
suffices to work on  each summand. On each such summand of dimension
$k$, there is a basis $e_1,
\dots, e_k$ such that
$Ye_i = e_{i-1}$ if $i>1$, and $Ye_1 = 0$. A holomorphic section of
$V(Y)$ corresponds to a vector of $C^\infty$-functions on $E$,
$(s_1,
\dots, s_k)$ such that 
$$ \frac{\partial s_i}{\partial \bar z} + s_{i+1} = 0, \quad\text{if $i<
k$,}\qquad
\text{and}\qquad
\frac{\partial s_k}{\partial \bar z} = 0.$$ First it follows that $s_k$
is a holomorphic function on $E$ and is thus a constant $c_k$. Next, we
can view  $s_{k-1}$  as a doubly periodic
$C^\infty$ function on the universal cover $\Cee$ of $E$ of the form
$-c_k\bar z + f_{k-1}(z)$, where
$f_{k-1}$ is holomorphic on $\Cee$. In this case $\partial s_{k-1} =
df_{k-1}$ is a holomorphic differential on $E$, and thus $f_{k-1}(z) =
b_{k-1}z + c_{k-1}$, where $b_{k-1}$ and $c_{k-1}$ are constant. An easy
argument shows that a function of the form $-c_k\bar z + b_{k-1}z +
c_{k-1}$ is doubly periodic if and only if $b_{k-1} = c_k = 0$. Thus
$s_k = 0$ and $s_{k-1}$ is a constant $c_{k-1}$. Continuing in this
way, we see by induction that $s_i= 0, i>1$ and that $s_1$ is a
constant. Thus the space of holomorphic sections of $V(Y)$ is exactly
the set of constant sections with values in $\Ker Y$. This concludes
the proof of Lemma 4.4.
\endproof

To prove Lemma 4.3, first note that it follows by comparing transition
functions that the vector bundle
$Hom (V_1(Y_1), V_2(Y_2))$ is given by $V(Y)$, where $V = \Hom (V_1,
V_2)$, and
$Y$ is the nilpotent endomorphism of $V$ defined by $Y(M) = Y_2M -
MY_1$. Thus, by Lemma 4.4, the space of holomorphic maps from
$V_1(Y_1)$ to $V_2(Y_2)$, or equivalently the space of global
holomorphic sections of $Hom (V_1(Y_1), V_2(Y_2))$, is the same as the
space  of
$M\in\Hom (V_1, V_2)$ such that $MY_1 = Y_2M$.
\endproof

Returning to the proof of Theorem 4.1,  suppose that $\varphi\: \xi \to
\xi'$ is an isomorphism of $G$-bundles. Then $\varphi$ induces an
isomorphism
$\hat\varphi\: V(\xi)
\to V(\xi')$ of the associated vector bundles, which by Claim 4.2 is
given by a constant linear isomorphism $M$ from $V$ to itself (with
respect to the given local trivializations) such that $MYM^{-1} = Y'$.
Since this is true for all irreducible representations $\pi$, we must
have $\varphi = g$ for some constant element $g\in G$, and
$\operatorname{Ad}g(X) = X'$. In particular,
$g\exp(f_{ij}X) = \exp(f_{ij}X')g$. But $g$ satisfies
$$gh_{ij}\exp(f_{ij}X)= h_{ij}'\exp(f_{ij}X')g$$ for all $i,j$, and
thus $gh_{ij} = h_{ij}'g$. Hence $g$ conjugates $\rho$ to
$\rho'$. In particular, if $\xi =\xi'$, then $g$ commutes with the
$h_{ij}$ and hence with $\operatorname{Im} \rho$, and furthermore
$\operatorname{Ad}g(X) = X$. This proves (i), as well as the first
statement of (ii). 

To see the last statement of (ii), it will suffice to prove that two
representations
$\rho$ and $\rho'$ of $\Gamma_\Ar$ into
$K$ are conjugate by an element of  $G$ if and only if they are
conjugate by an element of $K$. Following the argument of Ramanathan
\cite{27}, this is an easy consequence of the Cartan decomposition:
there exists a subset $Q\subseteq G$, invariant under conjugation by
$K$, such that every element $g$ of $G$ can be written uniquely as
$g=kq$, where $k\in K$ and $q\in Q$. Using this, it is a
straightforward argument to show that, if
$\rho$ and $\rho'$ are conjugate via $g=kq$, then they are already
conjugate via
$k$. This completes the proof of (ii). 
\endproof

We are now in a position to give the analogue for semistable
$G$-bundles of Atiyah's classification of vector bundles over an
elliptic curve.

\theorem{4.5} Let $\xi$ be a holomorphic principal $G$-bundle. Then
there exists a triple $(L,  \rho_L, X)$, where 
\roster
\item"{(i)}" $L$ is a reductive subgroup of $G$ with a maximal compact
subgroup $K_L$;
\item"{(ii)}" $\rho_L\: \Gamma_\Ar\to K_L$ is a central homomorphism and
$X$ is a nilpotent element of $\frak z_{\frak l}(\rho_L)$, where $\frak
l$ is the Lie algebra of $L$;
\item"{(iii)}" $L$ is saturated with respect to $\mu = d\rho_L(1)$in
the sense of
\rom{(2.12)};
\endroster such that, if $\xi_L$ is the semistable $L$-bundle
associated to the pair $(\rho_L, X)$, then  $\xi \cong \xi_L\times_LG$. 
Two triples $(L, \rho_L, X)$ and $(L', \rho_{L'}, X')$ with the above
properties lead to isomorphic $G$-bundles if and only if there is an
element
$g\in G$ which conjugates $L$ to $L'$, $\rho_L$ to $\rho_{L'}$
\rom(viewing these as homomorphisms to $L$ and $L'$ respectively\rom),
and $X$ to $X'$.
\qed
\endstatement

\ssection{4.2. Deformations of semistable bundles.}

Let $\xi$ be a semistable bundle. Methods very similar to the proof of
Theorem 4.1 show the following:

\theorem{4.6} Let $\xi$ be a semistable bundle corresponding to the
representation 
$\rho$ and the nilpotent element $X\in \frak z_{\frak g}(\rho)$. Then 
$$\align H^0(E; \ad \xi) & \cong \Ker\{\,\ad X: \frak z_{\frak g}(\rho)
\to \frak z_{\frak g}(\rho)\,\};\\ H^1(E; \ad \xi) & \cong
\operatorname{Coker}\{\,\ad X: \frak z_{\frak g}(\rho) \to
\frak z_{\frak g}(\rho)\,\}.\qed
\endalign$$
\endstatement

Of course, the first isomorphism is just a less precise version of
Theorem 4.1(i). The second isomorphism, in case $X=0$, was given in
Lemma 3.4(ii). As in the remarks before the statement of Lemma 3.4, if
we begin with normalized transition functions
$\{h_{ij}\exp(f_{ij}X)\}$ for $\xi$, given a $Y\in \frak z_{\frak
g}(\rho)$, we can form a new bundle with (not necessarily normalized)
transition functions
$\{h_{ij}\exp(f_{ij}(X+Y))\}$. By the proof of Lemma 3.4(ii), the
Kodaira-Spencer class of this family in the tangent direction
corresponding to $Y$ is
$\{f_{ij}Y\}$, viewed as a section of $H^1(E; \ad \xi)$. The
Kodaira-Spencer map then defines a homomorphism
$\frak z_{\frak g}(\rho) \to H^1(E; \ad \xi)$, and under the
identification of Theorem 4.6 this becomes the obvious surjection
$\frak z_{\frak g}(\rho) \to \frak z_{\frak
g}(\rho)/\operatorname{Im}\ad X$. In particular, a complement to the
image of $\ad X$ in $\frak z_{\frak g}(\rho)$ defines a locally
semiuniversal deformation of $\xi$ by the above construction.

\section{5. Local and global properties of the moduli space.}

\ssection{5.1. The moduli space as a complex variety.}

We begin by recalling the general properties of the moduli space of
semistable
$G$-bundles in the case of a  smooth curve $C$  of arbitrary genus. We
have the following theorem
\cite{28}:

\theorem{5.1} For each reductive group $G$ and each smooth projective
curve $C$, there exists a coarse moduli space $\Cal M_C(G)$ of
S-equivalence classes of semistable principal $G$-bundles over $C$. If
$G$ is simply connected, $\Cal M_C(G)$ is irreducible. In general, the
connected components of $\Cal M_C(G)$ are in one-to-one correspondence
with the topological types of $G$-bundles over $C$, in other words,
with $\pi _1(G)$. The components of $\Cal M_C(G)$ are normal projective
varieties.
\qed
\endstatement

Given an element $c\in \pi_1(G)$, we let $\Cal M_C(G,c)$ be the
corresponding component of the moduli space $\Cal M_C(G)$. Of course,
if $G$ is simply connected, then necessarily $\Cal M_C(G,c) = \Cal
M_C(G)$.

Strictly speaking, Theorem 5.1 is proved in \cite{28} in case $C$ has
genus at least $2$, but the proof applies in the special cases $g=0,1$
with minor modifications. The analogue of Theorem 5.1 for the open
subset of stable bundles is proved in
\cite{27}. Faltings
\cite{9} constructed a compactification of this space where the
boundary points were S-equivalence classes of strictly semistable
bundles. This result also does not apply in  case $C=E$ has genus one
and $G$ is simply connected, since there are no properly stable
$G$-bundles over an elliptic curve. However, the method of \cite{9}
suffices to handle this case as well. See also LePotier \cite{23}. 

  From now on, we restrict to the case $C=E$, and abbreviate $\Cal
M_E(G)$ by $\Cal M(G)$, and similarly for $\Cal M_E(G,c)$.

\ssection{5.2. The local structure of the moduli space.}

\theorem{5.2} Let $\xi_0$ be a Yang-Mills bundle, corresponding to the
representation
$\rho$. Let $[\xi_0]$ be the corresponding S-equivalence class, viewed
as a point of the coarse moduli space $\Cal M(G)$. Then there is a
neighborhood of $[\xi_0]$ in
$\Cal M(G)$ which is biholomorphic to a neighborhood of the origin in
the GIT quotient
$\frak z_{\frak g}(\rho)/\!\!/Z_G(\rho)$, where $Z_G(\rho)$ acts on
$\frak z_{\frak g}(\rho)$ via the adjoint action. Finally, the
universal bundle $\Xi_0$ over $E\times\frak z_{\frak g}(\rho)$
constructed in Lemma \rom{3.4} defines a $Z_G(\rho)$-invariant map 
from a neighborhood of the origin in
$\frak z_{\frak g}(\rho)/\!\!/Z_G(\rho)$ to $\Cal M(G)$ which induces
the above isomorphism.
\endstatement
\proof The proof will, for the most part, consist in recalling various
standard results. First, the construction of \cite{28} realizes $\Cal
M(G)$ as the GIT quotient of a smooth quasiprojective 
variety $R$ by a reductive group $G= GL_N(\Cee)$ for some $N$.
Furthermore, if $x\in R$ is some point lying over $[\xi_0]
\in \Cal M(G) = R/\!\!/G$, then the normal space $N$ to the orbit
$G\cdot x$ in $R$ is identified with the Zariski tangent space to the
deformation space of $\xi_0$, and this identification is equivariant
with respect to the action of the isotropy group
$G_x = \Aut \xi_0$. Now the assumption that $\xi_0$ is a Yang-Mills
bundle means, by (ii) of Theorem 3.1,  that
$G\cdot x$ is a closed orbit. A standard application of Luna's \'etale
slice theorem implies that there is an \'etale neighborhood of
$[\xi_0]$ in
$\Cal M(G)$ which is isomorphic to an \'etale neighborhood of the
origin in the GIT quotient $H^1(E; \ad \xi_0)/\!\!/\Aut \xi_0$. On the
other hand, it follows from  Lemma 3.4 that $H^1(E;
\ad \xi_0) \cong \frak z_{\frak g}(\rho)$, and it is easy to see that
this isomorphism is equivariant with respect to the identification of
$\Aut \xi_0$ with
$Z_G(\rho)$ given in Theorem 4.1 (for the case $X= 0$). Thus we have
identified $\Cal M(G)$, locally at
$[\xi_0]$, with $\frak z_{\frak g}(\rho)/\!\!/Z_G(\rho)$. The final
statement follows by Lemma 3.4(ii).
\endproof

To further analyze the local structure of the moduli space, we
introduce some notation. Let $Z^0_G(\rho)$ be the identity component of
$Z_G(\rho)$, and let $I =
\pi_0(Z_G(\rho))$ be the (finite) component group of $Z_G(\rho)$. Let
$\frak h(\rho)$ be a Cartan subalgebra of 
$Z^0_G(\rho)$. The adjoint quotient $\frak z_{\frak
g}(\rho)/\!\!/Z^0_G(\rho)$ can be identified with the finite quotient
$\frak h(\rho)/W(\rho)$, which is an affine space. Furthermore, there
is an action of $I$ on $\frak h(\rho)$ for which the quotient $\frak
z_{\frak g}(\rho)/\!\!/Z_G(\rho)$ is identified with $(\frak
h(\rho)/W(\rho))/I$. Thus we see:

\corollary{5.3} In the above notation, a neighborhood of $[\xi_0] \in
\Cal M(G)$ is biholomorphic to a neighborhood of the origin in $(\frak
h(\rho)/W(\rho))/I$.
\qed
\endstatement

To be able to describe this action explicitly, we shall use the
following result about adjoint quotients:

\lemma{5.4} Let $\frak g$ be a reductive Lie algebra and let $I$ be a
finite group of outer automorphisms of $\frak g$. 
\roster
\item"{(i)}" There exists a realization of $I$ as a group of
automorphisms of $\frak g$, and a Cartan subalgebra
$\frak h$ invariant under $I$. Thus, if $W$ is the Weyl group of $\frak
g$ with respect to $\frak h$, then  the group
$W \rtimes I$ acts on
$\frak h$. 
\item"{(ii)}" The quotient
$\frak h/W$ is an affine space on which $I$ acts. Moreover, the action
of $I$ on
$\frak h/W$ may be linearized in such a way that the corresponding
representation of $I$ on $\frak h/W$ is isomorphic to the
representation of $I$ on $\frak h$.
\endroster
\endstatement
\proof The proof of (i) is standard; see for example \cite{7}, Chap\.
7. The first statement of (ii) is also standard; see for example
\cite{7}, Chap\. 5. To see the final statement, first consider the
algebra $\Cee[\frak h^*]^W$ of Weyl invariant polynomials on $\frak h$.
This is a graded $\Cee$-algebra, on which $I$ acts, and the action of
$I$ preserves the grading. Moreover, $\Cee[\frak h^*]^W =\Cee[p_1,
\dots, p_r]$ is a polynomial algebra, where the $p_i$ are $W$-invariant
polynomials of degree $d_i$. But now an elementary argument using the
complete reducibility of the action of $I$ on  $\Cee[\frak h^*]^W$
shows that there exist polynomials $q_1, \dots, q_r$, with $\deg q_i
=\deg p_i$, such that $I$ acts linearly on the span on the $q_i$ and
for which $\Cee[q_1, \dots, q_r] =\Cee[p_1, \dots, p_r]$. In
particular, this shows that the action of $I$ on $\frak h/W$ can be
linearized.

By the proof of Th\'eor\`eme 3.4 in \cite{5} (see also \cite{7}, Chap\.
9, Corollaire on p\. 48), there exists a regular element $x$ in
$\frak h^I$, the fixed set of the action of $I$ on $\frak h$.  By
definition, the differential of the quotient map
$\frak h
\to
\frak h/W$ is then an $I$-equivariant isomorphism at $x$. But the
differential of the action of $I$, acting on the tangent space of
$\frak h$ at $x$, gives a representation of $I$ which is isomorphic to
the representation of $I$ acting on
$\frak h$. A similar statement holds for the differential of the action
of $I$ on the tangent space to $\frak h/W$ at the image of $x$. Since
these actions are identified, we see that representation of $I$ on
$\frak h/W$ is isomorphic to the representation of $I$ on $\frak h$.
\endproof

Applying the lemma above to the case where $\frak g = \frak z_{\frak
g}(\rho)$ and using Corollary 5.3, we have:

\corollary{5.5} With notation as in Corollary \rom{5.3}, a neighborhood
of  $[\xi_0]
\in \Cal M(G)$ is biholomorphic to a neighborhood of the origin in
$\frak h(\rho)/I$.
\qed
\endstatement

\ssection{5.3. The component group and its action on a Cartan
subalgebra.}

Thanks to Corollary 5.5, to describe the local structure of the moduli
space at a point $[\xi_0]$, it will suffice to understand the component
group $I =\pi_0(Z_G(\rho))$ and its action on a Cartan subalgebra. Here
we shall just recall the results of \cite{6}. We shall only consider
the case where $G$ is simple of rank
$r$. Let us introduce some basic notation. Suppose that $\alpha_1,
\dots, \alpha_r$ are the simple roots (with respect to some ordering)
of $G$, and let $\tilde \alpha$ be the highest root. Write the coroot
$\tilde \alpha\spcheck$ dual to $\tilde
\alpha$ as a positive linear combination of the simple coroots:
$$\tilde \alpha\spcheck = \sum _{i=1}^r g_i\alpha_i\spcheck.$$ Thus the
$g_i$ are positive integers. By convention we set $\alpha _0 = -\tilde
\alpha$ and $g_0 = 1$.

Next assume that $G$ is in addition simply connected. We then have the
following, which is proved in \cite{6}:

\theorem{5.6} Suppose that $G$ is simply connected. For a given
Yang-Mills bundle, let $\rho$ denote the corresponding representation.
\roster
\item"{(i)}" A Cartan subalgebra $\frak h$ of $\frak z_{\frak g}(\rho)$
is a Cartan subalgebra of $\frak g$. 
\item"{(ii)}" The group $I = \pi_0(Z_G(\rho))$ is cyclic, and its order
divides $g_i$ for some $i$.
\item"{(iii)}" The locus of points $[\xi_0]$ in $\Cal M(G)$ such that
$I$ is cyclic of order $d$ is a locally closed, irreducible analytic
subvariety of $\Cal M(G)$, and its complex dimension is equal to one
less than the number of $i$ such that
$d|g_i$.
\item"{(iv)}" If $I$ has order $d$, then the action of a generator of
$I$ on $\frak h$ has eigenvalues $\exp(2\pi \sqrt{-1}g_i/d)$ for each
$i$, $0\leq i\leq r$, such that
$d$ does not divide $g_i$. The multiplicity of the eigenvalue $1$ is
one less than the number of $i$ such that $d|g_i$. \qed
\endroster
\endstatement

Direct inspection of the numbers $g_i$ shows that, for $d> 1$, the
isotropy for the action of $I$ on $\frak h$ has codimension at least
two, and hence the locus where
$I$ has order $d$ for some $d>1$ is the singular locus of $\Cal M(G)$.

The corresponding results in case $G$ is not simply connected are
analogous, although a little more difficult to state. For the moment,
let us just assume that $G$ is semisimple. Let
$\tilde G$ be the universal cover of
$G$. Let $F\subseteq \Cal Z\tilde G$ be the kernel of the homomorphism
$\tilde G \to G$, where $\Cal Z\tilde G$  is the center of
$\tilde G$. Then from the exact sequence
$$0 \to F \to \tilde G \to G \to 0,$$ there is a first Chern class map
which associates to every $G$-bundle $\xi$ an element $c\in H^2(E; F)
\cong F$. Moreover, every $G$-bundle with first Chern class
$c$ lifts to a bundle over the group $\tilde G/\langle c\rangle$.
Finally, the moduli space of semistable $G$-bundles which lift to a
$\tilde G/\langle c\rangle$-bundle is a quotient of the moduli space of
semistable $\tilde G/\langle c\rangle$-bundles, which do not lift to
any further cover of $\tilde G/\langle c\rangle$, by the finite abelian
group $F/\langle c\rangle$. Thus, we shall assume that $G$ is of the
form $\tilde G/\langle c\rangle$, where 
$\tilde G$ is the universal cover of
$G$, and shall only consider bundles $\xi$ which do not lift to any
quotient of
$\tilde G$ by a proper subgroup of $\langle c\rangle$, or equivalently
such that the first Chern class of $\xi$ is a generator for $\langle
c\rangle$. The moduli space in the general case will then be a quotient
of such a moduli space by the action of the finite abelian group
$F/\langle c\rangle$.

In terms of representations, since $G$ is semisimple, a semistable
bundle $\xi$ is S-equivalent to a flat
$K$-bundle $\xi_0$, in other words to a homomorphism $\rho\: \pi_1(E)
\to K$, well-defined up to conjugation. To say that $\xi$, or
equivalently $\xi_0$, does not lift to a cover of $G$, in the case
where $G =\tilde G/\langle c\rangle$, says that the images of a pair of
generators of $\pi_1(E)$ lift to elements $x,y\in \tilde K$, the
universal cover of $K$, such that $[x,y] = c'$ for some generator $c'$
of
$\langle c\rangle$. The condition that the first Chern class of the
bundle $\xi_0$ is
$c$ easily implies that $c=c'$. Thus, we make the following definition:

\definition{Definition 5.7} Let $\tilde K$ be a compact, simply
connected group and let $c$ be an element of the center of $\tilde K$.
A {\sl $c$-pair\/} $(x,y) \in
\tilde K \times \tilde K$ is a pair $(x,y)$ such that $[x,y] = c$.
\enddefinition

We begin by showing that the classification of flat representations
$\rho$ and their centralizers in $K$ is essentially the same as the
classification of $c$-pairs $(x,y)$ and their centralizers in $\tilde
K$. 

\lemma{5.8} Let $\tilde K$ be a compact, simply connected group and let
$c$ be an element of the center of $\tilde K$. Let $K = \tilde
K/\langle c\rangle$. 
\roster
\item"{(i)}" Suppose that $\rho\: \pi_1(E) \to K$ and $\rho'\: \pi_1(E)
\to K$ are two homomorphisms, and that $(x,y)$ and $(x', y')$ are two
$c$-pairs in $\tilde K$ lifting the images of a pair of generators of
$\pi_1(E)$ under $\rho$, $\rho'$ respectively. Then $\rho$ and $\rho'$
are conjugate if and only if the pairs
$(x,y)$ and $(x', y')$ are conjugate by an element of $\tilde K$.
\item"{(ii)}" If $\rho\: \pi_1(E) \to K$ is a homomorphism, and $(x,y)$
is a
$c$-pair in $\tilde K$ lifting $\rho$ as above, then there is an exact
sequence
$$\{1\} \to \langle c\rangle \to Z_{\tilde K}(x,y) \to Z_K(\rho) \to
\langle c\rangle
\times \langle c\rangle \to \{1\}.$$ Thus there is an exact sequence
$$\langle c\rangle \to \pi_0(Z_{\tilde K}(x,y)) \to \pi_0(Z_K(\rho))
\to \langle c\rangle
\times \langle c\rangle \to \{1\}.$$ Similar statements hold when we
replace $K$ and $\tilde K$ by $G$ and $\tilde G$.
\item"{(iii)}" Suppose that $\tilde K\cong \prod _{i=1}^kK_i$ is a
product of simple and simply connected groups and that, under this
isomorphism, $c = (c_1, \dots, c_k)$. Then the moduli space of
$c$-pairs in $\tilde K$ modulo conjugation is isomorphic to the product
of the moduli spaces of of $c_i$-pairs in the groups
$K_i$, modulo conjugation in $K_i$.
\endroster
\endstatement
\proof To see (i), note that the statement that $\rho$ and $\rho'$ are
conjugate means that we can conjugate $(x,y)$ to $(c^ax', c^by')$ for
some choice of $a,b$. Since $(x', y')$ is a $c$-pair, it is easy to
check that the element
$(x')^{-b}(y')^a$ further conjugates $(c^ax', c^by')$ to $(x', y')$.
Conversely, if 
$(x,y)$ and $(x', y')$ are conjugate, then trivially $\rho$ and $\rho'$
are conjugate. 

To see (ii), let $\bar x$ and $\bar y$ be the images of $x,y$ in $K$.
An element
$\bar g$ of $K$ lies in $Z_K(\rho)$ if and only if $[\bar g, \bar x] =
1 = [\bar g,
\bar y]$. Lift $\bar g$ to an element $g$ of $\tilde K$. Then $[g,x] =
c^a$ and
$[g,y] = c^b$, where the elements $c^a, c^b$ are independent of the
choice of lift of $\bar g$. This gives a well-defined homomorphism
$\phi\: Z_K(\rho) \to  \langle c\rangle \times \langle c\rangle$. An
argument using powers $\bar x^{-b}\bar y^a$ as in the proof of (i)
shows that $\phi$ is surjective. Its kernel is the set of elements
$\bar g$ which lift to an element of $Z_{\tilde K}(x,y)$, yielding the
first exact sequence in (ii). From this, the exactness of the second
follows. Finally, (iii) is clear.
\endproof

In the second exact sequence in (ii), the element $c$ may or may not be
in the identity component of $Z_{\tilde K}(x,y)$, and hence the map
$\pi_0(Z_{\tilde K}(x,y)) \to \pi_0(Z_K(\rho))$ may or may not be
injective. However, one can show that the answer only depends on
$\tilde K$ and $c$, not on $x$ and $y$. In fact, in the notation
introduced below, the answer only depends on whether $c$ lies in the
identity component of $T^{w_c}$.

We return to the problem of describing the GIT quotient $\frak z_{\frak
g}(\rho)/\!\!/Z_G(\rho)$, or equivalently
$(\frak h(\rho)/W(\rho))/I$, in the non-simply connected case, and
shall assume that
$G$ is in fact simple. Note that the center of $\tilde G$ acts
trivially on $\frak z_{\frak g}(\rho)$ under the adjoint
representation, as do the elements $x$ and $y$. Hence 
$$\frak z_{\frak g}(\rho)/\!\!/Z_G(\rho) = \frak z_{\frak
g}(\rho)/\!\!/Z_{\tilde G}(x,y).$$ To analyze the action of $I$, we
introduce some more notation.  Fix a Cartan subalgebra $\frak h$ of
$\frak g$ and and a Weyl chamber in $\frak h$. Let $\alpha_1,
\dots, \alpha_r$ be the corresponding set of simple roots. There is a
well-defined homomorphism from the center of
$\tilde K$, or equivalently the center of $\tilde G$, to the Weyl group
$W$
\cite{6}. We denote the image of $c$ under this homomorphism by $w_c$.
If $\tilde K$ is simple, the Weyl element $w_c$ permutes the set $\tilde
\Delta =\{\alpha_0, \alpha_1, \dots, \alpha_r\}$, where
as usual
$\alpha _0 = -\tilde\alpha$.  Let
$T^{w_c}$ be the fixed subgroup of
$T$ under the action of
$w_c$. Thus
$T^{w_c}$ is an extension of a finite group by a torus. We define
$H^{w_c}$ and
$\frak h^{w_c}$ similarly. 

The action of $w_c$ on $\tilde \Delta$ divides $\tilde \Delta$ into
orbits $\bold o$. The number $g_i$ is independent of the choice of
$\alpha _i\in \bold o$, and depends only on $\bold o$. We write this
number as $g_{\bold o}$. For each orbit $\bold o$, let $n_{\bold o}$ be
the number of elements of $\bold o$. Finally, let $d_0$ be the gcd of
the numbers $n_{\bold o}g_{\bold o}$. We then have the following
generalization of Theorem 5.6, which is proved in
\cite{6}:

\theorem{5.9} Suppose that $G$ is simple. For a given Yang-Mills bundle,
let $\rho$ denote the corresponding representation.
\roster
\item"{(i)}" There exists a Cartan subalgebra for $\frak z_{\frak
g}(\rho)$ of the form $\frak h^{w_c}$, where $\frak h$ is a Cartan
subalgebra of $\frak g$.
\item"{(ii)}" The group $\tilde I = \pi_0(Z_{\tilde G}(x,y))$ is
cyclic, and its order divides
$n_{\bold o}g_{\bold o}$ for some $\bold o$. Moreover, $d_0$ divides
the order of
$\tilde I$ for every $\rho$.
\item"{(iii)}" Suppose that $d_0|d$ and that $d$ divides
$n_{\bold o}g_{\bold o}$ for some $\bold o$. The locus of points
$[\xi_0]$ in
$\Cal M(G)$ such that $\tilde I$ is cyclic of order $d$ is a locally
closed, irreducible analytic subvariety of
$\Cal M(G)$, and its complex dimension is equal to one less than the
number of $\bold o$ such that $d|n_{\bold o}g_{\bold o}$.
\item"{(iv)}" If $\tilde I$ has order $d$, then the action of a
generator of $\tilde I$ on
$\frak h^{w_c}$ has eigenvalues $\exp(2\pi \sqrt{-1}(n_{\bold
o}g_{\bold o})/d)$ for each orbit
$\bold o$ such that
$d$ does not divide $n_{\bold o}g_{\bold o}$. The multiplicity of the
eigenvalue $1$ is one less than the number of $\bold o$ such that
$d|n_{\bold o}g_{\bold o}$. \qed
\endroster
\endstatement

Unlike the simply connected case, it is possible for $\tilde I$ to be
nontrivial for every $\rho$. This happens exactly when $d_0>1$.
Moreover, the action of $\tilde I$ on $\frak h^{w_c}$ is not faithful
in this case. Another difference between the simply connected and the
non-simply connected case is that it is possible for there to be
isotropy in codimension one. Thus the locus where $\tilde I$ has more
than $d_0$ elements need not be equal to the singular locus of $\Cal
M(G,c)$.

\ssection{5.4. The global structure of the moduli space for a
semisimple group.}

We begin with the case where $G$ is semisimple and simply connected. Of
course, the same will be true for the maximal compact subgroup $K$. In
this case, the moduli space is described as follows:

\theorem{5.10} Let $K$ be a compact, simply connected group and let
$T\subset K$ be a maximal torus. Let $W=W(T,K)$ be the Weyl group of
the pair $(K,T)$. Then the moduli space of isomorphism classes of
representations $\rho\:\pi_1(E) \to K$ is isomorphic to 
$(T\times T)/W$ where the action of $W$ is the diagonal action.
\endstatement

\proof Choosing an isomorphism $\pi_1(E) \cong \Zee\oplus \Zee$, let
$x,y\in K$ be the images under $\rho$ of the two generators. Clearly,
$\rho$ is specified by $x$ and $y$ and the isomorphism class of $\rho$
is determined by the pair $(x,y)$ modulo simultaneous conjugation by
$K$. Furthermore $x$ and $y$ commute. We then have the following, which
is  a theorem of Borel \cite{5} or \cite{9, p. 48, Cor. 1(b)}: 

\lemma{5.11} Let
$K$ be a compact,  simply connected Lie group, and let
$x$ and $y$ be two commuting elements of $K$. Then there exists a
maximal torus $T$ in $K$ with $x,y\in T$.
\endstatement
\proof Given $y\in K$, let $T_1$ be a maximal torus in $K$ containing 
$y$ and let $K(y)$ be the centralizer of $y$ in
$K$. Then by \cite{5}  or \cite{9, p. 48, Th\'eor\`eme 1 (c)},
$K(y)$ is connected.  Since $y\in T_1$,
$T_1\subseteq K(y)$ is clearly a maximal torus of $K(y)$. It follows
that  every maximal torus for
$K(y)$ is also a maximal torus for $K$. Since
$x\in K(y)$, there is a maximal torus
$T$ of $K(y)$ containing $x$. Since $y$ is in the center of 
$K(y)$, $y$ is contained in every maximal torus of $K(y)$ and in 
particular in
$T$. Thus $T$ is a maximal torus of $K$ which contains both $x$ and $y$.
\endproof

To complete the proof of Theorem 5.10, after conjugation we may assume
that $x$ and
$y$ lie in a fixed maximal torus $T$ of $K$, and must consider the
equivalence relation induced by conjugation by  elements of $K$ on
$T\times T$. But an elementary argument
\cite{7, Chap\. 9, Corollaire 7 on p\. 10} shows that, if two ordered
pairs
$(x,y)$ and
$(x',y')$ in  $T\times T$  are conjugate by an element  of $K$, then
they are conjugate by an element of the Weyl group $W$ of $K$ with
respect to
$T$. Thus the moduli space of all commuting pairs up to conjugation is
identified with $(T\times T)/W$.
\endproof

There is another way to formulate Theorem 5.10.  Writing
$T$ invariantly as $U(1)\otimes_\Zee \Lambda$, where $\Lambda$ is the
coroot lattice of $G$, we have 
$$\Hom(\pi_1(E),T) = \Hom (\pi_1(E), U(1)) \otimes_\Zee \Lambda.$$ But
$\Hom (\pi_1(E), U(1))$ is the space of flat connections on $E$ and so
is canonically $\Pic^0E$, which we can identify with $E$ by the choice
of the base point $p_0$. Under this identification, given $\rho\:
\pi_1(E) \to T$, a character
$\alpha\colon T\to U(1)$ defines a homomorphism $\chi\: \pi_1(E)\to
U(1)$, namely
$\alpha \circ \rho$.  The flat $U(1)$-bundle corresponding to $\chi$
then determines a complex line bundle $\lambda_\chi$ of degree zero,
which is an element of $\Pic^0(E)\cong E$. Thus, given $\rho$, each
$\alpha\colon T\to U(1)$ determines a point $e_{\alpha(\rho)}\in E$. Of
course, $\alpha$ also induces a map $\Lambda\to \Zee$. Clearly the
element of
$E$ obtained by identifying $\rho$ with a point of $E\otimes _\Zee
\Lambda$ and evaluating the map $E\otimes _\Zee \Lambda \to E$ defined
by $\alpha$ at the point $\rho$ is just $e_{\alpha(\rho)}$. Note that
$\alpha (\rho)=1$ is the trivial character, if and only if
$\lambda_{\alpha(\rho)}=
\scrO_E$, if and only if $e_{\alpha(\rho)} = 0$.

Clearly, the natural action of the Weyl group on
$\Hom(\pi_1(E),T)$ goes over the natural action of $W$ on
$\Lambda$ in the tensor product $E\otimes _\Zee\Lambda$. Thus, we see:

\corollary{5.12} Let $E$ be an elliptic curve and let $K$ be a compact
simply connected group with complexification $G$. Let $\Lambda$ denote
the coroot lattice for $K$, i.e., the  fundamental group of a maximal
torus of $K$. Then  the moduli space of S-equivalence classes of
semistable
$G$-bundles is identified with the moduli space of flat
$K$-connections on $E$ which in turn is identified with
$(E\otimes _\Zee\Lambda)/W$, where the action of $W$ on $E\otimes
_\Zee\Lambda$ is via the natural action of $W$ on $\Lambda$.  Thus it
is a normal complex projective variety of complex dimension equal to
the rank of $G$.
\qed
\endstatement

A different proof of Corollary 5.12 has been given by Laszlo \cite{22}.
We shall compare the complex structure of $(E\otimes _\Zee\Lambda)/W$
with that of $\Cal M(G)$ in \S 5.6 below.

We turn now to the non-simply connected case, i\.e\.  to the
description of the moduli space of flat $K$-bundles corresponding to
first Chern class $c$. By (i) of Lemma 5.8, this is the same as the
moduli space of $c$-pairs $(x,y)$ up to conjugation. In order to
describe this moduli space, we first introduce some notation. Fix a set
of simple roots $\alpha_1, \dots, \alpha_r$ for
$\tilde K$. Then there is the element $w_c$ of $W$ defined above.  As
before, let
$T^{w_c}$ be the fixed set of $w_c$ acting on $T$. Let
$T_{w_c}$ be the group of coinvariants of $T$ under the action of
$w_c$: $T_{w_c} = T/(\Id - w_c)T$. Thus
$T_{w_c}$ is a torus, and in particular it is connected. 

Fix an element $x_0 \in T$ and a representative $v$ for $w_c$ such that
$[x_0, v] = c$; it is not difficult, using the construction of $w_c$,
to see that such $x_0$ exist. Then there is a map from $T^{w_c} \times
T$ to the space of $c$-pairs defined by
$$(t,s) \mapsto (x,y) = (tx_0, sv).$$ An easy calculation shows that
the assumptions $s\in T^{w_c}$, $t\in T$ imply that the $(x,y)$ above
is a $c$-pair. Here the projection of the image to the first factor is
clearly the set of $r\in T$ such that $w_c(r) = r+c$ (writing the group
law on $T$ additively). We write this set as $T^c$. Conjugation of
$(x,y)$ by
$u\in T$ corresponds to replacing $(t,s)$ by $(t, s+u-w_c(u))$, if we
write the group law on $T$ additively. Thus the map from $T^{w_c}
\times T$ to the space of $c$-pairs defined above factors through the
quotient map $T^{w_c} \times T\to T^{w_c} \times T_{w_c}$.  

We then have the following, whose proof is due to Schweigert \cite{29}
and is also given in
\cite{6}:

\theorem{5.13} Let $Z_W(w_c)$ be the centralizer of $w_c$ in the Weyl
group $W$. Then the moduli space of $c$-pairs in $\tilde K$ modulo
conjugation is equal to
$(T^c\times T_{w_c})/Z_W(w_c)$. Here the action of $Z_W(w_c)$ on the
first factor
$T^c$ is by an affine action, and it is the natural linear action on
the second factor.
\qed
\endstatement

\corollary{5.14} The moduli space of $c$-pairs has real dimension equal
to twice the dimension of $T^{w_c}$, or equivalently twice the
dimension of the vector space
$\frak t^{w_c}$.
\qed
\endstatement

Since $T^c$ need not be connected, it is not at all apparent from
Theorem 5.13 that the moduli space is connected. In fact, there are
other descriptions of the moduli space. One can show that the quotient
$(T^c
\times T_{w_{c}})/ Z_W(w_{c})$ is identified with the quotient
$(T_{w_{c}}\times T_{w_{c}})/ Z_W(w_{c})$, where the group $Z_W(w_{c})$
does not act faithfully, but acts via  group homomorphisms on
$T_{w_{c}}$, see
\cite{29} and \cite{6}. If
$\Lambda_{w_c} = \pi _1(T_{w_{c}})$, then we can identify this quotient
with
$(E\otimes \Lambda _{w_c})/Z_W(w_{c})$. Direct computation for the
classical groups identifies
$\Lambda _{w_c}$ with the coroot lattice of a simply connected group
and the action of $Z_W(w_{c})$ with the action of the corresponding
Weyl group, and a similar statement holds for the exceptional groups as
well. The correspondence, which is worked out in
\cite{29}, is as follows:
\roster
\item $\tilde G = SL_{n+1}(\Cee)$ and $c$ of order $f$ dividing $n+1$:
the moduli space of
$G$-bundles is the same as the moduli space of
$SL_{(n+1)/f}(\Cee)$-bundles.
\item $\tilde G = Spin(2n+1)$ and $c$ of order $2$: the moduli space of
$G$-bundles is the same as the moduli space of $Sp(2n-2)$-bundles.
\item $\tilde G = Sp(2n)$, $n$ even and $c$ of order $2$: the moduli
space of
$G$-bundles is the same as the moduli space of $Sp(n)$-bundles (which
has dimension $n/2$).
\item $\tilde G = Sp(2n)$, $n$ odd and $c$ of order $2$: the moduli
space of
$G$-bundles is the same as the moduli space of $Sp(n-1)$-bundles (which
has dimension $(n-1)/2$).
\item $\tilde G = Spin(2n)$, $c$ the element of order $2$ such that the
quotient
$Spin(2n)/\langle c\rangle$ is $SO(2n)$: the moduli space of
$G$-bundles is the same as the moduli space of $Sp(2n-4)$-bundles.
\item $\tilde G = Spin (2n)$, $n$ odd and $c$ of order $4$: the moduli
space of
$G$-bundles is the same as the moduli space of $Sp(n-3)$-bundles.
\item $\tilde G = Spin (2n)$, $n$ even and $c$ one of the exotic
elements of order
$2$: the  moduli space of $G$-bundles is the same as the moduli space of
$Spin(n+1)$-bundles.
\item $\tilde G = E_6$ and $c$ of order $3$: the moduli space
$G$-bundles is the same as the moduli space of $G_2$-bundles.
\item $\tilde G = E_7$ and $c$ of order $2$: the moduli space
$G$-bundles is the same as the moduli space of $F_4$-bundles.
\endroster

In spite of the fact that the moduli spaces are the same as the moduli
spaces of the corresponding simply connected groups, they can be given
by different weighted projective spaces, as we shall see below and in
Part II. This difference is  reflected in the structure of the
component groups of the automorphism groups of the corresponding flat
bundles.

\ssection{5.5. The adjoint bundle.}

We begin with the case of a simply connected group $G$, and identify
the Lie algebra
$\frak z_{\frak g}(\rho)$. 

\definition{Definition 5.15}  For a representation $\rho\: \pi_1(E)
\to T$, define
$$R(\rho ) = \{\,\alpha \in R: \alpha(\rho) = 1\,\}.$$ 
\enddefinition

We then clearly have the following:

\lemma{5.16} Let $\rho$ be a representation from $\pi_1(E)$ to $T$.
\roster
\item"{(i)}" If $\frak g^\alpha$ is the root space corresponding to the 
root 
$\alpha$, then
$$\frak z_{\frak g}(\rho) = \frak h \oplus \bigoplus _{ \alpha\in
R(\rho)}\frak g^\alpha.$$  In particular, $\frak z_{\frak g}(\rho)$ is
a reductive Lie algebra of rank
$r$. 
\item"{(ii)}"  If $\xi_0$ is the bundle corresponding to
$\rho$, then
$$\ad \xi_0 \cong (\frak h\otimes _\Cee\scrO_E) \oplus \bigoplus
_{\alpha\in R}\lambda_{\alpha(\rho)}\cong
\scrO_E^r \oplus \bigoplus _{\alpha\in R}\lambda_{\alpha(\rho)}.$$ 
Hence
$h^0(E;\ad \xi_0)\geq r$, and $h^0(E;\ad \xi_0) = r$ if and only if,
for  every 
$\alpha \in R$, $\alpha (\rho) \neq 1$.
\qed
\endroster
\endstatement

Thus there is an open and Zariski dense subset of flat bundles
$\xi _0$ for which
$h^0(E;\ad \xi_0) = r$, corresponding to the $\rho$ such that
$e_{\alpha(\rho)}\neq 0$ for every root $\alpha$.

Next we turn to a description of the adjoint bundle in the more general
case of
$c$-pairs, corresponding to the case where the group $K$ is not
necessarily simply connected. For a flat
$ K$-bundle with holonomy $\rho$ corresponding to the
${c}$-pair $(x,y)$, let
$\xi$ be the corresponding $G$-bundle. Then $\ad \xi$ is the flat
vector bundle with fiber $\frak g$ associated to the holonomy action of
$\rho$ on $\frak g$. We suppose that
$x$ and $y$ are chosen so that $x\in T^{w_{c}}x_0$ and $y\in N_G(T)$
corresponds to the Weyl element $w_{c}$. As usual, let $\frak h$ be the
Lie algebra of the Cartan subgroup of $G$ and let $\frak g^\alpha$ be
the root space corresponding to the root $\alpha$. We can decompose the
set of roots $R$ into orbits $\bold o$ for the action of $w_{c}$. Since
$w_{c}(x) = x{c}$, the value of a root $\alpha \in \bold o$ on $x$ is
independent of the choice of
$\alpha$. Fix the basis $e_1, e_2$ of $\pi_1(E)$ such that $\rho(e_1) =
\bar x$ and $\rho(e_2) = \bar y$. Let
$L_{x, \bold o}$ be the flat complex line bundle over $E$ with holonomy
given as follows:
$e_1$ acts as
$\alpha(x)$ for any choice of $\alpha \in \bold o$ and $e_2$ acts
trivially. Define the flat vector bundle $V_{y,\bold o}$ as follows: the
fiber of $V_{y,\bold o}$ is the direct sum 
$\bigoplus _{\alpha \in \bold o}\frak g^\alpha$. The holonomy is as
follows: $e_1$ acts trivially, and $e_2$ acts in the natural way, via
$\operatorname{Ad}(y)$. Thus, if $\bold o = \{\alpha, w_{c}(\alpha),
\dots, w_{c}^{d-1}(\alpha)\}$, then there are generators $X_i$ of
$\frak g^{w_{c}^i(\alpha)}$ such that
$\operatorname{Ad}(y) X_i = X_{i+1}$ for $0\leq i\leq d-2$, and 
$\operatorname{Ad}(y) X_{d-1} = c_{\bold o}(y)X_0$ for some $c_{\bold
o}(y)\in U(1)$, which in fact only depends on $\bold o$ and not on
$\alpha$. Note that
$V_{y,\bold o} \otimes L_{x, \bold o} = V_{\bold o}$ is the flat vector
bundle associated to the action of $\pi_1(E)$ on $\bigoplus _{\alpha
\in \bold o}\frak g^\alpha$. Finally,   let
$V_0$ be the corresponding flat vector bundle with fiber $\frak h$: here
$x$ acts trivially, and $y$ acts according to the representation of
$w_{c}$ on $\frak h$. In particular, as holomorphic vector bundles, we
have 
$$V_0 \cong (\frak h^{w_{c}} \otimes \scrO_E) \oplus V_0',$$ where
$V_0'$ is a direct sum of nontrivial torsion line bundles on $E$ and
depends only on $w_{c}$, not on $x$ and $y$.

Define $\bold O(\rho)$ to be the set of $w_c$-orbits $\bold o$ for the
action of
$w_c$ on $R$, such that, for every $\alpha \in \bold o$, $\alpha(x) =
1$ and such that, in the above notation, $c_{\bold o}(y) = 1$. For each
$\bold o\in \bold O(\rho)$, choose an $\alpha\in \rho$ and a nonzero
element 
$$X_{\bold o} = X_0 + \cdots + \operatorname{Ad}(y^{d-1}) X_0 \in
\bigoplus _{i=0}^{d-1}\frak g^{w_c^i(\alpha)},$$ where $d$ is the
number of elements of $\bold o$.

We then have:

\lemma{5.17} In the above notation, we have
$$\frak z_{\frak g}(\rho) = \frak h^{w_{c}} \oplus \bigoplus _{\bold
o\in \bold O(\rho)}\Cee\cdot X_{\bold o}.$$ Moreover, $\frak z_{\frak
g}(\rho)$ is a reductive Lie algebra of rank $r_c =\dim _\Cee \frak
h^{w_{c}} $ and $\frak h^{w_{c}}$ is a Cartan subalgebra of $\frak
z_{\frak g}(\rho)$. Finally,
$$\ad \xi \cong (\frak h^{w_{c}} \otimes \scrO_E) \oplus V_0'
\oplus \bigoplus _{\bold o}(V_{y,\bold o} \otimes L_{x, \bold o}).$$
\endstatement
\proof All of these statements follow from the definition, except the
fact that
$\frak z_{\frak g}(\rho)$ is reductive with rank $r_c$. For the proof
of this fact, we refer to \cite{6}.
\endproof

\corollary{5.18} Suppose that $x$ is regular, in other words, that no
root annihilates $x$. Then  $H^0(E; \ad \xi)\cong
\frak h^{w_{c}}$. Thus the connected component of $\Aut \xi$ is  an
algebraic torus of dimension 
$r_{c} = \dim _\Cee \frak h^{w_{c}}$.
\qed
\endstatement

Of course, it is possible for the conclusions of Corollary 5.18 to hold
without
$x$ being a regular element.

\ssection{5.6. Comparison of complex structures.}

In this subsection, we relate the global description of the moduli
space given by Lie theory to its complex-analytic structure. For
example, in the simply connected case, we have identified the moduli
space $\Cal M(G)$ as a set with the normal complex variety $(E\otimes
_\Zee\Lambda)/W$. Our goal now is to show that they are isomorphic as
complex varieties:

\theorem{5.19} Suppose that $G$ is simply connected, and let $\Lambda$
be the coroot lattice of $G$. Then, as normal projective varieties,
$\Cal M(G) \cong (E\otimes _\Zee\Lambda)/W$. More precisely, there is a
holomorphic map $(E\otimes _\Zee \Lambda)/W \to \Cal M(G)$ and it is an
isomorphism of normal projective varieties. More generally, suppose
that $G$ is semisimple and let $c$ be an element of the center of the
universal cover $\tilde G$ of $G$, such that $G =\tilde G/\langle
c\rangle$. Then, in the notation of Theorem
\rom{5.13}, there exists a complex structure on
$(T^c\times T_{w_c})/Z_W(w_c)$ for which it is an irreducible, normal
projective variety, and a holomorphic isomorphism from $(T^c\times
T_{w_c})/Z_W(w_c)$ to $\Cal M(G,c)$.
\endstatement
\proof We begin with the simply connected version.  The space
$(E\otimes _\Zee \Lambda)/W$ has a holomorphic structure, coming from
the holomorphic structure on $E\otimes _\Zee \Lambda \cong E^r$. Indeed,
$(E\otimes _\Zee \Lambda)/W$ is a normal projective variety. Let us
show that this holomorphic structure agrees with that of every coarse
moduli space of semistable
$G$-bundles modulo S-equivalence.  First we claim that there is a
universal tautological holomorphic bundle over $E\times (E\otimes _\Zee
\Lambda)$:

\lemma{5.20} There is a universal holomorphic $H$-bundle $\Xi_H$ over
$E\times (E\otimes _\Zee\Lambda)$, such that the restriction of $\Xi_H$
to the slice $E\times \{p\}$ is induced from the flat $T$-bundle
corresponding to the point $p\in E\otimes _\Zee\Lambda$.
\endstatement
\proof  One way to describe $\Xi_H$ is as follows: fix a basis
$\alpha_1, \dots, \alpha _r$ of simple roots, and use the corresponding
basis $\{\alpha_1\spcheck,
\dots, \alpha _r\spcheck\}$ of coroots to identify
$\Lambda$ with $\Zee^r$. Let $\Cal P_i$ be the Poincar\'e bundle over
$E\times E$, pulled back to $E\times (E\otimes _\Zee\Lambda)$ via the
projection
$E\otimes _\Zee\Lambda \cong E^r \to E$ defined by the  basis 
$\{\alpha_1\spcheck,
\dots, \alpha _r\spcheck\}$ of $\Lambda$. We then define $\Xi_H$ as an
$H$-bundle by the ordered $r$-tuple $(\Cal P_1, \dots, \Cal P_r)$. It
is routine to check that $\Xi_H$ has the desired property.
\endproof

Another way to think of the bundle $\Xi_H$ defined above is as follows:
there is a universal holomorphic $(\Cee^*)^r$-bundle over $E\times
E^r$, defined by taking
$r$ copies of the Poincar\'e bundle. The identification of $(\Cee^*)^r$
with $H$ given by the coroot basis $\{\alpha_1\spcheck,
\dots, \alpha _r\spcheck\}$ defines a homomorphism $(\Cee^*)^r \to G$
of algebraic groups with image $H$, and thus a universal $H$-bundle.

The universal bundle $\Xi_H$ above defines a holomorphic,
$W$-invariant morphism
$E\otimes _\Zee \Lambda \to \Cal M(G)$. By Theorem 5.10 and Corollary
5.12, the morphism
$(E\otimes _\Zee \Lambda)/W \to \Cal M(G)$ is a bijection on the set of
points. It follows from Zariski's main theorem that the map $(E\otimes 
_\Zee
\Lambda)/W \to \Cal M(G)$ is an isomorphism.

The non-simply connected case is similar, but slightly more involved.
First let us  define a natural complex structure on
$(T^c\times T_{w_c})/Z_W(w_c)$, depending of course on the complex
structure on $E$. To do so, let
$T_0 = (T^{w_{c}})^0$ be the identity  component of  $T^{w_{c}}$. There
is the natural quotient homomorphism $T_0 \to T_{w_{c}}$, and it is a
surjection. Thus there is a finite surjective homomorphism of real tori
from
$T_0\times T_0$ to $T_0 \times T_{w_{c}}$. If $\Lambda _0 =\pi_1(T_0)$,
then the identification $T_0 \cong \Lambda _0\otimes _\Zee U(1)$ and
the isomorphism
$$T_0\times T_0 \cong \Hom(\pi_1(E), T_0) \cong \Hom (\pi_1(E),
U(1))\otimes _\Zee
\Lambda _0 \cong E\otimes_\Zee \Lambda _0$$ show that $T_0\times T_0$
is naturally the abelian variety
$E\otimes_\Zee \Lambda _0$. Since $T_0 \times T_{w_{c}}$ is the
quotient of $T_0
\times T_0$ by a finite subgroup, it is also an abelian variety. The
identity component of
$T^{w_{c}} \times T_{w_{c}}$ has the structure of an abelian variety,
and thus $T^{w_{c}} \times T_{w_{c}}$ is a finite union of abelian
varieties. Finally,
$T^c \times T_{w_{c}}$ is isomorphic to $T^{w_{c}} \times T_{w_{c}}$
and so inherits a complex structure via this isomorphism. We claim that
the action of
$Z_W(w_{c})$ defined above, which is by affine maps of the real torus
to itself, is holomorphic with respect to this complex strucure. It is
enough to check the differential of the action, which corresponds to
the linear action of $Z_W(w_{c})$ on $E\otimes_\Zee \Lambda _0$ via its
action on $\Lambda_0$. Thus
$Z_W(w_{c})$ acts holomorphically, and the quotient variety
$(T^c\times T_{w_c})/Z_W(w_c)$ is an irreducible normal projective
variety.

Next we have the following analogue of Lemma 5.20, which constructs a
universal holomorphic $G$-bundle parametrized by  $T^c\times T_{w_c}$:

\lemma{5.21} There exists a universal holomorphic $G$-bundle $\Xi _1$
over the space $E\times (T^c\times T_{w_c})$, and hence the induced map
$$(T^c\times T_{w_c})/ Z_W(w_{c})
\to \Cal M(G,c)$$ is an isomorphism of complex spaces.
\endstatement
\noindent {\it Proof.}  We begin by describing a universal holomorphic
bundle on the connected set
$(T_0\cdot x_0)
\times (T_0\cdot v)$. Recall that we chose a fixed $c$-pair $(x_0, v)$
in the discussion prior to the proof of Theorem 5.13, and then were
able to write an arbitrary $c$-pair, up to conjugation, as $(tx_0,sv)$,
with $t\in T^{w_c}$ and $s\in T$. Let
$\xi_1$ be the flat
$G$-bundle corresponding to the choice of holonomy $\bar x_0, \bar v$,
where $\bar x_0, \bar v$ are the images of $x_0, v$ in $G$. Thus the
transition functions of $\xi_1$ live in the subgroup of $G$ generated
by $\bar x_0$ and $\bar v$. Let $r_c$ be the rank of
$\Lambda_0$, or in other words the dimension of
$T^{w_c}$. Let $H_0$ be the subtorus of $H$ corresponding to $T_0$, in
other words the connected component of the identity of the fixed set of
$v \in N_G(H)$ acting on $H$, and let
$\bar H_0$ be its image in $G$. Choosing an integral basis for
$\Lambda_0$ identifies
$H_0$ with $(\Cee^*)^{r_c}$ and $E\otimes_\Zee \Lambda_0$ with
$E^{r_c}$. For $i=1,
\dots, r_c$, let $\Cal P_i$ be the pullback to $E\times (E\otimes_\Zee
\Lambda_0)$ of the Poincar\'e bundle on the $i^{\text{th}}$ factor and
let $\Cal P$ be the
$H_0$-bundle corresponding to $(\Cal P_1, \dots, \Cal P_{r_c})$. As in
the proof of Lemma 5.20,  we can think of $\Cal P$ as corresponding to
the image of the universal topologically trivial
$(\Cee^*)^{r_c}$-bundle via the inclusion $(\Cee^*)^{r_c}
\cong H_0 \subset G$. There is an induced
$\bar H_0$-bundle $\overline{\Cal P}$ over $E\times (E\otimes_\Zee
\Lambda_0)$.

The images $\bar x_0,
\bar v$ of $x_0, v$ in $G$ commute with
$\bar H_0$ and with each other. Thus the structure group of
$\overline{\Cal P}$ centralizes the transition functions for
$\pi_1^*\xi _1$. It follows that there is a holomorphic $G$-bundle $\Xi
_1$ defined by taking the product of the transition functions, and it
is easily seen to be a universal bundle in the appropriate sense. A
similar construction handles the other components.

By the universal property of the moduli space, there is then an induced
holomorphic map from
$T^{w_{c}}\times T_0$ to
$\Cal M(G,c)$ which is holomorphic. Clearly, it factors through the
quotient map to
$T^{w_{c}}\times T_{w_{c}}$ and through the action of
$Z_W(w_{c})$, since these actions do not change the isomorphism class
of a flat bundle. By Theorem 5.13, the induced morphism $(T^c\times
T_{w_c})/ Z_W(w_{c})
\to \Cal M(G,c)$ is a holomorphic bijection between two normal
varieties and hence it is an isomorphism. Thus we have proved Lemma
5.21 and with it Theorem 5.19.
\endproof

Of course, the description of the moduli space above gives us another
way to describe its singularities. Let us just work out the case where
$G$ is simply connected for simplicity. Let $\rho$ correspond to a
point of $E\otimes _\Zee\Lambda$. Then by
\cite{6}, the isotropy group of $\rho$ is a semidirect product
$W(\rho)\rtimes I$. The local structure of the quotient is given by
looking at the action of the isotropy group on the tangent space to
$E\otimes _\Zee\Lambda$, which is canonically identified with $\frak
h$. Thus, as in Corollary 5.3, we see again that the moduli space is
locally given as
$(\frak h/W(\rho))/I$.

\ssection{5.7. Looijenga's theorem and its generalizations.}

In the case where $G$ is simply connected, we have described the moduli
space $\Cal M(G)$ as the complex variety $(E\otimes _\Zee\Lambda)/W$.
The precise nature of this variety has been determined in a theorem due
to Looijenga \cite{24} and Bernshtein-Shvartsman \cite{4}:

\theorem{5.22} Let $E$ be an elliptic curve and let
$\Lambda$ be the coroot lattice of a simple root system $R$ with Weyl
group $W$. Then $(E\otimes _\Zee
\Lambda)/W$ is the weighted projective space $\Pee(g_0, \dots, g_r)$,
in other words it is the quotient of $\Cee^{r+1} - \{0\}$ by the
diagonal action of $\Cee^*$ given by
$$\lambda \cdot (z_0, \dots, z_r) = (\lambda^{g_0}z_0, \dots,
\lambda^{g_r}z_r).
\qquad \qed$$
\endstatement

The proof of \cite{24} and \cite{4} makes use of formal theta functions
for a complexified affine Weyl group. We shall give a different proof
of Theorem 5.22 in Part II, based on a study of the deformations of
certain unstable bundles over special maximal parabolic subgroups of
$G$. Note that the holomorphic structure of the moduli space is
independent of the complex structure on $E$. Furthermore, the local
information concerning the singular points of the moduli space follows
immediately from Theorem 5.22.

In the non-simply connected case, we have the following generalization
of Theorem 5.22, which will also be proved in Part II:

\theorem{5.23} Let $G$ be a simple group and let $c$ be an element of
the center of the universal cover $\tilde G$ of $G$, such that $G
=\tilde G/\langle c\rangle$. Then $\Cal M(G,c)$ is
isomorphic to a weighted projective space with weights $n_{\bold
o}g_{\bold o}$. \qed
\endstatement

Note that, as opposed to the simply connected case, the knowledge of
the integers
$n_{\bold o}g_{\bold o}$ is more information than the isomorphism class
of the weighted projective space. For example, for the unliftable
$SO(2n)$-bundles, the integers $n_{\bold o}g_{\bold o}$ are all equal
to $2$, indicating that the weighted projective space is an ordinary
projective space but that the isotropy group, which acts trivially on
the Zariski tangent space $H^1(E; \ad \xi)$, has order $2$. A more
interesting example is given by the unliftable $SO(2n+1)$-bundles or
the unliftable $Sp(4n)/\{\pm\Id\}$-bundles. In this case, the weighted
projective space is a $\Pee(1, 2, \dots, 2)$, which is isomorphic as a
projective variety to $\Pee^n$. However, the isotropy behaves quite
differently. In these examples, as we have noted in the discussion
after Theorem 5.13 and Corollary 5.14, we can identify the moduli space
with $(E\otimes _\Zee\bar \Lambda)/\bar W$, where
$\bar \Lambda$ is the coroot lattice for a simply connected group of
smaller rank and $\bar W$ is the corresponding Weyl group. However, the
weighted projective space for the smaller group (which is of type
$C_n$) is a $\Pee(1, \dots, 1)$.

\ssection{5.8. The case where $G$ is not simple.}

We have given an explicit description of the moduli space $\Cal M(G)$
in case $G$ is simple and simply connected; it is a weighted projective
space. In case $G$ is simply connected but only semisimple, $G$ is a
product of simple and simply connected groups
$G_i$, and the moduli space is a product of the corresponding moduli
spaces. Of course, this is also clear from the description of the
moduli space as $(E\otimes _\Zee\Lambda)/W$, which did not need $G$ to
be simple, since in this case $\Lambda =
\bigoplus _i\Lambda _i$ and $W = \prod _iW_i$.

The next case to consider is the case where $G$ is semisimple. As in
the discussion prior to Definition 5.7, we begin by assuming that $G$
is of the form $\tilde G/\langle c\rangle$, where $\tilde G$ is simply
connected and $c$ is an element of the center of $\tilde G$, and  we
want to describe the component $\Cal M(G,c)$. The moduli space $\Cal
M(G,c)$ is the same as the set of conjugacy classes of $c$-pairs in
$\tilde K$. Suppose that $\tilde K =\prod_iK_i$, where the $K_i$ are
the simple factors of $\tilde K$, and that $c=(c_1, \dots, c_n)$, where
$c_i$ is an element of the center of $K_i$. Let $G_i$ be the complex
form of the compact and simply connected group $K_i$ and let $\bar G_i
= G_i/\langle c_i\rangle$. By (iii) of Lemma 5.8, the moduli space
$\Cal M(G,c)$ is a product of the corresponding moduli spaces $\Cal
M(\bar G_i,c_i)$. By Theorem 5.23, each space $\Cal M(\bar G_i,c_i)$ is
a weighted projective space. Thus $\Cal M(G,c)$ is a product of
weighted projective spaces.

In the case where $G$ is semisimple but not necessarily of the form
$\tilde G/\langle c\rangle$, the discussion prior to Definition 5.7
shows that the moduli space $\Cal M(G,c)$ is a quotient of the moduli
space $\Cal M(\tilde G/\langle c\rangle,c)$ by a finite abelian group.
Thus, $\Cal M(G,c)$ is a quotient of a product of weighted projective
spaces by a finite abelian group. Of course, the action of this group
can in principle be made explicit. Examples show that, even if $G$ is
simple, the moduli space $\Cal M(G,c)$ need not be a weighted
projective space.

Next suppose that $G = \Cal C$ is abelian, and is the complexification
of the compact torus
$T$. Let
$\Lambda (T) = \pi _1(T) = \pi_1(\Cal C)$. After choosing an
isomorphism $\Cal C
\to (\Cee^*)^d$, a holomorphic
$\Cal C$-bundle over
$E$ is the same thing as a collection of $d$ line bundles over $E$.
Thus there is a discrete invariant, the set of $d$ degrees,
corresponding to the choice of $c\in \pi_1(\Cal C)$. All components
are isomorphic to the identity component, which is isomorphic to
$(\Pic^0E)^d$. More invariantly, the identity component is given as
$E\otimes _\Zee\Lambda(\Cal C)$. Of course, we can see this from the
point of view of Yang-Mills connections as well. Since $T$ is abelian,
a homomorphism $\rho\:
\Gamma_\Ar \to T$ factors into a homomorphism $U(1)\times \pi_1(E) \to
T$, or equivalently a pair of homomorphisms $\rho_1\: U(1) \to T$ and
$\rho_2\: \pi_1(E)
\to T$. The first homomorphism $\rho_1$ records the multidegree, and a
fixed choice of a Yang-Mills connection for $\boldsymbol
\Gamma_\Ar$ defines a fixed holomorphic $\Cal C$-bundle of any given
multidegree. As for $\rho_2$, writing $T$ invariantly as
$U(1)\otimes_\Zee \Lambda(\Cal C)$, we have
$$\Hom (\pi_1(E), T) = \Hom (\pi_1(E),U(1)) \otimes_\Zee \Lambda(\Cal
C).$$ Now $\Hom (\pi_1(E),U(1))$ is the group of flat
$U(1)$-connections on $E$, and is thus identified with $\Pic^0E$. A
choice of base point of $E$ identifies $\Pic^0E$, as an algebraic
group, with $E$. Thus, from  the viewpoint of Yang-Mills connections,
we have again identified the identity component of the moduli space
with 
$E\otimes _\Zee\Lambda(\Cal C)$.

Finally suppose that $G$ is written as $\Cal C\times _FD$, where $\Cal
C$ is the identity component of the center of $G$, $D=DG$ is the
derived group and is semisimple, and 
$F$ is a finite central subgroup of $D$. Let $\Cal C_\Ar$ be the
maximal compact subgroup of $\Cal C$. There is an exact sequence
$$\{1\} \to \Cal C \to G \to \bar D\to \{1\},$$ where $\bar D = D/F$. A
$G$-bundle $\xi$ then produces a $\bar D$-bundle
$\bar \xi$, and $\xi$ is semistable if and only if $\bar \xi$ is
semistable. Moreover, since $\underline{\Cal C} \cong (\scrO_E^*)^c$,
$H^2(E;
\underline{\Cal C}) =0$, and so  every $\bar D$-bundle lifts to a
$G$-bundle. The choice of lift is well-defined up to an element in
$H^1(E; \underline{\Cal C})$, which, as we have seen above, is a
complex Lie group whose identity component is the complex torus
$E\otimes _\Zee\Lambda(\Cal C)$. In general, the action of
$H^1(E;
\underline{\Cal C})$ on the moduli space is not faithful, but has a
finite kernel. To see this, we again look at the problem from the point
of view of central representations of the group $\Gamma_\Ar$. Suppose
that $\xi$ is S-equivalent to the Yang-Mills bundle with holonomy
$\rho$, where $\rho\: \Gamma_\Ar\to K$ is a central representation.
Then 
$\bar \xi$ is S-equivalent to a flat $DK/F$-bundle, where $DK$ is the
derived group of $K$ and is the maximal compact subgroup of $D$,
corresponding to the homomorphism $\Gamma_\Ar \to DK/F$. The fiber of
the natural map from the set of isomorphism classes of central
representations
$\Gamma_\Ar\to K$ to the set of isomorphism classes of representations
$\pi_1(E)\to DK/F$ is acted on transitively by the group of
homomorphisms $\chi\: \Gamma_\Ar
\to \Cal C_\Ar$, which is just  $H^1(E; \underline{\Cal C})$.  To
determine the kernel of this action, let
$c = c_1(\bar \xi) \in H^2(E; F)
\cong F$.  Then the images under $\bar \rho$ of two generators for
$\pi_1(E)$ lift to elements $x,y\in DK$ with $[x,y] = c$. Hence $\rho$
itself lifts to a representation from $\Gamma_\Ar$ to $\Cal C_\Ar
\times _{\langle c\rangle}DK$, and $\xi$ lifts to a $\Cal C\times
_{\langle c\rangle}D$-bundle. As in the discussion of the semisimple
case, we assume that
$F=\langle c\rangle$, in other words that we have lifted $\xi$ as far
as possible, since the moduli space in the general case will be
finitely covered by the moduli space of such bundles. It is then an
easy exercise to check that
$\rho\chi$ is conjugate to
$\rho$ if and only if
$\operatorname{Im}\chi\subseteq \langle c\rangle$. In this way, the
moduli space for $G$-bundles becomes an \'etale fiber bundle over the
moduli space for $\bar D$-bundles, where the fibers are finite
quotients of copies of the complex torus 
$H^1(E; \underline{\Cal C})$.

We can summarize the results above as follows:

\corollary{5.24} Let $G$ be a reductive group and let $c \in \pi
_1(G)$. Then there is a finite cover of $\Cal M(G,c)$ which is an
\'etale fiber bundle over a product of weighted projective spaces, and
where the fiber is isogenous to a product of copies of $E$ and in
particular is a complex torus.
\qed
\endstatement

\section{6. Regular bundles.}

Our goal now will be to describe the bundles whose automorphism groups
have minimal possible dimension. We call such a bundle {\sl regular}.
It turns out that there is a unique regular bundle in each
S-equivalence class. For a Zariski open and dense subset of the moduli
space, the Yang-Mills representative is also regular, but they will be
different at special points. It is natural to study regular bundles for
many reasons. For example, in moduli questions, their deformation
theory is much closer to reflecting the local structure of the moduli
space than the Yang-Mills representative in general. Also, it is not
possible to fit together the Yang-Mills bundles to form a universal
bundle, even locally, but, as we shall see in Part II, local and in
many cases global universal bundles can be constructed using the
regular representatives. Finally, because they are more complicated,
regular bundles are interesting to study in their own right. To
illustrate this point, we work out the automorphism groups of regular
bundles in this section. In the final section, we shall determine the
regular bundles for the classical groups.

\ssection{6.1. Definition of a regular bundle.}

We begin by recalling some standard facts about nilpotent elements in
reductive Lie algebras \cite{7, Chap\. 8, VIII,
\S 11} or \cite{21}. Let
$\frak g$ be a complex semisimple Lie algebra of rank $r$ and let $X\in
\frak g$ be such that
$\ad X$ is nilpotent. Every such element can be completed to an $\frak
{sl}_2$-triple, in other words, a representation of $\frak {sl}_2$ on
$\frak g$ for which
$X$ is the image of a nonzero nilpotent element. For such $X$, $\dim
\Ker (\ad X)
\geq r$. If  $\dim \Ker (\ad X) = r$, then $X$ is called a {\sl
principal\/} nilpotent element. All principal $X$ are conjugate under
$G$. If $X$ is principal, as
$\frak {sl}_2$-modules, we have
$$\frak g = \bigoplus _{i=1}^r\Sym ^{2d_i-2}V(2)=\bigoplus _{i=1}^r
V(2d_i-1),$$ where $V(2)$ denotes the standard representation of $\frak
{sl}_2$, and, for an integer $k> 1$, $V(k) = \Sym ^{k-1}V(2)$ is the
unique irreducible representation of
$\frak {sl}_2$ of dimension $k$. Here $r$ is the rank of $\frak g$ and
the $d_i$ are the Casimir weights; $d_i = m_i + 1$, where the $m_i$ are
the exponents of
$\frak g$. For an element $X$ of $\frak g$ which is not principal, the
corresponding representation of $\frak {sl}_2$ on $\frak g$  breaks up
as 
$$\frak g = \bigoplus _{i=1}^sV(k_i),$$ where $s\geq r+2$ \cite{20,
Prop\. 4.11}. Equivalently, if $X$ is not principal, then 
$$\dim \Ker (\ad X) \geq r+2.$$
If $X$ is principal, the centralizer
$\frak a$ of $X$ in $\frak g$ is an abelian nilpotent Lie algebra.
Conversely, if the centralizer of $X$ is abelian, then $X$ is
principal. Clearly, all of the above extends  to the case where
$\frak g$ is only assumed to be reductive, and $\frak c$ is its center.
In this case, we have an isomorphism of
$\frak {sl}_2$-modules
$$\frak g \cong \frak c\oplus \bigoplus _{i=1}^r V(2d_i-1),$$ where the
$d_i$ are the Casimir weights of the semisimple part of $\frak g$, and
$\frak c$ is the center of $\frak g$, viewed as a trivial 
$\frak {sl}_2$-module.

\definition{Definition 6.1} A semistable $G$-bundle is said to be {\sl
regular} if it is given by normalized transition functions
$\{h_{ij}\exp(f_{ij}X)\}$ where $X$ is a principal nilpotent element of
$\frak z_{\frak g}(\rho)$. (Here, $\rho$ is the holonomy representation
of the S-equivalent Yang-Mills $K$-bundle.)
\enddefinition

\corollary{6.2}  Each S-equivalence class of semistable $G$-bundles has
a unique \rom(up to isomorphism\rom) regular representative. The
regular representative is also the Yang-Mills representative if and
only if $\frak z_{\frak g}(\rho)$ is abelian. Finally, there is a
Zariski open and dense subset of the moduli space $\Cal M(G,c)$ where
the regular representative is also the Yang-Mills representative. 
\endstatement

\proof The principal nilpotent elements of $\frak z_{\frak g}(\rho)$ are
all conjugate. Thus, the first statement follows immediately from
Theorem 4.1. As for the second statement, a regular representative is
Yang-Mills if and only if the principal nilpotent elememt $X$ of the
reductive Lie algebra $\frak z_{\frak g}(\rho)$ is zero, which is the
case if and only if $\frak z_{\frak g}(\rho)$ is abelian. The final
statement follows from Lemma 5.16 (ii) and Corollary 5.18.
\endproof

It follows from Lemma 5.17 that the rank $r_c$ of $\frak z_{\frak
g}(\rho)$ depends only on the topological type of the bundle $\xi$. For
example, if $G$ is simply connected, then $r_c$ is always equal to the
rank of $G$. In general, if $\xi$ lies in $\Cal M(G,c)$, and $w_c$ is
the Weyl group element corresponding to the element $c$, then $r_c$ is
equal to $\dim _{\Cee}\frak h^{w_c}$.  With this said, we have:

\corollary{6.3} Let $\xi$ be a semistable $G$-bundle associated to the
pair
$(\rho, X)$, where $\rho$ is the holonomy representation of the
Yang-Mills
$K$-bundle S-equivalent to $\xi$ and $X\in \frak z_{\frak g}(\rho)$ is
a nilpotent element used to define the normalized transition functions
of $\xi$. Let
$r_c$ be the rank of the reductive Lie algebra
$\frak z_{\frak g}(\rho)$. Then $\dim \Aut_G\xi = r_c$ if  $X$ is
regular, and if $X$ is not regular, then $\dim \Aut_G\xi \geq r_c+2$.
\qed
\endstatement

In principle, we can describe  the vector bundle $\ad \xi$ when $\xi$
is regular. For example, we have the following description of the
adjoint bundle associated to the trivial bundle:

\proposition{6.4} Let $\xi$ be the  unique regular bundle which is
S-equivalent to the trivial bundle. Then, as vector bundles over
$E$,
$$\ad \xi \cong \scrO_E^c\oplus \bigoplus _iI_{2d_i-1},$$ where $c$ is
the dimension of the center of $G$ and the $d_i$ are the Casimir
weights of
$G$.
\qed
\endstatement

Similar but more involved statements can describe in principle the
bundle $\ad
\xi$ for an arbitrary regular bundle $\xi$. To give one such statement
along these lines,  let
$\rho $ be be the holonomy representation of the  Yang-Mills bundle
S-equivalent to
$\xi$. Then $\frak z_{\frak g}(\rho)$ is a direct sum of its center
$\frak c(\rho)$ together with  $N$ simple  factors
$\frak g_i$.  Let $r_i$ be the rank of
$\frak g_i$, and let $d_{ij}, 1\leq j\leq r_i$ be the Casimir weights
of $\frak g_i$. Then the maximal subbundle of $\ad \xi$ which is
filtered by subbundles whose successive quotients are
$\scrO_E$ is
$$(\ad \xi)_{\scrO_E} = (\frak c(\rho)\otimes \scrO_E)\oplus \bigoplus
_{i=1}^N\bigoplus _{j=1}^{r_i}I_{2d_{ij}-1}.$$
  From this, it is possible in principle to give a complete description
of the vector bundle $\ad
\xi$.

\ssection{6.2. Automorphisms of regular bundles.}

We study in more detail the
structure of the automorphism group of a  regular semistable holomorphic
bundle $\xi$. Let $\rho$ be the holonomy representation corresponding
to  the unique Yang-Mills bundle S-equivalent to
$\xi$. Define $A(\xi)=\Aut _G(\xi)$, which according to Theorem 4.1 we
can view as a subgroup of  $Z_ G(\rho)$. Let
$A^0(\xi)$ be the component of the identity in $A(\xi)$. Let us analyze
the Lie algebra of $A^0(\xi)$. Let $D(\rho)$ be the derived
subgroup of the identity component $Z^0_G(\rho)$ of $Z_G(\rho)$
and let $\frak d(\rho)$ be its Lie algebra, so that $\frak d(\rho)$ is
the derived subalgebra of  $\frak z_{\frak
g}(\rho)$.

\lemma{6.5}  Let $\frak c(\rho)$ be the center of $\frak z_{\frak
g}(\rho)$ and let $\frak d(\rho)$ be its derived subalgebra,
so that 
$$\frak z_{\frak g}(\rho) \cong \frak c(\rho)\oplus \frak d(\rho).$$ 
Then $X\in \frak d(\rho)$ and its centralizer in
$\frak d(\rho)$ is an abelian nilpotent subalgebra $\frak n$.
Thus, the centralizer of $X$ in $\frak z_{\frak g}(\rho)$ is the
product $\frak c(\rho)\oplus \frak n$. In particular, the Lie algebra
of $A^0(\xi)$ is the product of a central semisimple subalgebra and an
abelian nilpotent algebra, and hence it is abelian.
\endstatement

\proof Everything in this statement is clear except the description of
the centralizer of $X$ in $\frak d(\rho)$, which follows
from the remarks above concerning the centralizers of principal
nilpotent elements.
\endproof

\corollary{6.6}
$A^0(\xi)$  is a connected abelian complex linear algebraic group of
the form
$\Cal C\times U$, where $\Cal C$ is an algebraic torus and $U$ is an
abelian unipotent group.  In fact,
$\Cal C$ is the identity component of the center  of $Z^0_ G(\rho)$,
$U$ is a unipotent subgroup of the derived subgroup $D(\rho)$ of
$Z_G^0(\rho)$, and the Lie algebra of $U$ is the centralizer of the
principal nilpotent element $X$ in  $\frak d(\rho)$.\qed
\endstatement

The next step in understanding $A(\xi)$ is to describe
$A(\xi)\cap Z_G^0(\rho) = A^1(\xi)$. Since $Z_G^0(\rho)$ is a
connected, reductive group, it decomposes as
$$Z_G^0(\rho)=\Cal C(\rho)\times _FD(\rho)$$ where $\Cal C(\rho)$,
the identity component of the center of
$Z_G^0(\rho)$, is an algebraic torus, $D(\rho)$, the derived
subgroup of $Z_G^0(\rho)$, is a complex semisimple group, and $F$ is a
finite group which embeds into the center of $D(\rho)$.

\lemma{6.7} Let $\xi$ be a regular semistable $G$-bundle, let $A(\xi)
=\Aut_G(\xi)$, which we view as a subgroup of $Z_G(\rho)$, and let
$A^1(\xi) = A(\xi) \cap Z_G^0(\rho)$. Then:
\roster
\item"{(i)}" $A^1(\xi) =\Cal  C(\rho) \times _F({\Cal
Z}D(\rho)\times U)$, where
${\Cal Z}D(\rho)$ is the center of $D(\rho)$. In particular,
$A^1(\xi)$ is abelian and $\pi_0(A^1(\xi)) = {\Cal Z}D(\rho)/F$.
\item"{(ii)}" The inclusion $A(\xi) \subseteq Z_G(\rho)$ induces  an
isomorphism
$$A(\xi)/A^1(\xi) @>{\cong}>> \pi _0(Z_G(\rho)).$$
\endroster
\endstatement

\proof The group $A^1(\xi)$, which is the  centralizer of $X$ in $Z_
G^0(\rho)$, is  of the form
$\Cal  C(\rho)\times _F D(\rho)(X)$, where $D(\rho)(X)$ is
the subgroup of $D(\rho)$ centralizing $X$. By standard results 
\cite{20,
\S  4.7}, if $X$ is a principal nilpotent element of
$\frak d(\rho)$, then
$D(\rho)(X) = {\Cal Z}D(\rho)\times U$, where ${\Cal
Z}D(\rho)$ is the (finite) center of
$D(\rho)$ and $U$ is the abelian unipotent subgroup of
$D(\rho)$ corresponding to the subalgebra of
$\frak d(\rho)$ centralizing $X$. Thus $A^1(\xi) = \Cal 
C(\rho)\times _F({\Cal Z}D(\rho)\times U)$. Clearly, $A^1(\xi)$ is
abelian, but need not be connected. 

To see (ii), let 
$g\in Z_G(\rho)$ be given. Then $g$ normalizes
$Z_G^0(\rho)$, and clearly
$(\operatorname{Ad}g)(X)$ is again a principal nilpotent element of
$\frak z_{\frak g}(\rho)$. By \cite{21}, every two principal nilpotent
elements of 
$\frak z_{\frak g}(\rho)$ are conjugate under $Z_G^0(\rho)$. Let $g'\in
Z_G^0(\rho)$ be such that $(\operatorname{Ad} (g'g))(X) = X$. Then
$g'g\in A(\xi)$, and clearly $g'g$ and
$g$ have the same image in $\pi_0(Z_G(\rho))$. Thus, the induced
homomorphism
$A(\xi) \to
\pi _0(Z_G(\rho))$ is surjective. Since its kernel is clearly
$A^1(\xi)$, this completes the proof.
\endproof

Let us give two simple examples of the structure of $A(\xi)$. If $\xi$
is the  trivial $G$-bundle, then $Z_G(\rho) = G$ and $A(\xi) =
A^1(\xi)$ is equal to ${\Cal Z} G\times U$, where
${\Cal Z} G$ is  the center of $G$ and $U$ is a connected, unipotent
abelian group. For another example, let
$G= SL_n(\Cee)$. Suppose that $\lambda_i$, $i=1, \dots, k$, are
pairwise distinct line bundles on
$E$ of degree zero, and let 
$\xi$ correspond to the regular vector bundle $\bigoplus
_{i=1}^kI_{d_i}(\lambda_i)$. Set $d =\gcd \{d_i\}$. Then $A^1(\xi) =
A(\xi)$, and it is easy to check that $A(\xi)/A^0(\xi) \cong
\Zee/d\Zee$.

\ssection{6.3. The action of the component group of $A(\xi)$ on its
  Lie algebra $\frak a(\xi)$.}

We fix a regular bundle $\xi$ and denote by $A(\xi)$ is group of
automorphisms. Let $\rho$ be the holonomy representation of the
Yang-Mills
$K$-bundle S-equivalent to $\xi$. As we have seen, the Lie algebra
$\frak a (\xi)$ of $A^0(\xi)$ is abelian. Thus, the adjoint action of
$A(\xi)$ on $\frak a(\xi)$ factors through an action of the component
group $A(\xi)/A^0(\xi)$. The purpose of this subsection is to determine
the action of the component group $A(\xi)/A^0(\xi)$ on $\frak a (\xi)$.
First notice that since $A^1(\xi)$ is abelian,
$A^1(\xi)/A^0(\xi)$ acts trivially on 
$\frak a (\xi)$ and there is an induced action of $A(\xi)/A^1(\xi)$ on 
$\frak a (\xi)$. 

In Lemma 6.5 we showed that 
$\frak a (\xi) \cong \frak c(\rho) \oplus \frak n$, where $\frak
c(\rho)$ consists of semisimple elements in $\frak z_{\frak  g}(\rho)$
and
$\frak n$ is a nilpotent subalgebra of $\frak z_{\frak g}(\rho)$. 
Clearly, the action of the group $A(\xi)/A^1(\xi)$ on the Lie algebra
must preserve this splitting. The inclusion of $A(\xi)\subseteq
Z_G(\rho)$ induces an isomorphism of
$A(\xi)/A^1(\xi)\cong \pi_0(Z_G(\rho))$. Thus, we have a realization of
$\pi_0(Z_G(\rho))$ as a group of automorphisms of $\frak a(\xi) = \frak
c(\rho)\oplus
\frak n$. The action of $\pi_0(Z_G(\rho))$ on $\frak c(\rho)$ is simply
the conjugation action of $\pi_0(Z_G(\rho))$ on the center $\frak
c(\rho)$ of the Lie algebra of $Z_G^0(\rho)$. 

Next we claim that there is also an action of $\pi_0(Z_G(\rho))$ on a
Cartan subalgebra $\frak h(\rho)$ of $\frak z_{\frak  g}(\rho)$. This
is a general statement about outer automorphisms:

\claim{6.8} Fix a Cartan subalgebra $\frak h(\rho)$ for $\frak z_{\frak
g}(\rho)$, with $H(\rho)$ the corresponding subgroup, and a Weyl chamber
$C_0\subset
\frak h(\rho)$. Then for an outer automorphism
$\gamma $ of $Z^0_G(\rho)$, there is a lift of $\gamma$ to an
automorphism $y_\gamma$ which normalizes $\frak h(\rho)$ and $C_0$.
Moreover, the element $y_\gamma$ is unique up to multiplication by some
$h\in H(\rho)$.
\endstatement

\proof Given an outer automorphism $\gamma$, let $y$ be some lift of 
$\gamma$ to an automorphism of $Z_G(\rho)$. Conjugation by $y$  sends
$\frak h(\rho)$ to  another Cartan subalgebra $\frak h'$ for $\frak
z_{\frak g}(\rho)$. Since $\frak h(\rho)$ and $\frak h'$ are conjugate
in
$Z_G^0(\rho)$, after multiplying $y$ by some element of $Z_G^0(\rho)$,
we can assume conjugation by $y$ normalizes $\frak h(\rho)$. Then
conjugation by $y$ sends  the Weyl chamber $C_0$ in $\frak h(\rho)$ to
some other Weyl chamber $C_1$. But $C_0$ and $C_1$ are conjugate by an
element of the Weyl group of
$Z^0_G(\rho)$, so after a further conjugation we can assume that
conjugation by
$y$ sends $C_0$ to itself. This produces the desired element
$y_\gamma$. Finally, two different choices for $y_\gamma$ differ by an
element of $Z_G(\rho)$ which normalizes $H(\rho)$ and fixes a Weyl
chamber, and hence lies in $H(\rho)$.
\endproof

Thus, we now assume that we have chosen the lift $y_\gamma$ of the outer
automorphism
$\gamma$ to normalize $\frak h(\rho)$ and $C_0$. This means that
$y_\gamma$ acts as a linear automorphism of $\frak h(\rho)$ which
normalizes the set of simple roots $\Delta(\rho)$ associated to $\frak
h(\rho)$ and
$C_0$. Thus, each outer automorphism  $\gamma$ acts as a diagram
automorphism $\sigma_\gamma$ of the Dynkin diagram of $Z_G^0(\rho)$. It
is easy to see  that this in fact defines an action of the group of
outer automorphisms  as a group of diagram automorphisms of the Dynkin
diagram of $Z_G^0(\rho)$. (The point here is that, by the uniqueness
part of Claim 6.8, 
$y_{\gamma_1}y_{\gamma_2}$ is equal to $y_{\gamma_1\gamma_2}$ up to an
element of
$H(\rho)$. Hence conjugation by $y_{\gamma_1}y_{\gamma_2}$ and by
$y_{\gamma_1\gamma_2}$ agree  on
$\frak h(\rho)$ and on the set of roots.) Note that the action of the
group of outer automorphisms on the Dynkin diagram then defines a linear
action on $\frak h(\rho)$, since we can identify the vertices of the
Dynkin diagram with a basis of $\frak h(\rho)$, and this linear action
is the same as the action of $y_\gamma$ defined above. Applying this
discussion in particular to the group $\pi_0(Z_G(\rho))$, which acts as
a group of outer automorphisms of $Z^0_G(\rho)$, we see that there is
an action of
$Z^0_G(\rho)$ on $\frak h(\rho)$.

Our goal for the remainder of this subsection is to prove the following
theorem:

\theorem{6.9} There is a linear isomorphism from $\frak a(\xi)$ to
$\frak h(\rho)$ which is equivariant with respect to the actions of
$\pi_0(Z_G(\rho))$.
\endstatement
\proof We begin by showing that there is a good choice of the principal
nilpotent element $X$ used to define $\xi$, namely one which is
equivariant for the action of $\pi_0(Z_G(\rho))$ in an appropriate
sense.  Fix a Cartan subalgebra $\frak h(\rho)$ for $\frak z_{\frak
g}(\rho)$ and a Weyl chamber $C_0\subset \frak h(\rho)$. Let $H(\rho)$
be the  Cartan subgroup of $Z^0_G(\rho)$ corresponding to $\frak
h(\rho)$. Fix lifts
$y_\gamma\in Z_G(\rho)$ for
$\gamma\in
\pi_0(Z_G(\rho))$ as in the statement of Claim  6.8. Let $\Delta(\rho)$
be the set of simple roots associated to the Weyl chamber $C_0$. Then
the action of $\pi_0(Z_G(\rho))$ on $\frak h(\rho)$ induces an action on
$\Delta(\rho)$. 

\claim{6.10} There exist elements $h_\gamma\in H(\rho)$ and a principal
nilpotent element $X_+$ of $\frak z_{\frak g}(\rho)$ such that, for
every $\gamma \in
\pi_0(Z_G(\rho))$, if we set $y'_\gamma=y_\gamma h_\gamma^{-1}$, then
$X_+$ is invariant under $\operatorname{Ad}y'_\gamma$. 
\endstatement

\proof For each  root $\alpha\in
\Delta(\rho)$, choose a non-zero element $X^\alpha$ in the root space
$\frak z^\alpha$ associated to $\alpha$.  For each $\gamma$  and each
$\alpha\in \Delta$, let $\mu(\alpha,\gamma)$ be the non-zero complex
number with the property that
$$\operatorname{Ad}(y_\gamma)(X^\alpha)=
\mu(\alpha,\gamma)X^{\sigma_\gamma(\alpha)}.$$  Since the $\alpha\in
\Delta(\rho)$ are a linearly independent set of characters on 
$H(\rho)$, it follows that there is an element $h_\gamma\in H(\rho)$
such that for each 
$\alpha\in \Delta(\rho)$ we have $\alpha(h_\gamma)=\mu(\alpha,\gamma)$.
Now replace $y_\gamma$ with the element $y_\gamma'=y_\gamma
h^{-1}_\gamma$. Conjugation by the elements $y_\gamma$ and $y'_\gamma$ 
agree on $\frak h(\rho)$. But in addition  for each $\alpha\in
\Delta(\rho)$ we have
$$\operatorname{Ad}(y'_\gamma)X^\alpha=
\operatorname{Ad}(y_\gamma)(\alpha(h_\gamma^{-1})X^\alpha)= 
\operatorname{Ad}(y_\gamma)(\mu(\alpha,\gamma)^{-1}X^\alpha)=
X^{\sigma_\gamma(\alpha)}.$$   Thus, we see that the element
$$X_+=\sum_{\alpha\in \Delta(\rho)}X^\alpha$$ is invariant under each
of the $y_\gamma'$. Since, for each $\alpha \in \Delta(\rho)$,
$X^\alpha\not= 0$, it follows from
\cite{7, Chap\. 8} that $X_+$ is a principal  nilpotent element. 
\endproof

Replacing the $y_\gamma$ by the $y_\gamma'$, let us recap what we have
managed to establish so far.

\corollary{6.11} Fix a Cartan subalgebra $\frak h(\rho)$ and a Weyl
chamber
$C_0\subset \frak h(\rho)$. Then for each element $\gamma\in
\pi_0(Z_G(\rho))$ there is a lift $y_\gamma$ in $Z_G(\rho)$ such
that\rom:
\roster
\item"{(i)}"  For each $\gamma$ the element $y_\gamma$ normalizes
$\frak h(\rho)$ and
$C_0$.
\item"{(ii)}" Associating to each $\gamma\in \pi_0(Z_G(\rho))$ the
automorphism of $\frak h(\rho)$ given by  conjugation by $y_\gamma$
yields a linear action of $\pi_0(Z_G(\rho))$ on
$\frak h(\rho)$ preserving $C_0$ and hence preserving the set of simple
roots
$\Delta(\rho)$ associated to $\frak h(\rho)$ and $C_0$.
\item"{(iii)}" Viewing the simple roots as a basis for $\frak h(\rho)$,
the linear action  of
$\pi_0(Z_G(\rho))$ on
$\frak h(\rho)$ is compatible with an action of $\pi_0(Z_G(\rho))$ by
diagram automorphisms on the Dynkin diagram for $Z_G^0(\rho)$.
\item"{(iv)}" There is a principal nilpotent element $X_+$ which is
invariant under $\operatorname{Ad}(y_\gamma)$ for each $\gamma$. \qed
\endroster
\endstatement

We return now to the proof of Theorem 6.9. Let $\xi$ be the regular
semisimple
$G$-bundle  determined by the holonomy representation $\rho$ and the
nilpotent element $X_+$ given in Claim 6.10. Then the automorphism
group  $A(\xi)\subseteq Z_G(\rho)$ of $\xi$ is the subgroup of elements
in $Z_G(\rho)$ which centralize $X_+$. In particular, the elements
$y_\gamma$ lie in $ A(\xi)$. Conjugation by the elements $y_\gamma$
normalizes the Lie algebra
$\frak a(\xi)$ of $A(\xi)$ and determines the action of
$A(\xi)/A^1(\xi)=\pi_0(Z_G(\rho))$ on $\frak a(\xi)$. Thus, it will
suffice to exhibit an isomorphism from $\frak a(\xi)$ to $\frak h(\rho)$
which is equivariant with respect to conjugation by the $y_\gamma$ for
all $\gamma
\in \pi_0(Z_G(\rho))$.

Let $h^0\in\frak h(\rho)$ be the unique element contained in the derived
subalgebra $\frak d(\rho)$ which has the property that
$\langle \alpha,h^0\rangle=2$ for  all $\alpha\in
\Delta(\rho)$.  Clearly, $h^0$ is invariant by the action of
$\pi_0(Z_G(\rho))$ on
$\frak h(\rho)$ and hence by conjugation by each of the $y_\gamma$.
Further note that
$[h^0, X_+] = 2X_+$ (and in fact $h^0$ is the unique element of $\frak
h(\rho)$ with this property). By the Jacobson-Morozov theorem, we can
complete
$h^0$ and the principal nilpotent element of $X_+$ of $\frak z_{\frak
g}(\rho)$ to an
$\frak{sl}_2$-triple
$(X_+, h^0, X_-)$. Since $X_+$ is principal, it is easy to see that
$X_-$ is uniquely determined by $X_+$ and $h^0$. Thus, since $X_+$ and
$h^0$ are invariant under conjugation by the
$y_\gamma$ for all $\gamma$, the same is true for $X_-$. Hence the
$\frak{sl}_2$ triple $(X_+, h^0, X_-)$ is invariant under conjugation
by the
$y_\gamma$ for all $\gamma$. Break the Lie algebra $\frak z_{\frak
g}(\rho)$ into  irreducible summands for this action of
$\frak{sl}_2$. Let the rank of $\frak z_{\frak g}(\rho)$ be $s$.
 Then since $X_+$ is  principal there are exactly $s$ irreducible
summands in this decomposition. On the other hand,  by the
classification of irreducible representations of $\frak{sl}_2$,
$\dim
\Ker \ad h^0$ is equal to the number of irreducible summands $V(k)$
which have odd dimension. Since $\frak h(\rho) \subseteq \Ker
\ad h^0$ and $\frak h(\rho)$ is of dimension at least as large as $\Ker
\ad h^0$,  we  see that every irreducible summand must have odd
dimension and so be of the  form $V(2d-1)$ for some integer $d$, and
that $\Ker
\ad h^0 = \frak h(\rho)$. Now $(\ad X_-)^d\:  V(2d-1) \to V(2d-1)$
defines an isomorphism from $\Ker \ad X_+$ on
$V(2d-1)$ to $\Ker \ad h^0$ on
$V(2d-1)$. Thus, given the $\frak{sl}_2$-triple
$(X_+, h^0, X_-)$, there is a canonical linear isomorphism from the
centralizer
$\frak a(\xi)$ of $X_+$ to $\frak h(\rho)$, and since the
$\frak{sl}_2$-triple
$(X_+, h^0, X_-)$ is invariant under conjugation by the $y_\gamma$ for
all
$\gamma$,  it is easy to check that this isomorphism is equivariant with
respect to the conjugations by the $y_\gamma$ for all $\gamma$. Thus we
have constructed a  $\pi_0(Z_G(\rho))$-equivariant isomorphism from 
$\frak a(\xi)$  to $\frak h(\rho)$, as claimed in the statement of
Theorem 6.9.
\endproof

\ssection{6.4. The deformation space of a regular bundle.}

Our goal now will be to relate the deformation theory of a regular
bundle
$\xi$ to the local structure of the moduli space at $[\xi]$. There is a
universal family of semistable $G$-bundles over  a small neighborhood
$U$ of the origin in
$H^1(E; \ad \xi)$. Note that the finite
group $I=A(\xi)/A^1(\xi)$ acts on $H^1(E; \ad \xi)$ and, choosing $U$
to be $I$-invariant, the action lifts to an action on the total space
of this family. The   family induces an $I$-invariant  map from
$U$ to the moduli space
$\Cal M(G)$. Our aim now is to show that, under this map, $U/I$
is identified with a neighborhood of the S-equivalence class 
$[\xi]$ of $\xi$, viewed as  a point of $\Cal M(G)$. This result
does not follow in general from the \'etale slice theorem, since the
stabilizer of a point corresponding to $\xi$ in a GIT space may  fail
to be reductive.

\theorem{6.12} Let $\xi$ be a regular semistable bundle and let
$I=A(\xi)/A^1(\xi)$. Let $U$ be a sufficiently small
$I$-invariant neighborhood of the origin in
$H^1(E; \ad \xi)$ with a universal family as constructed above. Then the
induced map
$U \to \Cal M(G)$ defines an analytic isomorphism of $U/I$ onto a
neighborhood of
$[\xi]$ in ${\Cal M}(G)$.
\endstatement
\proof Let $\xi_0$ be the Yang-Mills bundle S-equivalent to $\xi$ and
let $\rho$ be the corresponding holonomy representation. We choose
normalized transition functions $\{h_{ij}\exp(f_{ij}X)\}$ for
$\xi$, where
$X\in\frak z_{\frak g}(\rho)$ is a regular nilpotent element. The fixed
$1$-cocycle $\{f_{ij}\}$ gives us an identification of $H^1(E;
\scrO_E)$ with $\Cee$. As we have seen in Lemma
3.4, $H^1(E;\ad \xi_0 )=H^1(E;{\Cal O}_E)\otimes_\Cee \frak z_{\frak
g}(\rho)\cong \frak z_{\frak g}(\rho)$ and the Lie algebra
$\frak z_{\frak g}(\rho)$ parametrizes a locally semiuniversal family
$\Xi_0$ of semistable bundles whose central member is $\xi_0$. The
parametrization is given in terms of transition functions (not
necessarily normalized) by
$$Z\mapsto \{h_{ij}\exp (f_{ij}Z)\}.$$
By Theorem 5.2 and Corollary 5.3, a local model for the moduli space
near
$[\xi]=[\xi_0]$  is given by
$$\frak z_{\frak g}(\rho)/\!\!/Z_G(\rho).$$ This GIT quotient is
identified with
$(\frak h(\rho)/W(\rho))/I$, where ${\frak h}(\rho)$ is a Cartan
subalgebra for $\frak z_{\frak g}(\rho)$ and $W(\rho)$ is the Weyl
group.

By Theorem 4.6,
$$H^1(E;\ad\xi)\cong H^1(E;{\Cal O}_E )\otimes
(\frak z_{\frak g}(\rho)/\operatorname{Im}\ad X)\cong \frak z_{\frak
g}(\rho)/\operatorname{Im}\ad X.$$
The tangent space at $X$ to the orbit
of $X$ in $\frak z_{\frak g}(\rho)$ under the adjoint action is the
image of $\ad  X$. Consider a linear  slice
$N$ in $\frak z_{\frak g}(\rho)$  normal to the tangent space at $X$ to
the orbit of
$X$. Then $N$  parametrizes a family $\Xi$ of
bundles: the point
$Y\in N$ determines the bundle with (not necessarily normalized)
transition functions
$$\{h_{ij}\exp (f_{ij}(X+Y))\}.$$
Note that the origin $Y=0$ of $N$ corresponds to the bundle $\xi$.
The Kodaira-Spencer map of $\Xi$ at the central point
$\xi$ defines  a linear isomorphism from $N$ onto $H^1(E;\ad\xi)$, and
thus $\Xi$ is a locally semiuniversal deformation of $\xi$. 

As we saw in Theorem 4.1, the group $A(\xi)$ is identified with the
subgroup of
$Z_G(\rho)$ which
fixes the point $X$. Since $A(\xi)$ is an extension of a finite group by
a commutative algebraic group, it has a Levi factor
$R(\xi)$ which is a reductive subgroup, in fact  an extension of a
finite group by an algebraic torus. The subgroup
$R(\xi)$ maps onto the component group of $A(\xi)$.  We can choose the
linear slice $N$ to be invariant under $R(\xi)$. Clearly, this action
lifts to an action on the universal family $\Xi$ constructed above. By
Lemma 6.7(ii), the intersection of
$R(\xi)$ with
$A^1(\xi)$ is contained in the center of
$Z_G^0(\rho)$ and hence acts trivially on $\frak z_{\frak g}(\rho)$ and
on $\Xi$. Thus, the action of $R(\xi)$ on the slice $N$ factors through
the finite quotient
$I=A(\xi)/A^1(\xi)$. For this choice of slice $N$, the differential of
the Kodaira-Spencer map from $N$ to $H^1(E;\ad \xi )$ is a linear
isomorphism which is  equivariant with respect to the actions of 
$A(\xi)/A^1(\xi)$.

Now let us consider the map from the family $N$ to the moduli space 
${\Cal M}(G)$.  By the remarks at the beginning of the proof, the
restriction of the universal family $\Xi_0$ for $\frak z_{\frak
g}(\rho)$ to the affine slice $\{X\}+N$ is the universal family $\Xi$
for $\xi$. The induced map from
$N$ to $\Cal M(G)$ defined by $\Xi$ is then
induced by the inclusion of $\{X\}+N$ into $\frak z_{\frak g}(\rho)$
followed first by the natural map to the quotient $\frak z_{\frak
g}(\rho)/\!\!/Z_G^0(\rho)$, and then by the induced map to $\Cal M(G)$.

\claim{6.13}  The differential at $0\in N$ of the composition
$$N \to \{X\}+N\to \frak z_{\frak g}(\rho)\to \frak z_{\frak
g}(\rho)/\!\!/Z_G^0(\rho)$$ is an isomorphism.
\endstatement

\proof Quite generally, let $G$ be a connected, complex reductive group
with maximal torus $H$. Then the GIT adjoint quotient of $G$ is
identified with
$H/W$, the quotient of $H$  by the Weyl group. This quotient is a
smooth variety of dimension equal to the rank of $G$. According to a
theorem of Steinberg (see  Section 4.20  of \cite{20}), the differential
of the adjoint quotient map $G\to H/W$ at any regular element $x\in G$
is surjective. Now let $X$ be a regular nilpotent element in the Lie
algebra $\frak g$ of $G$. Then $x=\exp  X$ is a regular unipotent
element in $G$. Furthermore, since $X$ is conjugate to an element in a 
small neighborhood of the origin, we see that the differential of the
exponential mapping at $X$ is an isomorphism. Let $\frak h$ be the
Cartan subalgebra of $\frak g$ tangent to $H$. The map from the adjoint
quotient of $\frak g$ to the adjoint quotient of $G$ is identified with
the natural map $\frak h/W\to H/W$ and hence it is a local
diffeomorphism at the origin.
 Thus, it follows that the differential of the adjoint quotient mapping
$\frak g\to \frak h/W$ is surjective at $X$. Of course, the kernel of
this differential contains the tangent to the orbit space through $X$.
Since $X$ is regular, the normal to the  orbit is the same dimension as
$\frak h/W$, and hence we conclude that the restriction of this
differential to the normal space is an isomorphism. Applying the above
to the case
$G=Z_G^0(\rho)$ establishes the claim.
\endproof

The inclusion $\{X\}+N\subset \frak z_{\frak g}(\rho)$ is
equivariant with respect to the homomorphism
$$I=A(\xi)/A^1(\xi)\to \pi_0(Z_G(\rho))$$ 
 induced by the inclusion of $A(\xi)\subseteq Z_G(\rho)$. There is thus
an equivariant isomorphism from
$H^1(E;\ad \xi )$ with its  $I$-action to 
$\frak z_{\frak g}(\rho)/\!\!/Z_G^0(\rho)=\frak h(\rho)/W(\rho)$ with
its 
$\pi_0(Z_G(\rho))$-action. We furthermore have a commutative diagram
$$\CD
H^1(E;\ad \xi )/I @>>> (\frak z_{\frak g}(\rho)/\!\!/Z_G^0(\rho))/I\\
@VVV @VVV\\
\Cal M(G) @= \Cal M(G).
\endCD$$
Theorem 6.12 now follows from Theorem 5.2 and Corollary 5.3.
\endproof

\section{7. The description of regular bundles in the classical cases.}

In this section, we shall explicitly describe the regular $G$-bundles
when $G$ is one of the classical groups, by elementary methods. The
cases $G=GL_n$ and
$G=SL_n$ were described in detail in Section 1. In this section, we
describe the regular bundles corresponding to the groups
$Sp(2n)$ and $SO(n)$. Here, the study of $SO(n)$-bundles divides into
two cases,  depending on whether or not the bundle is liftable to
$Spin(n)$. However, in case the bundle is liftable, we do not
explicitly analyze the set of liftings. Finally, we also consider
unliftable bundles on the non-simply connected groups $SO(2n)/\{\pm
\Id\}$ and $Sp(2n)/\{\pm \Id\}$. In this case, it is better to work
with conformal bundles, i\.e\. with vector bundles together with a
nondegenerate form with values in a line bundle. The set of such
bundles corresponding to regular bundles is described at the end of the
section.

\ssection{7.1. Existence of nondegenerate pairings.}

We are interested in the case where $V$ is a semistable rank $n$ bundle
with a nondegenerate form (with values in $\scrO_E$), which is either
symmetric or alternating. Such a form defines an isomorphism from $V$
to $V\spcheck$. We have the following straightforward result for
regular vector bundles $V$:

\proposition{7.1} Suppose that $V=\bigoplus _iI_{d_i}(\lambda _i)$ is a
regular rank $n$ semistable vector bundle with $\det V = \scrO_E$ such
that
$V\cong V\spcheck$. Then the
$\lambda _i$ appearing as the Jordan-H\"older constituents of $V$ have
the following property. If $\lambda _i \neq \lambda _i^{-1}$, then
there is a unique $j$ such that $\lambda _j = \lambda _i^{-1}$, and in
this case $d_i = d_j$. If $\lambda _i=\lambda _i^{-1}$, so that
$\lambda _i$ is a line bundle  of order $2$, then every isomorphism
$V\to V\spcheck$ restricts to an isomorphism
$I_{d_i}(\lambda _i) \to I_{d_i}(\lambda _i)\spcheck = I_{d_i}(\lambda
_i)$. Finally, suppose that $\lambda _i$ is a line bundle  of order $2$
and that $d_i$ is odd, so that in particular $d_i\neq 0$. Then for
every $j$ such that $\lambda _j$ is a nontrivial line bundle  of order
$2$, $d_j$ is odd as well. \qed
\endstatement

Next we turn to existence of nondegenerate symmetric or alternating
forms on vector bundles, and to the corresponding question of when such
forms give regular symplectic or orthogonal bundles. Suppose that 
$V$ is either an orthogonal or a symplectic bundle. Writing
$V= \bigoplus _iI_{d_i}(\lambda _i)$, we see that if $\lambda_i \neq
\lambda_i ^{-1}$ we can pair up $I_{d_i}(\lambda _i)$ with a dual space
$I_{d_i}(\lambda _i^{-1})$. Moreover, the only symmetric or alternating
forms on
$I_{d_i}(\lambda _i)\oplus I_{d_i}(\lambda _i)^{-1}$ are given as
follows: take an isomorphism $\varphi\: I_{d_i}(\lambda _i)^{-1}\to 
I_{d_i}(\lambda _i)\spcheck $ and define
$$\langle s_1, s_2\rangle = \varphi (s_2)s_1\pm {}^t\varphi
(s_1)(s_2),$$ with the sign chosen according to whether we want the
form to be symmetric or alternating. Note that every two such forms are
then conjugate under the action of $\Aut (I_{d_i}(\lambda _i)\oplus
I_{d_i}(\lambda _i)^{-1})$, and the group of orthogonal or symplectic
automorphisms is just a copy of $\Aut I_{d_i}(\lambda _i) \cong \Aut
I_{d_i}$. In particular it is a connected abelian group of dimension
$d_i$. Thus the crucial case to work out is the case of
$I_d(\lambda)$ where $\lambda$ is a line bundle of order two, and in
fact it is clearly enough to work out the case of $I_n$ itself:

\theorem{7.2} \roster
\item"{(i)}" There exists a nondegenerate alternating form on $I_n$,
i\.e\. a nondegenerate skew symmetric pairing $I_n\otimes I_n \to
\scrO_E$, if and only if $n$ is even. Every two such pairings are
conjugate under the action of $\Aut I_n$.
\item"{(ii)}" There exists a nondegenerate symmetric form on $I_n$,
i\.e\. a nondegenerate symmetric pairing $I_n\otimes I_n \to \scrO_E$,
if and only if
$n$ is odd. In this case every two such pairings are conjugate under
the action of
$\Aut I_n$. 
\endroster
\endproclaim
\demo{Proof} Fix a homomorphism $\varphi\: I_d \to I_d\spcheck \cong
I_d$. Given the canonical filtration $F^i$ on
$I_d$, we have the dual decreasing filtration $(F^i)^{\perp}$ of
$I_d\spcheck$, as well as the canonical increasing filtration
$(F^i)\spcheck$ of
$I_d\spcheck$. Since both filtrations are canonical, we must have
$(F^i)^{\perp} = (F^{d-i})\spcheck$. Now $\varphi$ must preserve the
filtrations, and so
$\varphi (F_i)\subseteq (F^{d-i})^\perp$. It follows that $\varphi$
induces a well-defined pairing
$$\bar{\varphi _i}\: \left(F^i/F^{i-1}\right) \otimes
\left(F^{d-i+1}/F^{d-i}\right)\to \scrO_E.$$ Either $\varphi$ is an
isomorphism or it factors through some homomorphism
$I_d \to I_{d-k} \cong (F^{d-k})\spcheck = (F^k)^\perp$ for some $k>0$.
In the second case, $\varphi (F^i)
\subseteq (F^{i-k})\spcheck = (F^{d-i+k})^\perp$, and hence, for the
induced pairing on $I_d$, $F^i$ is annihilated by $F^{d-i+1}$ But then
$\bar{\varphi _i} = 0$ for every $i$. Conversely, if $\varphi$ is an
isomorphism, then so is $\bar{\varphi _i}$ for every $i$. We see that
$\varphi$  is an isomorphism if and only if $\bar{\varphi _i}$ is an
isomorphism for every
$i$, if and only if there exists an $i$ such that $\bar{\varphi _i}\neq
0$.

Now given $\varphi\: I_d \to I_d\spcheck$, we can take the transpose
${}^t\varphi\: I_d \to I_d\spcheck$. Then $\varphi$ is symmetric if 
$\varphi = {}^t\varphi$, i\.e\. if $\varphi - {}^t\varphi =0$, and it is
alternating if $\varphi + {}^t\varphi = 0$. Note that, given $\varphi$,
we can consider $\overline{{}^t\varphi}_i = {}^t\bar \varphi _{d-i+1}$.
First consider the case where
$d = 2n+1$ is odd. In this case taking $i = n+1$ we have
$$(\overline{\varphi + {}^t\varphi})_{n+1} = \bar{\varphi}_{n+1} +
{}^t\bar{\varphi}_{n+1}.$$ But $\bar{\varphi}_{n+1}$ is a pairing
$(F^{n+1}/F^n)\otimes (F^{n+1}/F^n) \to
\scrO_E$, and is thus equal to its transpose. Thus, if $\varphi$ is an
arbitrary isomorphism, then $\varphi + {}^t\varphi$ is an isomorphism
as well, since $(\overline{\varphi + {}^t\varphi})_{n+1} =
2\bar{\varphi}_{n+1}\neq 0$. It follows that $\varphi + {}^t\varphi$ is
a symmetric isomorphism. Moreover,  if we begin with a symmetric
$\varphi$, then $\frac12(\varphi + {}^t\varphi) =
\varphi$, and so every symmetric isomorphism arises from this
construction. Since every two $\varphi$ are conjugate under the
multiplication action of $\Aut I_{2n+1}$, every two symmetric forms
$\varphi$ are conjugate under the corresponding action of $\Aut
I_{2n+1}$. Clearly there are no nondegenerate alternating forms in this
case, since $d = 2n+1$ is odd. We could also see this as follows: if
$\varphi$ is alternating and nondegenerate, then  so is $\varphi -
{}^t\varphi = 2\varphi$. But restricting to $F^{n+1}/F^n$ and using the
fact that $\bar{\varphi}_{n+1} = {}^t\bar{\varphi}_{n+1}$, we find
instead that $(\overline{\varphi - {}^t\varphi})_{n+1} = 0$, so that
$\varphi$ must in fact be degenerate.

Now assume that $d=2n$ is even. We claim that in this case there are no
symmetric nondegenerate forms. Indeed, suppose that $\varphi\: I_{2n}
\to I_{2n}\spcheck$ is a symmetric homomorphism. Then $\varphi$ induces
a pairing
$$\tilde{\varphi}\: (F^{n+1}/F^{n-1})\otimes (F^{n+1}/F^{n-1}) \to
\scrO_E,$$ whose associated diagonal components are just $\bar\varphi
_{n+1}$ and
$\bar\varphi _n$. We claim that $\tilde{\varphi}$ is degenerate, or in
other words that every symmetric form on $I_2$ is degenerate. Now on
$\Cee ^2$ with the standard nondegenerate form, there are exactly two
isotropic lines. Thus if there were a nondegenerate form on $I_2$, the
isotropic subbundle $F^1$ would have a unique isotropic complement, and
thus $I_2$ would be a direct sum of line bundles. This is a
contradiction. Thus $\varphi$ is degenerate.

Now begin with an arbitrary isomorphism $\varphi\: I_{2n} \to I_{2n}$.
We claim that $\varphi - {}^t\varphi$ is an isomorphism as well, thus
giving a nondegenerate alternating form on $I_{2n}$. Indeed, we must
have $\varphi + {}^t\varphi$ degenerate, and thus $(\overline{\varphi +
{}^t\varphi})_i=0$ for every $i$. It follows that $\bar{\varphi}_i = -
{}^t\bar{\varphi}_{2n-i}$ for every $i$, and thus that
$(\overline{\varphi - {}^t\varphi})_i$ is nonzero for every $i$. We
conclude that $\varphi -  {}^t\varphi$ defines a nondegenerate
alternating form on $I_{2n}$. As in the symmetric case, if $\varphi$ is
alternating to start with, then $\varphi -  {}^t\varphi = 2\varphi$, so
that all alternating forms arise in this way, and every two such are
conjugate under the automorphisms of $I_{2n}$. This concludes the proof
of the theorem.
\endproof

If $V$ is a vector bundle with an alternating form $Q$, we denote by 
$\Aut^QV$ the automorphisms of $V$ which preserve the form; these
automatically have determinant $1$ If instead $Q$ is symmetric, we
define $\Aut^QV$ similarly, but also impose the requirement that the
determinant be $1$. 

\proposition{7.3} \roster
\item"{(i)}" Let $I_{2n}$ be given a nondegenerate alternating form,
and equip
$F^{2n-1}/F^1$ with the induced nondegenerate form. Then the map
$$\Aut^QI_{2n} \to \Aut^Q(F^{2n-1}/F^1)$$ is surjective, and its kernel
is isomorphic to $\Cee$.  Moreover $\dim
\Aut^QI_{2n} = n$.
\item"{(ii)}" Let $I_{2n+1}$ be given a nondegenerate symmetric form,
and give $F^{2n}/F^1$  the induced nondegenerate form. Then the map
$$\Aut^QI_{2n+1} \to \Aut^Q(F^{2n}/F^1)$$ is surjective, and its kernel
is isomorphic to $\Cee$. Moreover $\dim
\Aut^QI_{2n+1} = n$.
\item"{(iii)}" For every $n\geq 1$, 
$$\gather
\dim H^0(\bigwedge^2 I_{2n}) = \dim H^0(\Sym^2 I_{2n}) = \dim
H^0(\bigwedge^2 I_{2n+1})  =n;\\
\dim H^0(\Sym^2 I_{2n+1})
 = n+1.
\endgather$$
\endroster
\endstatement
\demo{Proof} Consider first the symplectic case. We argue by induction
on $n$. In case $n=1$, $\Aut^QI_2$ is the same as the elements of
determinant one in $\Cee[T]/(T^2)$, which is isomorphic to $\{\pm
1\}\times
\Cee$. In general, we always have an exact sequence
$$\{1\} \to T^{2n-2}\Cee [T]/(T^{2n}) \to \Aut I_{2n} \to \Aut
(F^{2n-1}/F^1)
\to \{1\}.$$ The image of $\Aut^QI_{2n}$ in $\Aut (F^{2n-1}/F^1)$
clearly lies in $\Aut^Q(F^{2n-1}/F^1)$. If we denote the alternating
form on
$I_{2n}$ by $\langle \cdot, \cdot
\rangle$, and let $e_1$ be a nowhere vanishing section of $F^1$, then
the Picard-Lefschetz reflection
$$s\mapsto s+ t\langle s, e_0\rangle e_0$$ is a symplectic automorphism
of
$I_{2d}$, and a calculation with a local basis shows that these are the
unique symplectic elements in $\Ker(\Aut I_{2n} \to \Aut
(F^{2n-1}/F^1))$. More intrinsically, $\Ker(\Aut I_{2n} \to \Aut
(F^{2n-1}/F^1))$ is identified with $\Hom (F^{2n}, F^1)$ via
$$\alpha \in \Hom (F^{2n}, F^1) \mapsto (s\mapsto s+\alpha (s)).$$ Thus
$\dim \Aut^QI_{2n} \leq n$, and equality holds only if the map
$$\Aut^QI_{2n} \to \Aut^Q(F^{2n-1}/F^1)$$  is surjective. (Note that by
induction $\Aut^Q(F^{2n-1}/F^1)$ is a vector group times
$\{\pm 1\}$ and $-\Id \in \Aut^QI_{2n}$ maps onto the $-1\in \{\pm
1\}$.)

The set of nondegenerate symplectic forms in  $H^0(\bigwedge^2 I_{2n})$
is an open dense subset of $H^0(\bigwedge^2 I_{2n})$, and $\Aut I_{2n}$
acts transitively on this set, with stabilizer $\Aut^QI_{2n}$. Since
$\dim
\Aut I_{2n} = 2n$, we see that  $h^0(\bigwedge^2 I_{2n}) = 2n - \dim
\Aut^QI_{2n} \geq n$, with equality if and only if
$\dim\Aut^QI_{2n} = n$. We further claim that, if $A$ is a form in
$H^0(\bigwedge^2 I_{2n+1})$, then $F_1$ is in the radical of $A$, and
so there is an induced map $H^0(\bigwedge^2 I_{2n+1}) \to
H^0(\bigwedge^2 I_{2n})$ which is clearly an isomorphism. Indeed, given
any homomorphism $\varphi\: I_d\to I_d\spcheck$, $\Ker \varphi$ is a
subbundle of
$I_d$ and so has degree $\leq 0$. On the other hand $\operatorname{Im}
\varphi
$ has degree $=-\deg \Ker \varphi$, and also has degree zero since it
injects into the semistable bundle $I_d\spcheck$. It follows that $\deg
\Ker \varphi =0$ and thus that $\Ker \varphi = F^i$ for some $i$. In
particular, if $\Ker
\varphi \neq 0$, then it contains $F^1$. Applying this remark to the
necessarily degenerate form $A$, we see that the radical of $A$
contains $F^1$. In conclusion, then, we can at least say that
$h^0(\bigwedge^2 I_{2n+1}) = h^0(\bigwedge^2 I_{2n}) \geq n$.

Next we consider the orthogonal case. For $n=0$, the orthogonal
automorphisms of
$\scrO_E$  of determinant one  are just the identity. For $n > 0$, we
can again consider the kernel of the natural map $\Aut^QI_{2n+1} \to
\Aut^Q(F^{2n}/F^1)$. We claim that again this kernel is isomorphic to
$\Cee$. In the orthogonal case, the analogue of the Picard-Lefschetz
reflections giving unipotent orthogonal automorphisms are the
Eichler-Siegel transformations. In our case, suppose that $\psi \in \Ker
(\Aut^QI_{2n+1} \to \Aut^Q(F^{2n}/F^1)$. Then $(\psi -\Id)(F^{2n})
\subseteq F^1$, and we can write $\psi (v) = v+ \alpha (v)e_0$ for a
well-defined homomorphism $\alpha \: F^{2n} \to \scrO_E$. A local
calculation shows that, given $\alpha$, there is a unique way to define
an orthogonal
$\psi$ with the given retsriction to $F^{2n}$. In this way, we have
identified
$\Ker (\Aut^QI_{2n+1} \to \Aut^Q(F^{2n}/F^1)$ with $\Hom (F^{2n}, F^1)
\cong \Cee$.  Note that all such have determinant one. It follows that
$\Aut^QI_{2n+1}$ is a vector group of dimension at most $n$, with
equality only if $\Aut^QI_{2n+1} \to \Aut^Q(F^{2n}/F^1)$ is surjective.

Arguing as in the symplectic case, we see that $h^0(\Sym^2 I_{2n+1})
\geq n+1$, with equality if and only if $\Aut^QI_{2n+1} \to 
\Aut^Q(F^{2n}/F^1)$ is surjective, and that there is a natural
isomorphism $H^0(\Sym^2 I_{2n+2}) \to H^0(\Sym^2 I_{2n+1})$. Thus for
all
$n\geq 0$, $h^0(\Sym^2 I_{2n+2}) \geq n+1$, with equality if and only
if $\Aut^QI_{2n+1} \to 
\Aut^Q(F^{2n}/F^1)$ is surjective. On the other hand,
$$2n = h^0(I_{2n} \otimes I_{2n}) = h^0(\bigwedge^2 I_{2n}) + h^0(\Sym^2
I_{2n}) \geq n+n = 2n,$$ so that equality must hold for all $n$. We
have thus proved all of the statements of (7.3).
\endproof

To handle orthogonal bundles in the even rank case whose Jordan-H\"older
constituents are all $\scrO_E$, we must use the bundle $I_{2d-1}\oplus
\scrO_E$. The next lemma shows that this bundle has the right number of
automorphisms:

\lemma{7.4} Suppose that $n > 1$. Let $I_{2n-1}\oplus \scrO_E$ have a
nondegenerate symmetric form given by choosing a nondegenerate
symmetric form on $I_{2n-1}$ and taking the orthogonal direct sum with
a  nondegenerate form on $\scrO_E$. Then $\Aut^Q( I_{2n-1} \oplus
\scrO_E)$ is abelian of dimension $n$. In fact, the identity component
of $\Aut^Q( I_{2n-1}
\oplus \scrO_E)$ is $(\Aut^Q I_{2n-1})\times \Cee$.
\endstatement
\proof There is a natural inclusion of $\Aut^Q I_{2n-1}$ in
$\Aut^Q( I_{2n-1} \oplus \scrO_E)$. Now let $\psi \in  
\Aut^Q( I_{2n-1} \oplus \scrO_E)$. we can write
$$\psi (s, \lambda) = (a\lambda e_0 + Ts, t\alpha (s) + c\lambda),$$
where $T\in \Aut I_{2n-1}$ and $\alpha \: I_{2n-1} \to \scrO_E$ is a
homomorphism, well-defined up to scalars, which we may as well take to
be
$\langle s, e_0\rangle$. The condition
$$s^2 + \lambda ^2 = \psi (s, \lambda)^2$$ for all $s,\lambda$ implies
that $c= \pm 1$ (taking $s=0$).  We have (setting 
$\lambda = 0$)
$$(Ts)^2 + t^2\alpha (s)^2 = s^2,$$ and the remaining condition is
$$2a\alpha (Ts) + 2tc\alpha (s) = 0.$$  If $s\in F^{2n-2}$, then
$\alpha (s) = 0$ and thus $(Ts)^2 = s^2$. Now $\Aut^Q I_{2n-1}$ maps
onto $\Aut^Q (F^{2n-2}/F^1)$, and so after modifying by an element of
$\Aut^Q I_{2n-1}$ we can assume that $\pm Ts = s + b\alpha (s)$ for all
sections $s$ of $I_{2n-1}$ and thus that $\alpha (Ts) = \alpha (s)$. A
calculation shows that, if we require that the determinant be one, the
unique choices are
$$\align
\psi _t^+(s, \lambda)  &= (s-(t\lambda + \frac{t^2}2\langle s,
e_0\rangle)e_0, t\langle s, e_0\rangle + \lambda);\\
\psi _t^-(s, \lambda)  &= (-s+(t\lambda + \frac{t^2}2\langle s,
e_0\rangle)e_0, -t\langle s, e_0\rangle - \lambda).
\endalign$$ Note that $\psi _t^- = \psi _t^+\circ (-\Id) = (-\Id) \circ
\psi _t^+$. Finally, to see that $\Aut^Q( I_{2n-1}\oplus \scrO_E)$ is
abelian, it suffices to show, for all $A\in \Aut^Q I_{2n-1}$, that
$A\circ \psi ^+_t = \psi ^+_t\circ A$ and likewise for $\psi ^-_t$. A
calculation shows that it is equivalent to show that, for all $A\in
\Aut^Q I_{2n-1}$ and all sections $s$ of $I_{2n-1}$, $Ae_0 = e_0$ and
$\langle As, e_0 \rangle = \langle s, e_0 \rangle$. Since $A$ is
orthogonal it will suffice to check that $Ae_0 = e_0$. If $Ae_0 =
ce_0$, then $c^{-1}$ is also an eigenvalue of $A$, and since $I_{2n-1}$
is regular the only possibility is $c=\pm 1$ and that $A-c\Id$ is
nilpotent. Taking determinants we see that
$c=1$.  \endproof

Of course, in case $n=1$, the bundle $\scrO_E\oplus \scrO_E$ with the
induced diagonal form has automorphism group $SO(2)$, which has
dimension one and is abelian again.

\ssection{7.2. Description of regular symplectic and orthogonal
bundles.}

We can now describe the bundles with minimal automorphism group.

\proposition{7.5} Let $V$ be a symplectic bundle of rank $2n$ which is a
semistable bundle with trivial determinant. Then $\dim \Aut ^QV\geq n$.
Moreover equality holds if and only if $V$ is regular as a vector
bundle, and in this case the symplectic form on
$V$ is unique up to conjugation by $\Aut V$.
\endstatement
\proof We can write $V = \bigoplus _iI_{d_i}(\lambda _i)$, where if
$\lambda _i
\neq \lambda _i^{-1}$ the summand $I_{d_i}(\lambda _i)$ must be paired
with a summand $I_{d_i}(\lambda _i^{-1})$, whereas if $\lambda _i =
\lambda _i^{-1}$ this is not necessarily the case. Note that, for
$\lambda _i
\neq \lambda _i^{-1}$, if there are two summands $I_{d_i}(\lambda _i)$
and
$I_{d_j}(\lambda _i)$, then 
$$\Aut^Q(I_{d_i}(\lambda _i)\oplus I_{d_j}(\lambda _i) \oplus
I_{d_i}(\lambda _i^{-1}) \oplus I_{d_j}(\lambda _i^{-1})$$ at least
includes the group $\Aut (I_{d_i}(\lambda _i)\oplus I_{d_j}(\lambda
_i))$, which has dimension $\geq d_i + d_j$, with equality only if one
of $d_i, d_j$ is zero. Now consider the case of a $\lambda _i$ of order
$2$. If
$I_{d_i}(\lambda _i)$ is paired with another direct summand isomorphic
to
$I_{d_i}(\lambda _i)$, then $\Aut^Q(I_{d_i}(\lambda _i)\oplus
I_{d_i}(\lambda _i))$ contains $\Aut I_{d_i}(\lambda _i)$ acting on the
factors by $(\varphi, {}^t\varphi)$. But it also contains an upper
triangular group of the form
$$(s_1, s_2) \mapsto (s_1 + \varphi (s_2), s_2),$$ as long as $\varphi$
is symmetric. As we have seen, $h^0(\Sym ^2I_{d_i}) \neq 0$ as long as
$d_i > 0$. Thus in this case $\dim
\Aut^Q(I_{d_i}(\lambda _i)\oplus I_{d_i}(\lambda _i)) > d_i$. 

In case the form pairs $I_{d_i}(\lambda_i)$ with itself, we must have
$d_i= 2d'_i$ even and $\dim \Aut^QI_{2d'_i}(\lambda _i) =d'_i$. Suppose
now that $V$ contains two summands $I_{2d'_i}(\lambda _i)$ and
$I_{2d'_j}(\lambda _i)$, each paired with itself by the symplectic
form. Then
$\Aut^Q(I_{2d'_i}(\lambda _i)\oplus I_{2d'_j}(\lambda _i))$ contains
$\Aut^Q(I_{2d'_i}(\lambda _i))\times  \Aut^Q(I_{d'_j}(\lambda _i))$ and
so has dimension at least $d'_i + d'_j$. However, if both $d'_i$ and
$d'_j$ are positive the dimension of the automorphism group must be
larger. It suffices to consider the case $\lambda _i = \scrO_E$. 
Choose a section $e_0$ of $I_{2d'_i}$ and a section $f_0$ of
$I_{2d'_j}$. The for all $t$ we have the symplectic automorphism
$$(s_1, s_2) \mapsto (s_1 + t\langle s_2, f_0\rangle e_0, s_2 +
t\langle s_1, e_0\rangle f_0).$$ So $\dim \Aut^Q(I_{2d'_i}(\lambda
_i))\times  \Aut^Q (I_{d'_j}(\lambda _i)) > d'_i + d'_j$ unless there
is just one factor. We see that we have shown that $\dim \Aut ^QV\geq
n$, with  equality holding if and only if $V$ is regular. In this case,
it follows from  the discussion of the cases
$I_{d_i}(\lambda _i)\oplus I_{d_i}(\lambda _i)^{-1})$, $\lambda _i \neq
\lambda _i^{-1}$ prior to (7.2) and (7.2) (i) that a nondegenerate 
symplectic form on
$V$ exists and is unique up to conjugation by $\Aut V$.
\endproof

\definition{Definition 7.6} A semistable symplectic bundle $V$ of rank
$2n$ with trivial determinant such that
$\dim \Aut ^QV = n$ will be called a {\sl regular symplectic bundle}.
Note that $V$ is regular if and only if the underlying vector bundle is
regular.
\enddefinition

In the orthogonal case, we have similar but more complicated results.
Let $\eta _1, \eta _2, \eta _3$ denote the three distinct line bundles
of order $2$ on
$E$. 

\proposition{7.7} Let $V$ be a semistable vector bundle with trivial
determinant and a nondegenerate symmetric form. 
\roster
\item"{(i)}" If no Jordan-H\"older constituent of $V$ is a line bundle
of order
$2$, then necessarily $V$ has even rank $2n$. In this case,  $\dim \Aut
^QV\geq n$. Moreover equality holds if and only if $V$ is regular, and
in this case the orthogonal form on
$V$ is unique up to conjugation by $\Aut V$.
\item"{(ii)}" Suppose given $\lambda _i \neq \lambda _i^{-1}, 1\leq
i\leq k$, positive integers $d_i$, and nonnegative integers $a_1,a_2,
a_3$. Consider the bundle
$$V = \left(\bigoplus _{i=1}^kI_{d_i}(\lambda _i)\oplus I_{d_i}(\lambda
_i^{-1})\right) \oplus I_{2a_1+1}(\eta _1) \oplus I_{2a_2+1}(\eta _2) 
\oplus I_{2a_3+1}(\eta _3) .$$ Then $V$ is a regular vector bundle of
rank $2n +1$, where $n = \sum _id_i + a_1+a_2+a_3 +1$, and $V$ carries
a nondegenerate orthogonal form. Moreover every two such forms on $V$
are conjugate under $\Aut V$. Finally, $\dim \Aut ^QV = n-1$. If
conversely $V$ is an $SO(2n+1)$-bundle, then $\dim \Aut ^QV \geq n-1$,
and equality holds if and only if either $V$ is as described above or
$V$ is of the form $V'\oplus I_{2a+1}\oplus \scrO_E$, where $V'$ is as
described above and $I_{2a+1}\oplus \scrO_E$ has a nondenerate diagonal
form as in Lemma \rom{7.4}.
\item"{(iii)}" Suppose given $\lambda _i \neq \lambda _i^{-1}, 1\leq
i\leq k$, positive integers $d_i$, and nonnegative integers
$a_0,a_1,a_2, a_3$. Consider the bundle
$$V = \left(\bigoplus _{i=1}^kI_{d_i}(\lambda _i)\oplus I_{d_i}(\lambda
_i^{-1})\right) \oplus I_{2a_0+1} \oplus I_{2a_1+1}(\eta _1) \oplus
I_{2a_2+1}(\eta _2)  \oplus I_{2a_3+1}(\eta _3) .$$ Then $V$ is a
regular vector bundle of rank $2n$, where $n = \sum _id_i +
a_0+a_1+a_2+a_3  + 2$, and $V$ carries a nondegenerate orthogonal form.
Moreover every two such forms on $V$ are conjugate under $\Aut V$.
Finally,
$\dim \Aut ^QV = n-2$. If conversely $V$ is an $SO(2n)$-bundle, then
$\dim \Aut ^QV \geq n-2$, and equality holds if and only if $V$ is
regular, if and only if $V$ is as described above.
\endroster
\endstatement
\proof The proofs are very similar to the symplectic case. If $V$
contains a summand $I_{d_i}(\lambda _i)\oplus I_{d_i}(\lambda _i^{-1})$
with $\lambda_i \neq \lambda _i^{-1}$, the argument in the symplectic
case shows that such a summand always contributes a group of dimension
$d_i$ to
$\Aut ^QV$, and two summands  $I_{d_i}(\lambda _i)\oplus I_{d_i}
(\lambda _i^{-1})$ and  $I_{d_j}(\lambda _i)\oplus I_{d_j}(\lambda
_i^{-1})$ with the same $\lambda _i$ will contribute a group of
dimension greater than $d_i+d_j$ unless one of $d_i, d_j$ is zero. If
$I_d(\eta)$ is not paired with itself (here $\eta$ is a line bundle of
order $2$), then $\Aut ^Q(I_d(\eta) \oplus I_d(\eta))$ contains $\Aut
^Q I_d(\eta)$ as well as the  maps $(s_1, s_2) \mapsto (s_1+ \varphi
(s_2), s_2)$, as long as $\varphi$ is alternating. This will force $\Aut
^Q(I_d(\eta) \oplus I_d(\eta))$ to have dimension greater than $e$,
except in case $d=1$ in which case it has dimension exactly $1$. Thus
$I_d(\eta)$ must pair with itself, and hence  $d=2e+1$ is odd.  Keeping
track of the parity of the rank gives the statement of the theorem.
\endproof

Cases (ii) and (iii) above describe the regular unliftable
$SO(k)$-bundles. One difference between the two cases is as follows. In
Case (iii), the automorphism group of $V$ as an $SO(2n)$-bundle is
always abelian. In case (ii), if the bundle $V$ has a factor isomorphic
to $\scrO_E$, then the automorphism group is nonabelian, and
conversely. For example, in the simplest case $V$ has a summand of the
form $\scrO_E\oplus \scrO_E$ as well as the summand $\eta_1\oplus
\eta_2\oplus \eta _3$. In this case, we can embed the automorphism
group $O(2)$ of the factor $\scrO_E\oplus \scrO_E$ into $\Aut^QV$,
using an automorphism of the form
$\pm 1$ on one of the $\eta_i$ factors to make the determinant $1$. The
fact that the automorphism group is not abelian in this case is
reflected in the fact that, as a weighted projective space, the moduli
space is of the form
$\Pee(1, 2, \dots, 2)$, as we shall see in Part II.

In Cases (ii) and (iii) of Proposition 7.7, the $SO(2n+1)$ or $SO(2n)$
bundles have nonzero
$w_2$, and so do not lift to $Spin$ bundles. The statement for liftable
bundles is as follows:

\theorem{7.8} Let $V$ be a vector bundle of rank $2n$ over
$E$ with a nondegenerate symmetric form, and suppose that $V$ can be
lifted to a  principal
$Spin (2n)$-bundle.Then the dimension of the group of orthogonal
automorphisms of $V$ is always at least
$n$. If this dimension is exactly $n$, then
$V$ is isomorphic to 
$$\bigoplus _i\left(I_{d_i}(\lambda _i) \oplus I_{d_i}(\lambda
_i^{-1})\right)\oplus \bigoplus _{j\in
S}\left(I_{2a_j+1}(\eta_j)\oplus\eta _j\right),$$  where the $\lambda
_i$ are line bundles of degree zero, not of order two, such that, for
all
$i\neq j$, $\lambda _i\neq \lambda _j^{\pm1}$,  $\eta _0 =
\scrO_E, \eta_1,
\eta _2, \eta _3$ are the four distinct line bundles of order two on
$E$, and the second sum is over some subset $S$
\rom(possibly empty\rom) of $\{0,1,2,3\}$. Conversely, every such
vector bundle
$V$ has a nondegenerate symmetric form, all such forms have a group of
orthogonal automorphisms of dimension exactly $n$, and every two
nondegenerate symmetric forms on $V$ are equivalent  under the action
of $\Aut V$.

The case of $SO(2n+1)$ is similar, except that the summand
$I_{2a_0+1}\oplus\scrO_E$ is replaced by the odd rank summand 
$I_{2a_0+1}$, which must always be present.
\qed
\endstatement

Here the symmetric form on $I_{2a_0+1}\oplus\scrO_E$ consists of the
orthogonal direct sum of the  nondegenerate form on the factor
$I_{2a_0+1}$ given by (7.2)(ii), together with the obvious form on
$\scrO_E$, and similarly for the summands
$I_{2a_i+1}(\eta _i)\oplus \eta _i$. Moreover, not all of the summands
$I_{2a_j+1}(\eta _j)\oplus\eta _j$ need be present in
$V$.

\remark{Remark 7.9} We have shown that, if $\xi$ is a regular bundle
for $G = SL_n(\Cee)$ or $Sp(2n)$, then $\Aut _G\xi$ is abelian. If
$\xi$ is a liftable
$SO(2n)$-bundle and $\xi$ contains at least two factors of the form
$\eta_i\oplus \eta_i$ , then there is a copy of $O(2)$ inside the
$SO(2n)$-automorphisms of $\xi$ and thus this group is not abelian. In
general, however, the component group of the $SO(2n)$-automorphism
group acts nontrivially on the set of the four $Spin$-liftings of
$\xi$. If however $\xi$ contains a summand of the form 
$$(\eta_1\oplus \eta_1) \oplus (\eta_2\oplus \eta_2) \oplus
(\eta_3\oplus
\eta_3) \oplus (\scrO_E\oplus \scrO_E),$$ then inside the component
group, which is $\cong (\Zee/2\Zee)^4$, there is a subgroup
$(\Zee/2\Zee)^3$ where the product is trivial, corresponding to the
subgroup contained in $SO(2n)$. There is a subgroup $\cong \Zee/2\Zee$
of this $(\Zee/2\Zee)^3$ which fixes a given $Spin$-lifting and thus
acts on the lift. Hence the $Spin$-automorphism group is nonabelian in
this case.  
\endremark

\ssection{7.3. Regular conformal bundles.}

In this subsection, we consider models for the regular bundles
associated to the groups $Sp(2n)/\{\pm \Id\}$ and $SO(2n)/\{\pm \Id\}$.
For example, given the symplectic group $Sp(2n)$, we have the conformal
symplectic group $GSp(2n) =
\Cee^*\times _{\Zee/2\Zee}Sp(2n)$, with center $\Cee^*$, and the
quotient
$PSp(2n) = GSp(2n)/\Cee^* = Sp(2n)/\{\pm \Id\}$. The group $GSp(2n)$ is
the set of all $A\in GL_{2n}(\Cee)$ such that, if $Q$ is the standard
symplectic form, then for all $v, w\in \Cee^{2n}$, $Q(Av, Aw) =
\delta(A)Q(v,w)$ for some fixed
$\delta(A) \in \Cee^*$ (depending on $A$). The map $A\mapsto \delta(A)$
is a homomorphism $\delta \: GSp(2n) \to \Cee^*$ with kernel $Sp(2n)$,
and the restriction of
$\delta$ to the center is raising to the power $2$. Thus, the
restriction of the determinant to $GSp(2n)$ is equal to $\delta^n$. 

Let $\xi$ be a $PSp(2n)$-bundle. Since $H^2(E; \scrO_E^*) = 0$, $\xi$
lifts to a $GSp(2n)$-bundle, giving a vector bundle $V$ of rank $2n$, a
line bundle
$\lambda$ on $E$ corresponding to the character $\delta$, and a
nondegenerate skew-symmetric form
$Q\: V\otimes V
\to
\lambda$. Twisting $V$ by a line bundle $\mu$ on $E$ has the effect of
replacing $\lambda$ by $\lambda \otimes \mu^{\otimes 2}$. Thus
$\lambda$ is determined up to $2\Pic E$, and so we can assume that
$\lambda$ is either trivial or $\scrO_E(p_0)$. The bundle $\xi$ lifts
to $Sp(2n)$ if and only if
$\lambda$ is trivial, and so we will assume henceforth that $\lambda =
\scrO_E(p_0)$. In this case, the $GSp(2n)$-bundle $V$ lifting $\xi$ is
determined up to twisting by $2$-torsion in $\Pic E$. (In fact, for the
semistable bundles, one can see directly that the $2$-torsion line
bundles act trivially on the set of lifts.) It also follows that
$\det V = \scrO_E(np_0)$ and that $V\cong V\spcheck\otimes
\scrO_E(p_0)$.

If $V$ is a semistable vector bundle of rank $2n$ and determinant
$\scrO_E(np_0)$, then it follows from (ii) of Proposition 1.6 that $V$
is a direct sum of bundles of the form $I_d(W_2(1, \lambda))$, where
$\lambda$ is a line bundle on $E$ of degree $1$ and the product of the
$\lambda ^d$ is
$\scrO_E(np_0)$. For brevity we set $W_2(1, \lambda ) = W_2(q)$, where
$q$ is the unique point of $E$ such that $\lambda \cong \scrO_E(q)$.
Thus in particular $W_2 = W_2(p_0)$. Clearly $W_2(q)\spcheck \otimes
\scrO_E(p_0) \cong W_2(-q)$, where $-q$ is the inverse of $q$ in the
group law for $E$ where $p_0$ is the origin. Similarly
$(I_d(W_2(q)))\spcheck \otimes
\scrO_E(p_0)\cong I_d(W_2(-q))$. We then clearly have the following:

\lemma{7.10} Let $V = \bigoplus _iI_{d_i}(W_2(q_i))$ be a regular
vector bundle, and suppose that
$\Phi\: V\spcheck \otimes \scrO_E(p_0)\to V$ is an isomorphism. If
$q_i\neq -q_i$, then there exists a unique $j$ such $\Phi$ restricts to
an isomorphism from
$I_{d_i}(W_2(q_i))\spcheck \otimes
\scrO_E(p_0)$ to $I_{d_j}(W_2(q_j))$. Necessarily
$d_i = d_j$ and $q_i = -q_j$. \qed
\endstatement

If $q_i = -q_i$ in the above notation, then it is possible for $\Phi$
to define an isomorphism from the summand $I_{d_i}(W_2(q_i))\spcheck
\otimes
\scrO_E(p_0)$ to $I_{d_i}(W_2(q_i))$. The next lemma determines when
such an isomorphism can correspond to an alternating or symmetric form.
Note the parity change in the cases depending on whether $q=p_0$ or
$q\neq p_0$, as well as the contrast in parity with Theorem 7.2:

\lemma{7.11} Let $q\in E$.
\roster
\item"{(i)}" If $q\neq -q$, $H^0(E; Hom(I_d(W_2(q))\spcheck \otimes
\scrO_E(p_0), I_d(W_2(q)) ) = 0$. 
\item"{(ii)}" If $q=p_0$ and $d$ is odd, there exists a nondegenerate
alternating form $I_d(W_2(q)) \otimes I_d(W_2(q)) \to \scrO_E(p_0)$,
but no such nondegenerate symmetric form.
\item"{(iii)}" If $q=p_0$ and $d$ is even, there exists a nondegenerate
symmetric form $I_d(W_2(q)) \otimes I_d(W_2(q)) \to \scrO_E(p_0)$, but
no such nondegenerate alternating form.
\item"{(iv)}" If $q=-q,q\neq p_0$ and $d$ is odd, there exists a
nondegenerate symmetric form $I_d(W_2(q)) \otimes I_d(W_2(q)) \to
\scrO_E(p_0)$, but no such nondegenerate alternating form.
\item"{(v)}" If $q=-q,q\neq p_0$ and $d$ is even, there exists a
nondegenerate alternating form $I_d(W_2(q)) \otimes I_d(W_2(q)) \to
\scrO_E(p_0)$, but no such nondegenerate symmetric form. 
\endroster In all cases where a nondegenerate form exists, every two
such forms are conjugate under the action of $\Aut I_d(W_2(q))$. \qed
\endstatement
\proof First consider the case of $W_2(q)$ itself. If $q\neq -q$, then
there is no nondegenerate form
$W_2(q)\otimes W_2(q) \to \scrO_E(p_0)$. In case $q=p_0$, so that
$W_2(q) =W_2$, then there is a unique map $W_2\otimes W_2\to
\scrO_E(p_0)$ up to a nonzero scalar, the determinant, and it is
alternating. If $q=-q$ but $q\neq p_0$, then the isomorphism $W_2(q)
\cong W_2(q)\spcheck\otimes \scrO_E(p_0)$ defines a map
$Q\: W_2(q)\otimes W_2(q) \to \scrO_E(p_0)$, unique up to scalars since
$W_2(q)$ is simple. Hence $Q(v,w) = t Q(w,v)$ for some $t\in \Cee^*$,
necessarily $\pm 1$, so that
$Q$ is either alternating or symmetric. But if $Q$ were alternating, it
would induce a nonzero homomorphism $\bigwedge ^2W_2(q) = \scrO_E(q) \to
\scrO_E(p_0)$. This is impossible for $q\neq p_0$. Thus $Q$ is
symmetric. 

Now in general, every form $Q\: I_d(W_2(q)) \otimes I_d(W_2(q)) \to
\scrO_E(p_0)$ defines a homomorphism from $I_d(W_2(q))$ to
$I_d(W_2(q))\spcheck
\otimes \scrO_E(p_0)$, and conversely. Fixing a nondegenerate form
$R_0\: W_2(q) \otimes W_2(q) \to \scrO_E(p_0)$ defines an isomorphism
from 
$I_d(W_2(q))\spcheck \otimes \scrO_E(p_0)$ to $I_d\spcheck \otimes
W_2(q)$. By Lemma 1.15, the inclusion $\Hom (I_d, I_d\spcheck) \to \Hom
(I_d(W_2(q)),I_d(W_2(q))\spcheck
\otimes \scrO_E(p_0))$ is an isomorphism. Hence every form $Q$ is
induced by a corresponding form
$Q_0$ on $I_d$; more precisely, $Q = Q_0\otimes R_0$. Clearly $Q$ is
nondegenerate if and only if
$Q_0$ is nondegenerate. Moreover $Q$ is alternating if and only if
either $Q_0$ is alternating and
$R_0$ is symmetric or $Q_0$ is symmetric and $R_0$ is alternating.
Thus, a nondegenerate alternating form exists if and only if $q=p_0$
and $d$ is odd or  $q=-q, q\neq p_0$, and $d$ is even. The existence of
nondegenerate symmetric forms is similar. Finally, the last statement
follows from the corresponding statement in Theorem 7.2.
\endproof 

Using the above, we can give a description of regular bundles in the
conformally symplectic case. Before stating the general result, let us
work out the case of $GSp(4)$ as a concrete example. For $q\neq -q$,
there is a unique alternating form on $W_2(q)\oplus W_2(-q)$, unique up
to automorphisms of the bundle. Its group of conformally symplectic
automorphisms, modulo the image of
$\Cee^*$ acting by scalars, is $\Cee^*/\{\pm \Id\}$, where $t\in
\Cee^*$ acts via $t$ on the first factor and by $t^{-1}$ on the second.
For $q = p_0$, since
$I_2(W_2)$ does not have an alternating form, we are forced to use
$W_2\oplus W_2$ with its natural symplectic form. A calculation shows
that the conformally symplectic automorphisms of this bundle, modulo
$\Cee^*$, is given by $\Cee^*
\rtimes (\Zee/2\Zee)$, where the action of $\Zee/2\Zee$ is by $t\mapsto
t^{-1}$. In particular, the automorphism group is not abelian
(reflecting the fact that $O(2)$ is not abelian). For bundles
containing a factor of the form
$W_2(q)$, where $q=-q$ but
$q\neq p_0$, there is a symplectic form on $I_2(W_2(q))$, and its
automorphism group is abelian. (On the other hand, for $q=-q$ but
$q\neq p_0$, the conformal automorphism group of the bundle
$W_2(q) \oplus W_2(q)$ with the natural symplectic form is $GL_2(\Cee)$,
and so such bundles are not regular.)

More generally, we have the following: 

\theorem{7.12} Let $V$ be a vector bundle of rank $2n$, and suppose
that there is an alternating nondegenerate form $Q\: V\otimes V \to
\scrO_E(p_0)$.
\roster
\item"{(i)}" If $n$ is odd, then the minimal possible dimension for the
automorphisms of $V$ as a $PSp(2n)$-bundle is $(n-1)/2$. In this case, 
necessarily $V$ is isomorphic to
$$\bigoplus _i(I_{d_i}(W_2(q_i)) \oplus I_{d_i}(W_2(-q_i)))\oplus
\bigoplus _{j=1}^3I_{2e_j}(W_2(r_j)) \oplus I_{2a+1}(W_2),$$ where
$d_i$ are positive integers, $q_i\neq -q_i$, the $e_j$ and $a$ are
nonnegative integers, and the $r_j$ are the three points such that $r_j
= -r_j, r_j \neq p_0$. The automorphism group of the associated
$PSp(2n)$-bundle is always abelian.
\item"{(ii)}" If $n$ is even, then the minimal possible dimension for
the automorphisms of $V$ as a $PSp(2n)$-bundle is $n/2$. In this case, 
necessarily $V$ is isomorphic either to
$$\bigoplus _i(I_{d_i}(W_2(q_i)) \oplus I_{d_i}(W_2(-q_i)))\oplus
\bigoplus _{j=1}^3I_{2e_j}(W_2(r_j))$$ or to
$$\bigoplus _i(I_{d_i}(W_2(q_i)) \oplus I_{d_i}(W_2(-q_i)))\oplus
\bigoplus _{j=1}^3I_{2e_j}(W_2(r_j)) \oplus (I_{2a+1}(W_2)\oplus W_2),$$
where $d_i$ are positive integers, $q_i\neq -q_i$, the $e_j$ and $a$ are
nonnegative integers, and the $r_j$ are the three points such that $r_j
= -r_j, r_j \neq p_0$. The automorphism group of the associated
$PSp(2n)$-bundle is abelian if and only if the summand
$I_{2a+1}(W_2)\oplus W_2$ is not present.
\qed
\endroster
\endstatement
\proof  Let
$I$ be a vector bundle of the form
$\bigoplus _iI_{k_i}$. If $V$ has a summand of the form
$I\otimes W_2(q)$ with $q\neq -q$, the the argument is straightforward.
Consider summands of $V$ of the form
$I\otimes W_2(q)$ with $q = -q$. By Lemma 1.15,  the inclusion $\Hom
(I,I') \to \Hom (I\otimes W_2(q), I'\otimes W_2(q))$ is an isomorphism.
In particular,
$\Aut (I\otimes W_2(q)) \cong \Aut I$ under the natural map. Moreover,
viewing a form $Q$ on
$I\otimes W_2(q)$ with values in $\scrO_E(p_0)$ as a homomorphism from
$I\otimes W_2(q)$ to
$I\spcheck \otimes W_2(q)\spcheck \otimes \scrO_E(p_0)$ and using the
fixed form $R_0$ to identify $W_2(q)\spcheck \otimes \scrO_E(p_0)$ with
$W_2(q)$, we see that $Q$ determines a homomorphism from $I$ to
$I\spcheck$ and thus a form $Q_0$ on $I$. Hence $Q$ is of the form
$Q_0\otimes R_0$ for some form
$Q_0$ on $I$. It  follows that, if
$\Aut^{\Cee^*Q}$ denotes the group of conformal automorphisms of the
form $Q$, and similarly for $Q_0$, then
$\Aut^{\Cee^*Q}(I\otimes W_2(q)) \cong \Aut^{\Cee^*Q_0}(I)\cong
\Cee^*\times _{\Zee/2\Zee}\Aut^{Q_0}(I)$.  The theorem then follows
easily from Proposition 7.3 and Lemma 7.4.
\endproof

In Case (i) above, we see that the moduli space is given by $(n-1)/2$
points of
$E$ moduli sign change and permutation, which is a $\Pee^{(n-1)/2}$. In
Case (ii), the moduli space is likewise given by $n/2$ points of $E$
modulo sign change and permutation, and is a $\Pee^{n/2}$. 

Very similar results hold for the symmetric case. In this case, we
consider
$SO(2n)/\{\pm \Id\}$-bundles $\xi$ which do not lift to
$SO(2n)$-bundles. Note that if $n$ is odd, such bundles lift to no
quotient of $Spin (2n)$ by a proper subgroup of the center, whereas if
$n$ is even such bundles will always lift to the quotient of $Spin
(2n)$ by an exotic central subgroup of order $2$.

\theorem{7.13} Let $V$ be a vector bundle of rank $2n$, and suppose
that there is an symmetric nondegenerate form $Q\: V\otimes V \to
\scrO_E(p_0)$.
\roster
\item"{(i)}" If $n$ is odd, then the minimal possible dimension for the
automorphisms of $V$ as a $PSO(2n)$-bundle is $(n-3)/2$. In this case, 
necessarily $V$ is isomorphic to
$$\bigoplus _i(I_{d_i}(W_2(q_i)) \oplus I_{d_i}(W_2(-q_i)))\oplus
\bigoplus _{j=1}^3I_{2e_j+1}(W_2(r_j)) \oplus I_{2a}(W_2),$$ where
$d_i$ are positive integers, $q_i\neq -q_i$, the $e_j$ and $a$ are
nonnegative integers, and the $r_j$ are the three points such that $r_j
= -r_j, r_j \neq p_0$. The automorphism group of the associated
$PSO(2n)$-bundle is always abelian.
\item"{(ii)}" If $n$ is even, then the minimal possible dimension for
the automorphisms of $V$ as a $PSO(2n)$-bundle is $n/2$. In this case, 
necessarily $V$ is isomorphic to
$$\bigoplus _i(I_{d_i}(W_2(q_i)) \oplus I_{d_i}(W_2(-q_i)))\oplus
\bigoplus _{j\in S}(I_{2e_j}(W_2(r_j))\oplus W_2(r_j)) \oplus
I_{2a}(W_2),$$ where $d_i$ are positive integers, $q_i\neq -q_i$, the
$e_j$ and $a$ are nonnegative integers,  the $r_j$ are the three points
such that $r_j = -r_j, r_j \neq p_0$, and $S$ is a subset of
$\{1,2,3\}$, possibly empty.  \qed
\endroster
\endstatement

In both cases, the moduli space is a projective space. However, in Case
(ii) the bundles lift to $Spin (2n)/(\Zee/2\Zee)$, where the
$\Zee/2\Zee$ is the exotic subgroup of order $2$. As we shall see in
Part II, the corresponding moduli space for $Spin
(2n)/(\Zee/2\Zee)$-bundles is in fact a weighted projective space
$\Pee(1,1,1,2,2,\dots, 2)$. Likewise, the automorphism group as a $Spin
(2n)/(\Zee/2\Zee)$-bundle is not always abelian.


\Refs

\ref \no  1\by M. Atiyah \paper Vector bundles over an elliptic curve
\jour Proc. London Math. Soc. \vol 7\yr 1957 \pages 414--452\endref

\ref \no  2\by M. Atiyah and R. Bott \paper The Yang-Mills  equations
over Riemann surfaces \jour Phil. Trans. Roy. Soc. London A\vol 308\yr
1982\pages 523--615\endref




\ref \no  3\by A. Beauville and Y. Laszlo  \paper Conformal blocks and
generalized theta functions \jour Comm. Math. Phys.
\vol 164 \yr 1994 \pages 385--419 
\endref

\ref \no 4 \by I.N. Bernshtein and O.V. Shvartsman \paper Chevalley's
theorem for complex crystallographic Coxeter groups \jour Funct. Anal.
Appl. \vol 12 \yr 1978 \pages 308--310 \endref

\ref \no 5 \by A. Borel \paper Sous-groupes commutatifs et torsion des
groupes de Lie compactes \jour T\^ohoku Math. Jour. \vol 13 \yr 1961
\pages 216--240 \endref

\ref \no 6 \by A. Borel, R. Friedman, and J. W. Morgan \toappear
\endref


\ref \no  7\by N. Bourbaki \book Groupes et Alg\`ebres de  Lie 
\bookinfo Chap\. 4, 5, et 6 \publ Masson \publaddr Paris \yr 1981
\moreref \bookinfo Chap\. 7 et 8 \publ Hermann \publaddr Paris
\yr 1975 
\moreref \bookinfo Chap\. 9 \publ Masson \publaddr Paris \yr 1982\endref


\ref \no 8 \by S.K. Donaldson \paper A new proof of a theorem of
Narasimhan and Seshadri \jour J. Differential Geometry \vol 18 \yr 1983
\pages 269--277
\endref



\ref \no 9 \by G. Faltings \paper Stable $G$-bundles and projective
connections \jour J. Algebraic Geometry \vol 2\yr 1993 \pages 507--568
\endref

\ref \no 10 \bysame \paper A proof for the Verlinde formula \jour J.
Algebraic Geometry \vol 3\yr 1994 \pages 347--374
\endref

\ref \no  11\by R. Friedman \paper Rank two vector bundles over regular
elliptic surfaces \jour Inventiones Math. \vol 96 \yr 1989 \pages
283--332
\endref

\ref \no 12 \bysame \paper Vector bundles and
$SO(3)$-invariants for elliptic surfaces \jour J. Amer. Math. Soc. \vol
8 \yr 1995 
\pages 29--139 \endref

\ref \no 13\by R. Friedman and J.W. Morgan \book Smooth  Four-Manifolds
and Complex Surfaces \bookinfo Ergebnisse der Mathematik und ihrer
Grenzgebiete 3. Folge
\vol 27 \publ Springer-Verlag \publaddr Berlin Heidelberg New York
\yr 1994
\endref

\ref \no 14 \by R. Friedman, J.W. Morgan and E. Witten \paper  Vector
bundles and
$F$ theory \jour Commun. Math. Phys. \vol 187 \yr 1997 \pages 679--743
\endref

\ref \no 15 \bysame \paper  Principal $G$-bundles over elliptic curves
\jour Math. Research Letters \vol 5 \yr 1998 \pages 97--118
\endref

\ref \no 16 \bysame \paper  Vector bundles over elliptic fibrations
\jour J. Algebraic Geometry \toappear
\endref

\ref \no 17 \by A. Grothendieck \paper A general theory of fiber spaces
with structure sheaf \paperinfo University of Kansas Report Number
$4$\yr 1955
\endref

\ref \no 18 \bysame \paper Sur la classification des fibr\'es
holomorphes sur la sph\`ere de Riemann \jour Amer. J. Math \vol 79 \yr
1956 \pages 121--138 \endref



\ref \no 19 \by J. Humphreys  \book Linear Algebraic Groups
\bookinfo Graduate Texts in Mathematics \vol 21 \publ Springer-Verlag
\publaddr New York-Heidelberg-Berlin \yr 1975
\endref

\ref \no 20 \bysame  \book Conjugacy Classes in Semisimple Algebraic
Groups 
\bookinfo Mathematical Surveys and Monographs \vol 43\publ Amer. Math.
Soc. \publaddr Providence \yr 1995 \endref

\ref \no 21\by B. Kostant \paper The principal three-dimensional
subgroup and the Betti numbers of a complex simple Lie group
\jour Amer. J. Math. \vol 81 \yr 1959 \pages 973--1032
\endref

\ref \no 22\by Y. Laszlo \paper About $G$-bundles over elliptic curves
\jour Annales Inst. Fourier \vol 48 \yr 1998
\pages 413--424\endref

\ref \no 23\by J. Le Potier \book Lectures on Vector Bundles
\bookinfo Cambridge Studies in Advanced Mathematics \vol 54
\publ Cambridge University Press \publaddr Cambridge \yr 1997 \endref



\ref \no 24\by E. Looijenga  \paper Root systems and elliptic curves
\jour Invent. Math. \vol 38 \yr 1976
\pages 17--32
\endref



\ref \no  25\by M. S. Narasimhan and C. S. Seshadri \paper Stable and
unitary bundles on a compact Riemann surface \jour Annals of Math. \vol
82
\yr 1965 \pages 540--567 \endref

\ref \no 26\by S. Ramanan and A. Ramanathan\paper Some remarks on the
instability flag \jour T\^ohoku Math. Jour. 
\vol 36 \yr 1984 \pages 269--291
\endref

\ref \no 27\by A. Ramanathan\paper Stable principal bundles on a
compact Riemann surface \jour Math. Ann.
\vol 213 \yr 1975 \pages 129--152
\endref


\ref \no 28\bysame\paper Moduli for principal bundles over algebraic
curves I,II
\jour Proc. Indian Acad. Sci. Math. \vol 106 \yr 1996 \pages 301--328,
421--449
\endref

\ref \no 29\by C. Schweigert \paper On moduli spaces of flat
connections with non-simply connected structure group \jour Nucl. Phys.
B \vol 492 \yr 1997 \pages 743--755 \endref


\endRefs
\enddocument